\documentclass[eqthmnum,nocolour,skinny]{jt-calcs}

\usepackage[euler-digits,euler-hat-accent]{eulervm}
\usepackage{mathrsfs}
\usepackage{verbatim}
\usepackage{tikz-cd}

\usepackage[backend=bibtex,style=alphabetic,sorting=nyt,maxalphanames=6,maxbibnames=99]{biblatex}
\bibliography{galois-autos}

\usepackage{kbordermatrix}
\usepackage{booktabs}
\usepackage{array}

\usepackage{mathdots}

\newcommand{\der}{\mathrm{der}}
\newcommand{\mon}{\mathrm{mon}}
\newcommand{\grp}{\mathrm{grp}}
\newcommand{\spr}{\mathrm{spr}}
\newcommand{\Clu}{\mathrm{Cl_{uni}}}
\newcommand{\rE}{\mathrm{E}}
\newcommand{\rD}{\mathrm{D}}
\newcommand{\rN}{\mathrm{N}}
\newcommand{\rO}{\mathrm{O}}
\newcommand{\rK}{\mathrm{K}}
\newcommand{\rC}{\mathrm{C}}
\newcommand{\rZ}{\mathrm{Z}}

\DeclareMathOperator{\Spr}{Spr}

\DeclareMathOperator{\Aso}{Aso}
\DeclareMathOperator{\Tor}{Tor}

\DeclareMathOperator{\Cusp}{Cusp}
\DeclareMathOperator{\Frob}{Frob}

\setlist[enumerate,1]{label={\normalfont (\roman*)}}

\newcommand{\ssymb}[2]{\ensuremath{\bigl[\begin{smallmatrix} #1 \\ #2 \end{smallmatrix}\bigr]}}

\newcommand{\symb}[2]{\ensuremath{\begin{bmatrix} #1 \\ #2 \end{bmatrix}}}

\newcommand{\wt}[1]{\widetilde{#1}}
\newcommand{\wc}[1]{\widecheck{#1}}
\newcommand{\wh}[1]{\widehat{#1}}
\newcommand{\bg}[1]{\mathbf{#1}}
\newcommand{\gal}{\mathcal{G}}
\newcommand{\galh}{\mathcal{H}}
\newcommand{\Q}{\mathbb{Q}}
\newcommand{\la}{\lambda}
\newcommand{\legendre}[2]{\genfrac{(}{)}{}{}{#1}{#2}}

\crefname{enumerate}{}{}

\title{Galois Automorphisms and Classical Groups}
\author{A.~A. Schaeffer Fry and Jay Taylor}
\address{J. Taylor, Department of Mathematics, The University of Manchester, Oxford Road, Manchester, M13 9PL, UK}
\email{jay.taylor@manchester.ac.uk}
\address{A.~A. Schaeffer Fry, Department of Mathematical and Computer Sciences, Metropolitan State University of Denver, Denver, CO 80217, USA}
\email{aschaef6@msudenver.edu}
\mscno{2010}{20C30}{20C15}
\keywords{Classical groups, character fields, Galois automorphisms, Harish-Chandra theory, generalised Gelfand--Graev characters, unipotent elements}

\begin{document}
\begin{abstract}
In a previous work, the second-named author gave a complete description of the action of automorphisms on the ordinary irreducible characters of the finite symplectic groups. We generalise this in two directions. Firstly, using work of the first-named author, we give a complete description of the action of Galois automorphisms on irreducible characters. Secondly, we extend both descriptions to cover the case of special orthogonal groups. As a consequence, we obtain explicit descriptions for the character fields of symplectic and special orthogonal groups.
\end{abstract}

\section{Introduction}

\begin{pa}
Given a finite group $G$, one has natural actions of its automorphism group $\Aut(G)$ and the absolute Galois group $\Gal(\overline{\mathbb{Q}}/\mathbb{Q})$ on the set of its ordinary irreducible characters $\Irr(G)$. Many questions in character theory concern these actions. For instance, determining the \emph{character field} $\mathbb{Q}(\chi) = \mathbb{Q}(\chi(g) \mid g \in G)$ of $\chi \in \Irr(G)$ is equivalent to understanding the stabiliser of $\chi$ in $\Gal(\overline{\mathbb{Q}}/\mathbb{Q})$. It is the purpose of this article to study, in some detail, these actions in the case of finite symplectic and special orthogonal groups. Our results sharpen, in these cases, the general statements on character fields obtained by Geck \cite{geck:2003:character-values-schur-indices-and-char-sheaves} and Tiep--Zalesski \cite{tiep-zalesski:2004:unipotent-elements} and extend previous results in type $\A$ \cite{turull01,SFV19}.

The need for such precise information regarding these actions has become increasingly more relevant owing to recent developments regarding the McKay--Navarro Conjecture, sometimes referred to as the Galois--McKay Conjecture. This conjecture has just been reduced to a problem about quasi-simple groups by Navarro--Sp\"ath--Vallejo \cite{nsv}. The resulting problem involves understanding explicitly the actions of automorphisms and Galois automorphisms on irreducible characters. 

Checking these conditions seems to be extremely complicated, with it already being a challenge in the case of $\SL_2(q)$. In the verification of the McKay conjecture for the prime $2$ (recently completed by Malle--Sp{\"a}th \cite{malle-spaeth:2016:characters-of-odd-degree}), special consideration was needed for the case of symplectic groups, see \cite{malle:2008:extensions-of-unipotent-characters}.  It stands to reason that the same will be true for the McKay--Navarro conjecture.  In fact, since the initial submission of the current paper, the results here have been used in \cite{SF21, RSF21} to make significant progress on proving the required conditions from \cite{nsv} for quasi-simple groups and the prime $2$.

To give our first results, let us fix a pair $(\bG^{\star},F^{\star})$ dual to $(\bG,F)$, where $\bG$ is a connected reductive algebraic group and $F$ is a Frobenius morphism on $\bG$ yielding an $\mathbb{F}_q$-rational structure. Given any semisimple element $s \in \bG^{\star F^{\star}}$, we have a corresponding (rational) Lusztig series $\mathcal{E}(\bG^F,s) \subseteq \Irr(\bG^F)$. We use $\Inn(\bG,F)$ to denote the group of inner diagonal automorphisms, defined precisely in \ref{pa:aso} below, and $\Inn(\bG,F)_{\chi}$ for the stabiliser of $\chi\in\Irr(\bG^F)$ under the natural action.
\end{pa}

\begin{mthm}\label{mainthmgalSO}
Assume $\bG = \SO(V)$ is a special orthogonal group defined over $\mathbb{F}_q$ with $q$ odd. Let $F_p,\gamma : \bG \to \bG$ be field and graph automorphisms, respectively, with $\gamma$ the identity when $\dim(V)$ is odd. Assume $s \in \bG^{\star F^{\star}}$ is a quasi-isolated semisimple element, i.e., $s^2 = 1$, and that $s$ is not contained in the centre of $\bG^{\star F^{\star}}$. Then for any $\chi \in \mathcal{E}(\bG^F,s)$, the following hold:
%%%%
\begin{enumerate}
	\item $\mathbb{Q}(\chi) = \mathbb{Q}$,
	\item $\chi^{F_p} = \chi$,
	\item $\chi^{\gamma} = \chi$ if and only if $\Inn(\bG,F)_{\chi} = \Inn(\bG^F)$.
\end{enumerate}
\end{mthm}

\begin{pa}
The case when $s$ is central is very well known. In that case (i) and (ii) also hold but (iii) does not. When $\dim(V)$ is odd or $q$ is even, the group $\SO(V)$ has a trivial, hence connected, centre. In these cases, the statements in Theorem~\ref{mainthmgalSO} are all easy consequences of the unicity of multiplicites in Deligne--Lusztig virtual characters, see \cite[Prop.~6.3]{digne-michel:1990:lusztigs-parametrization} or \cite[Thm.~4.4.23]{geck-malle:2020:the-character-theory-of-finite-groups-of-lie-type}. However, we include the odd-dimensional case in our analysis as it is little extra effort and some results used on the way to prove Theorem~\ref{mainthmgalSO} may be of independent interest.

Our next main result gives the analogue of (i) of Theorem~\ref{mainthmgalSO} in the case of symplectic groups. The character fields themselves are a little more complicated in this case, but the behaviour is still uniform across a Lusztig series.
\end{pa}

\begin{mthm}\label{mainthmgalSp}
Assume $\bG = \Sp(V)$ is a symplectic group defined over $\mathbb{F}_q$ with $q$ odd. Assume $s \in \bG^{\star}$ is quasi-isolated and $\chi\in \mathcal{E}(\bG^F,s)$. Then $\mathbb{Q}(\chi) \subseteq \mathbb{Q}(\sqrt{\omega p})$,  where $\omega = (-1)^{\frac{p-1}{2}}$ is the unique element of $\{\pm1\}$ for which $p \equiv \omega \pmod{4}$. Moreover, the following hold:
%%%%
\begin{enumerate}
	\item if $q$ is a square or $s=1$, then $\mathbb{Q}(\chi) = \mathbb{Q}$;
	\item if $q$ is not a square and $s\neq 1$, then $\mathbb{Q}(\chi) = \mathbb{Q}(\sqrt{\omega p})$. In particular, if $\sigma \in \mathrm{Gal}(\overline{\Q}/\Q)$, then $\chi^{\sigma} = \chi$ if and only if $\sqrt{\omega p}^{\sigma} = \sqrt{\omega p}$.
\end{enumerate}
%%%%
\end{mthm}

\begin{pa}
On the way to Theorem~\ref{mainthmgalSp}, we give further details for the action of the Galois group on the characters in this case. Moreover, in Sections~\ref{sec:actionH} and \ref{sec:actionH-Sp}, we describe the action of the specific Galois automorphisms relevant for the McKay--Navarro conjecture, which will be useful in studying the conjecture for finite reductive groups. Note that, whilst seemingly only covering special cases, Theorems~\ref{mainthmgalSO} and \ref{mainthmgalSp} do allow one to completely determine the field of values for all irreducible characters of $\bG^F$ when $\bG=\SO(V)$ or $\Sp(V)$. Specifically, we obtain the following:
\end{pa}

\begin{mthm}\label{thm:completecharfields}
Assume $\bG = \Sp(V)$ or $\SO(V)$ is defined over $\mathbb{F}_q$ with $q$ odd and let $s \in \bG^{\star F^{\star}}$ be a semisimple element. Assume $s$ has order $d > 0$ and let $\zeta_d\in\overline{\Q}^\times$ be a primitive $d$th root of unity. We let $\mathbb{Q}_s = \overline{\mathbb{Q}}^{\gal_s} \subseteq \mathbb{Q}(\zeta_d)$ be the subfield fixed by the subgroup $\mathcal{G}_s \leqslant \mathcal{G}$ stabilising the series $\mathcal{E}(\bG^F,s)$, which can be described in terms of conjugacy in $\bG^{\star F^{\star}}$, as in Remark~\ref{rem:gal-stab-L-series}.
%%%%
\begin{enumerate}
	\item If $\bG = \SO(V)$, then any $\chi \in \mathcal{E}(\bG^F,s)$ has character field $\mathbb{Q}(\chi) = \mathbb{Q}_s$.
	\item If $\bG = \Sp(V)$, then any $\chi \in \mathcal{E}(\bG^F,s)$ has character field $\Q(\chi)$, where:
	\begin{itemize}
	\item $\mathbb{Q}(\chi) = \mathbb{Q}_s$ if $q$ is a square or $-1$ is not an eigenvalue of $s$;
	\item $\Q_s \subsetneq \mathbb{Q}(\chi) = \mathbb{Q}_s(\sqrt{\omega p})$ otherwise, where $\omega=(-1)^{(p-1)/2}$.
	\end{itemize} 
\end{enumerate}
\end{mthm}

\begin{pa}
We note that, assuming $q$ is sufficiently large, it is theoretically possible to compute the character values of the characters of symplectic groups using a result of Waldspurger \cite{waldspurger:2004:lusztigs-conjecture}. However, no such statement is known for the special orthogonal groups, as the results in \cite{waldspurger:2004:lusztigs-conjecture} concern the (disconnected) orthogonal group. Our arguments, which are based on Kawanaka's theory of Generalised Gelfand--Graev Characters (GGGCs) and the Howlett--Lehrer parametrisation of Harish-Chandra series, require no restriction on $q$ and do not require the precise determination of character values.

We now wish to describe two consequences of our description of the action of $\Aut(\bG^F)$ on $\Irr(\bG^F)$ when $\bG = \SO(V)$. For this, let us recall that the automorphism group $\Aut(\bG^F)$ has a subgroup $\Aso(\bG,F) \leqslant \Aut(\bG^F)$ generated by inner, diagonal, field, and graph automorphisms. This group is defined precisely in \ref{pa:aso}. There is a complement $\Gamma(\bG,F) \leqslant \Aso(\bG,F)$ to the subgroup $\Inn(\bG,F) \leqslant \Aso(\bG,F)$ such that $\Gamma(\bG, F)$ consists of field and graph automorphisms. We choose this explicitly in Lemma~\ref{lem:asogeny-orthogonal-fin}. Recalling the notion of a Jordan decomposition from \cite[\S11.5]{digne-michel:2020:representations-of-finite-groups-of-lie-type}, we will show the following.
\end{pa}

\begin{mthm}\label{thm:equiv-J-decomp-no-D}
Assume $\bG = \SO(V)$ is a special orthogonal group defined over $\mathbb{F}_q$ with $q$ odd and $\dim(V)$ even. Moreover, let $\gamma : \bG \to \bG$ and $\gamma^{\star} : \bG^\star \to \bG^\star$ be dual graph automorphisms.
%%%%
\begin{enumerate}
	\item If $s \in \bG^{\star F^{\star}}$ is $\gamma^{\star}$-fixed, then there exists a Jordan decomposition $J_s^{\bG} : \mathcal{E}(\bG^F,s) \to \mathcal{E}(\rC_{\bG^{\star}}(s)^{F^{\star}},1)$ satisfying $J_s^{\bG}(\chi^{\gamma})^{\gamma^{\star}} = J_s^{\bG}(\chi)$ for all $\chi \in \mathcal{E}(\bG^F,s)$.
	\item Each $\chi \in \Irr(\bG^F)$ satisfies the condition $\Aso(\bG,F)_{\chi} = \Inn(\bG,F)_{\chi} \rtimes \Gamma(\bG,F)_{\chi}$.
\end{enumerate}
\end{mthm}

\begin{pa}
The statement in (i) of Theorem~\ref{thm:equiv-J-decomp-no-D} was claimed in the proof of a result of Aubert--Michel--Rouquier \cite[Prop.~1.7]{aubert-michel-rouquier:1996:correspondence-de-howie}. Their result states that the orthogonal group $\GO(V)$ also has a Jordan decomposition when $\dim(V)$ is even. The existence of such a Jordan decomposition has recently been used by Vinroot in his work on totally orthogonal groups \cite{vinroot:2020:totally-orthogonal}, who noticed the arguments in \cite{aubert-michel-rouquier:1996:correspondence-de-howie} are insufficient. We thank him for bringing this issue to our attention.

The statement in (ii) of Theorem~\ref{thm:equiv-J-decomp-no-D} was shown to hold in the case $\bG = \Sp(V)$ in \cite{taylor:2018:action-of-automorphisms-symplectic} and independently, by a completely different argument, in \cite{cabanes-spaeth:2017:inductive-mckay-type-C}. Such a statement is of interest owing to Sp\"ath's criterion for the Inductive McKay Conditions \cite{spaeth:2012:inductive-mckay-defining}. We also remark that, using the known actions of $\Aut(\bG^F)$ on irreducible characters, Li \cite[Thm.~5.9]{li:2019:an-equivariant-bijection} has shown the existence of an equivariant Jordan decomposition when $\bG = \Sp(V)$. Similar arguments can be applied to our results here to establish such a statement when $\bG = \SO(V)$. However, we have decided not to pursue this.
\end{pa}

\begin{pa}
Given Theorem~\ref{thm:completecharfields}, one might expect that when $\bG = \SO(V)$, there is a Jordan decomposition that is equivariant with respect to Galois automorphisms, as in the case of connected centre established in \cite{SrinivasanVinroot-1}. Of course, such a statement cannot hold in general, as it will not hold when $\bG = \Sp(V)$. The case of $\SL_2(q) \cong \Sp_2(q)$ already provides a counterexample, as pointed out in \cite[\S5]{SrinivasanVinroot-1}. Given our work here, the only issue remaining to prove such a statement for $\SO(V)$ is to determine the rationality of unipotent characters of the (disconnected) orthogonal group $\GO(V)$.

With regards to this, we obtain an important partial result, see Theorem~\ref{prop:cuspidalisol}. Namely, we prove that cuspidal unipotent characters of $\GO(V)$ are rational-valued. Extending the results of \cite{schaeffer-fry:2019:galois-automorphisms} to disconnected groups would easily yield the complete description. We intend to consider this in future work.
\end{pa}

\begin{pa}\label{pa:intro-weird-aut}
We end with a comment on the automorphism group of $\bG^F$ when $\bG = \SO(V)$ and $\dim(V)$ is even. Under certain conditions on $q$, we have $\rZ(\bG^F) \leqslant \rO^{p'}(\bG^F)$. When this happens, the group $\bG^F$ admits an automorphism that is not contained in the group $\Aso(\bG,F)$, see Proposition~\ref{prop:aut-grp-quasisimp}. The action of this automorphism is closely related to the study of central characters. To give a complete description for the action of $\Aut(\bG^F)$ on $\Irr(\bG^F)$, we must also describe the action of this automorphism. This is dealt with in Section~\ref{sec:notaso}.
\end{pa}

\subsection*{Structure of the Paper}
\begin{pa}
In Sections~\ref{sec:auts} to \ref{sec:actionH} we give some generalities on automorphisms and Galois automorphisms, in addition to recalling a few results from the literature to be used in further sections. In this paper, we use Kawanaka's GGGCs to translate the automorphism and Galois actions to questions on unipotent conjugacy classes. To understand the Galois actions, we study the permutation on conjugacy classes induced by power maps in Sections~\ref{sec:power-maps-1} and \ref{sec:power-maps-2}. This allows us to conclude, for instance, that all GGGCs of special orthogonal groups are rational-valued.

To understand the actions of $\Aut(\bG^F)$ and $\Gal(\overline{\mathbb{Q}}/\mathbb{Q})$ on irreducible characters of $\SO(V)$, it is necessary to understand such actions on extensions of these characters to $\GO(V)$. In Section~\ref{sec:GGGC-ext}, we define extensions of GGGCs for disconnected groups and study the relevant automorphism and Galois actions on these characters. To show (i) of Theorem~\ref{thm:equiv-J-decomp-no-D}, we need to describe the action of automorphisms on unipotent characters of certain disconnected groups, which we do in Section~\ref{sec:aut-uni-char}.

We study the case of quasi-isolated Lusztig series in Sections~\ref{sec:cusp-chars} to \ref{sec:gal-qiso-sp}. This is broken up by studying individual Harish-Chandra series. We first treat the case of cuspidal characters of special orthogonal groups in Section~\ref{sec:cusp-chars}, which extends results in \cite{taylor:2018:action-of-automorphisms-symplectic}. For another approach to the action of $\Aut(\bG^F)$ on cuspidal characters, see \cite{malle:2017:cuspidal-characters-and-automorphisms}. Hecke algebra techniques and the results of \cite{schaeffer-fry:2019:galois-automorphisms} are then used in Sections~\ref{sec:galoisSO} and \ref{sec:gal-qiso-sp} to study the members of a Harish-Chandra series, first in the special orthogonal group case and then in the symplectic group case. The extra automorphism of $\SO(V)$ mentioned in \ref{pa:intro-weird-aut} is studied in Section~\ref{sec:notaso}. Finally, our reductions to the quasi-isolated case, as well as the proof of Theorems~\ref{thm:completecharfields} and \ref{thm:equiv-J-decomp-no-D}, are contained in Section~\ref{sec:red-qi-case}.
\end{pa}

\subsection*{Notation}

\begin{pa}
For a set $X$ and elements $x,y \in X$, we denote by $\delta_{x,y} \in \{0,1\}$ the usual Kronecker delta. Assume $G$ is a group. If $G$ acts on a set $X$, then we denote by $G_x \leqslant G$ the stabiliser of $x \in X$. We denote by $\Aut(G)$ the automorphism group of $G$ as an abstract group, ignoring any additional structure on $G$. For any $g \in G$, we will use $\iota_g : G \to G$ to denote the inner automorphism defined by $\iota_g(x) = {}^gx = x^{g^{-1}} = gxg^{-1}$ and $\Inn(G) \leqslant \Aut(G)$ to denote the inner automorphism group.

{If $S \subseteq G$ is a subset of a group $G$ then we denote by $\langle S\rangle_{\mon}$ or $\langle S\rangle_{\grp}$ the submonoid or subgroup generated by $S$. We typically denote $\langle S\rangle_{\grp}$ by $\langle S\rangle$ if no extra clarity is needed. Of course, if $\langle S\rangle_{\mon}$ is finite then $\langle S\rangle_{\mon} = \langle S\rangle_{\grp}$.}

Assume $S$ is an arbitrary set. If $\phi : G \to H$ is a group isomorphism, then for any function $f : G \to S$, respectively $f : H\to S$, we let ${}^{\phi}f = f\circ \phi^{-1} : H \to S$, respectively $f^{\phi} = f\circ\phi : G \to S$. Taking $H=G$ this defines left and right actions of $\Aut(G)$ on functions $G \to S$. If $g \in G$, then we also let ${}^gf = {}^{\iota_g}f$ and $f^g = f^{\iota_g}$. Throughout, we will usually use right actions, although we will use the left action when it is notationally convenient.

We fix a prime $p$ and denote by $\mathbb{F} = \overline{\mathbb{F}}_p$ an algebraic closure of the finite field $\mathbb{F}_p$. Any algebraic group $\bG$, usually denoted in bold, is assumed to be affine and defined over $\mathbb{F}$. We denote by $\bG^{\circ}$ the connected component containing the identity. For any $x \in \bG$, we denote by $A_{\bG}(x)$ the component group $C_{\bG}(x)/C_{\bG}^{\circ}(x)$ of the centraliser. Note that if $x \in \bG^{\circ}$, then $A_{\bG^{\circ}}(x) \leqslant A_{\bG}(x)$ is a subgroup, since $C_{\bG}^{\circ}(x) = C_{\bG^{\circ}}^{\circ}(x)$.

We fix a prime $\ell \neq p$ and let $\Ql$ be an algebraic closure of the $\ell$-adic numbers. We identify $\overline{\mathbb{Q}}$ with a subfield of $\Ql$. If $G$ is a finite group, then $\Irr(G)$ is the set of $\Ql$-irreducible characters of $G$. We denote by $1_G \in \Irr(G)$ the trivial character. If $G$ is cyclic, we denote by $\varepsilon_G \in \Irr(G)$ the unique character with kernel $\{x^2 \mid x \in G\}$. For each $\chi \in \Irr(G)$, we have the central character $\omega_{\chi} : \rZ(G) \to \Ql$ defined by $\omega_{\chi}(z) = \chi(z)/\chi(1)$. For a subgroup $H\leq G$ and characters $\chi\in\Irr(G)$ and $\varphi\in\Irr(H)$, we will write $\Res_H^G(\chi)$ and $\Ind_H^G(\varphi)$ for the corresponding restricted and induced characters, respectively. 

Throughout, we let $\gal:=\mathrm{Gal}(\overline{\Q}/\Q)$ denote the absolute Galois group of $\Q$. Note that if $G$ is a finite group and $\chi \in \Irr(G)$, then $\chi(g) \in \overline{\mathbb{Q}}$. For any $\xi \in \overline{\Q}$ and $\sigma \in \gal$, we let $\xi^{\sigma}$ be the image of $\xi$ under $\sigma$. We have an action of $\gal$ on $\Irr(G)$ given by $\chi^{\sigma}(g) = \chi(g)^{\sigma}$. For each character $\chi$ of $G$, we fix a $\Ql G$-module $M_{\chi}$ affording $\chi$. We then write $\End_G(\chi)$ to denote the endomorphism algebra $\End_{\Ql G}(M_{\chi})$.
\end{pa}

\section{Automorphisms}\label{sec:auts}
\begin{pa}\label{pa:central-auts}
For a group $G$, the group $\rC_{\Aut(G)}(\Inn(G)) = \{\varphi \in \Aut(G) \mid g^{-1}\varphi(g) \in \rZ(G)$ for all $g \in G\}$ is a normal subgroup of $\Aut(G)$ whose elements are referred to as \emph{central automorphisms} in the literature. If $\alpha \in \Hom(G,\rZ(G))$ is a homomorphism, then we define $\tau_{\alpha} : G \to G$ by setting $\tau_{\alpha}(g) = g\cdot\alpha(g)$. This is also a homomorphism. Note that if $\varphi \in \rC_{\Aut(G)}(\Inn(G))$ is a central automorphism then $\varphi = \tau_{\alpha}$ with $\alpha \in \Hom(G,\rZ(G))$ defined by $\alpha(g) = g^{-1}\varphi(g)$.

Clearly there will be homomorphisms $\alpha \in \Hom(G,\rZ(G))$ for which $\tau_{\alpha}$ is not an automorphism. Simply take $G = \rZ(G)$ and $\alpha$ the inversion map. However, if $\tau_{\alpha}$ is an automorphism then it is clearly central. In other words, the central automorphisms are precisely those homomorphisms $\tau_{\alpha}$ that are automorphisms. The following easy observation gives a way to check when $\tau_{\alpha}$ is an automorphism.
\end{pa}

\begin{lem}\label{lem:autos-triv-on-derived}
Assume $G$ is a finite group.
%%%%
\begin{enumerate}
	\item Given a homomorphism $\alpha \in \Hom(G,\rZ(G))$, we have $\tau_{\alpha}$ is an automorphism if and only if $\tau_{\alpha}|_{\rZ(G)} = \mathrm{Id}_{\rZ(G)}+(\alpha|_{\rZ(G)})$ is an automorphism of $\rZ(G)$. Here we consider $\End(\rZ(G))$ as an additive abelian group so that for any $\phi,\psi \in \End(G)$ we have $(\phi+\psi)(z) = \phi(z)\psi(z)$ for all $z \in \rZ(G)$.

	\item If $\chi\in\Irr(G)$ and $\tau_{\alpha} \in \rC_{\Aut(G)}(\Inn(G))$, then $\chi^{\tau_\alpha} = (\omega_{\chi}\circ \alpha)\chi$, where $\omega_{\chi}\circ \alpha \in \Irr(G)$ is a linear character, and $\omega_{\chi^{\tau_{\alpha}}} = \omega_{\chi}^{\tau_{\alpha}}$.
\end{enumerate}
\end{lem}

\begin{proof}
(i). It suffices to know when $\tau_{\alpha}$ is injective. If $g \in \Ker(\tau_{\alpha})$ then $\alpha(g) = g^{-1} \in \rZ(G)$ so $\Ker(\tau_{\alpha}) \subseteq \rZ(G)$. The statement follows easily from this.

(ii). This is a simple consequence of the fact that $\chi(zg) = \omega_{\chi}(z)\chi(g)$ for any $z \in \rZ(G)$ and $g \in G$.
\end{proof}

\begin{rem}\label{rem:cent-in-der}
We will be interested in the case where $\rZ(G)$ is of order $2$. Let $\rE^2(G) \lhd G$ be the normal subgroup minimal with respect to the condition that the quotient $G/\rE^2(G) \cong C_2\times\cdots\times C_2$ is an elementary abelian $2$-group.  We then have
%%%%
\begin{equation*}
\Hom(G,\rZ(G)) \cong \Hom(G/\rE^2(G),\rZ(G)) \cong G/\rE^2(G)
\end{equation*}
%%%%
has order $2^s$ for some integer $s \geqslant 0$. Indeed, viewing $G/\rE^2(G)$ as an $\mathbb{F}_2$-vector space we have this is the dual space.

In this case $\tau_{\alpha}$ is an automorphism of $G$ if and only if $\rZ(G) \leqslant \Ker(\alpha)$. Hence, inflating $\Hom(G/\rZ(G)\rE^2(G),\rZ(G))$ gives those $\alpha \in \Hom(G,\rZ(G))$ for which $\tau_{\alpha}$ is an automorphism. Thus, there are $2^s$ such automorphisms if $\rZ(G) \leqslant \rE^2(G)$ and $2^{s-1}$ otherwise.

Note that this discussion, together with (ii) of Lemma~\ref{lem:autos-triv-on-derived}, shows that any $\tau_{\alpha} \in \Aut(G)$ fixes the central character of any $\chi \in \Irr(G)$ when $|\rZ(G)|=2$.
\end{rem}

\begin{pa}\label{pa:aso}\label{pa:alg-auts}
Assume $\bG$ is an affine algebraic group over $\mathbb{F}$ (not necessarily connected). A bijective homomorphism of algebraic groups $\varphi : \bG \to \bG$ will be called an \emph{asogeny}, and we denote by $\Aso(\bG) \leqslant \Aut(\bG)$ the \emph{submonoid} consisting of all asogenies.

Of interest to us will be the following subgroups of $\Aut(\bG)$. Firstly, we denote by $\Aut_0(\bG) \leqslant \Aut(\bG)$ the subgroup consisting of those automorphisms $\varphi \in \Aut(\bG)$ for which both $\varphi$ and $\varphi^{-1}$ are asogenies. We then also define $\Aut_1(\bG) \leqslant \Aut(\bG)$ to be the subgroup generated by the asogenies $\Aso(\bG)$.

For any Frobenius root $F : \bG \to \bG$, in the sense of \cite[Def.~4.2.4]{digne-michel:2020:representations-of-finite-groups-of-lie-type}, we have $F \in \Aso(\bG)$. Restricting to $\bG^F$, we obtain a natural monoid homomorphism $\rC_{\Aso(\bG)}(F) \to \Aut(\bG^F)$. The subgroup of $\Aut(\bG^F)$ generated by the image of this map will be denoted by $\Aso(\bG,F)$. We will denote by $\Inn(\bG,F) \leqslant \Aso(\bG,F)$ the image of $\rC_{\Inn(\bG)}(F)$, which is a subgroup containing $\Inn(\bG^F)$. We define $\Out(\bG,F) = \Aso(\bG,F)/\Inn(\bG^F)$.

Recall that if $G$ is a finite group, then $\rO^{p'}(G) \lhd G$ is the unique normal subgroup minimal with respect to the condition that the quotient $G/\rO^{p'}(G)$ is a $p'$-group.  Equivalently, this is the subgroup generated by all $p$-elements of $G$. Using the notation of Lemma~\ref{lem:autos-triv-on-derived}, we define a subgroup
%%%%
\begin{align*}
\rK(G) &= \{\varphi \in \rC_{\Aut(G)}(\Inn(G)) \mid \varphi|_{\rO^{p'}(G)} = \mathrm{Id}_{\rO^{p'}(G)}\} \leqslant \Aut(G).
\end{align*}
%%%%
This should be denoted $K_p(G)$ but as $p$ is implicit we omit this from the notation.
\end{pa}

\begin{rem}
Suppose we are in the setting of Remark~\ref{rem:cent-in-der}. Then $K(G)$ consists of those $\tau_{\alpha}$ where $\alpha$ is obtained by inflation from $\Hom(G/\rZ(G)\rO^{p'}(G)\rE^2(G),\rZ(G))$.
\end{rem}

The following is well known to finite group theorists, but we include some details in the proof as it is difficult to extract this exact statement from the literature.

\begin{prop}\label{prop:aut-grp-quasisimp}
Assume $\bG$ is a quasisimple algebraic group with $F : \bG \to \bG$ a Frobenius root. Then the map $\rC_{\Aso(\bG)}(F) \to \Aso(\bG,F)$ is surjective and $\Aut(\bG^F) =  \rK(\bG^F) \rtimes \Aso(\bG,F)$.
\end{prop}

\begin{proof}
The subgroup $\Aut_1(\bG) \leqslant \Aut(\bG)$ is described in \cite[Thm.~1.15.7]{gorenstein-lyons-solomon:1998:classification-3}. Picking a set of Steinberg generators for $\bG$, we have a splitting $\Aut_1(\bG) = \Inn(\bG) \rtimes \Gamma(\bG)$, where $\Gamma(\bG)$ is generated by field and graph automorphisms. In fact $\Gamma(\bG)$ is generated, as a group, by graph automorphisms whose inverse is also an asogeny and an element $\psi$, which is either a standard Frobenius $F_p$ or $\psi^2 = F_p$ is such. In any case, one sees that the submonoid of $\Aut(\bG)$ generated by $\Aso(\bG)$ and $F_p^{-1} \in \Aut(\bG)$ is the subgroup $\Aut_1(\bG)$.

Up to conjugacy, we can assume $F \in \Gamma(\bG)$, so that $C_{\Aut_1(\bG)}(F) = C_{\Inn(\bG)}(F) \rtimes C_{\Gamma(\bG)}(F)$. We have $F_p \in C_{\Gamma(\bG)}(F)$ and so $\Aso(\bG,F)$ is generated by the image of $C_{\Aso(\bG)}(F)$ and the image of $F_p^{-1}$ in $\Aut(\bG^F)$. However, $\langle F_p \rangle \leqslant \Aut(\bG^F)$ is finite, so $C_{\Aso(\bG)}(F)$ contains a preimage of the image of $F_p^{-1}$. This implies $C_{\Aso(\bG)}(F) \to \Aso(\bG,F)$ is surjective.

According to \cite[2.5.14(g)]{gorenstein-lyons-solomon:1998:classification-3} the restriction map
%%%%
\begin{equation*}
\rC_{\Aso(\bG)}(F) \to \Aut(\rO^{p'}(\bG^F))
\end{equation*}
%%%%
is surjective. Therefore, so is the map $\Aso(\bG,F) \to \Aut(\rO^{p'}(\bG^F))$. We wish to now show that $\Aut(\bG^F) = \rK(\bG^F)\Aso(\bG,F)$. To see this it suffices to show that $\rC_{\Aut(\bG^F)}(\Inn(\bG^F))$ contains the kernel of the natural restriction map $\Aut(\bG^F) \to \Aut(\rO^{p'}(\bG^F))$.

Assume $\varphi \in \Aut(\bG^F)$ is such that $\varphi|_{\rO^{p'}(\bG^F)} = \mathrm{Id}_{\rO^{p'}(\bG^F)}$. Then for any $g \in \bG^F$ and $x \in \rO^{p'}(\bG^F)$ we have ${}^gx = \varphi({}^gx) = {}^{\varphi(g)}x$ so $g^{-1}\varphi(g) \in \rC_{\bG^F}(\rO^{p'}(\bG^F))$. By the proof of \cite[Lem.~6.1]{bonnafe:2006:sln}, there exist two unipotent elements $u,v \in \bG^F$ such that $C_{\bG}(u) \cap C_{\bG}(v) = Z(\bG)$. These are clearly contained in $\rO^{p'}(\bG^F)$, so we have $C_{\bG^F}(\rO^{p'}(\bG^F)) \leqslant Z(\bG) \cap \bG^F = Z(\bG)^F$. Hence $\rC_{\bG^F}(\rO^{p'}(\bG^F)) = \rZ(\bG)^F = \rZ(\bG^F)$ so $\varphi$ is a central automorphism.

We now need only show that $\rK(\bG^F) \cap \Aso(\bG,F)$ is trivial. For this, let $\pi : \bG \to \bG_{\ad}$ be an adjoint quotient of $\bG$ and let $F : \bG_{\ad} \to \bG_{\ad}$ be the Frobenius root such that $F\circ\pi = \pi\circ F$. By \cite[Thms.~1.15.6, 2.5.14]{gorenstein-lyons-solomon:1998:classification-3}, restriction through $\pi$ defines injective homomorphisms $\rC_{\Aso(\bG)}(F) \to \rC_{\Aso(\bG_{\ad})}(F)$ and $\Aut(\rO^{p'}(\bG^F)) \to \Aut(\rO^{p'}(\bG_{\ad}^F))$. Suppose $\gamma \in \rC_{\Aso(\bG)}(F)$ restricts to the identity on $\rO^{p'}(\bG^F)$ and let $\bar{\gamma} \in \Aut(\rO^{p'}(\bG_{\ad}^F))$ be the unique automorphism such that $\bar{\gamma}\circ\pi = \pi \circ\gamma$. Then $\bar{\gamma}$ is the identity on $\rO^{p'}(\bG_{\ad}^F)$ so $\bar{\gamma} \in \langle F\rangle$ by \cite[Lem.~2.5.7]{gorenstein-lyons-solomon:1998:classification-3}. Therefore $\gamma \in \langle F \rangle$ so $\gamma$ is the identity on $\bG^F$.
\end{proof}

After (ii) of Lemma~\ref{lem:autos-triv-on-derived}, we see that to understand the action of $\rK(\bG^F)$ on irreducible characters, it is crucial to have a good understanding of central and linear characters. For convenience, we recall a few results from the literature regarding this.

As in \cite{taylor:2018:action-of-automorphisms-symplectic} we denote by $\mathcal{C}(\bG,F)$ the set of pairs $(\bT,\theta)$ consisting of an $F$-stable maximal torus $\bT \leqslant \bG$ and an irreducible character $\theta \in \Irr(\bT^F)$. Dually we denote by $\mathcal{S}(\bG,F)$ the set of all pairs $(\bT,s)$ consisting of an $F$-stable maximal torus $\bT \leqslant \bG$ and a semisimple element $s \in \bT^F$.

\begin{prop}\label{prop:cent-chars}
Assume $\bG$ is a connected reductive algebraic group with Frobenius root $F : \bG \to \bG$. Let $(\bG^{\star},F^{\star})$ be dual to $(\bG,F)$.
%%%%
\begin{enumerate}
	\item For any semisimple element $s \in \bG^{\star F^{\star}}$, there exists a unique character $\omega_s \in \Irr(\rZ(\bG)^F)$ such that $\omega_s = \Res_{\rZ(\bG)^F}^{\bG^F}(\theta)$ for any pair $(\bT,\theta) \in \mathcal{C}(\bG,F)$ corresponding to $(\bT^{\star},s) \in \mathcal{S}(\bG^{\star},F^{\star})$ under duality. Moreover, $\omega_{\chi} = \omega_s$ for all $\chi \in \mathcal{E}(\bG^F,s)$.
	\item There is a group isomorphism $\rZ(\bG^{\star})^{F^{\star}} \to \Irr(\bG^F/\rO^{p'}(\bG^F))$, denoted by $z \mapsto \hat{z}$, such that $-\otimes\hat{z} : \mathcal{E}(\bG^F,s) \to \mathcal{E}(\bG^F,sz)$ is a bijection.
	\item We have a well-defined group isomorphism $\bG^{\star F^{\star}}/\rO^{p'}(\bG^{\star F^{\star}}) \to \Irr(\rZ(\bG)^F)$ given by $s\rO^{p'}(\bG^{\star F^{\star}}) \mapsto \omega_s$, where $s \in \bG^{\star F^{\star}}$ is a semisimple element representing the coset.
\end{enumerate}
\end{prop}

\begin{proof}
(i) is \cite[11.1(d)]{bonnafe:2006:sln} and (ii) is \cite[Prop.~11.4.12, Rem.~11.4.14]{digne-michel:2020:representations-of-finite-groups-of-lie-type}. Dualising the arguments in \cite{digne-michel:2020:representations-of-finite-groups-of-lie-type} gives (iii), see also \cite[Lem.~4.4]{navarro-tiep:2013:characters-relative-degree}.
\end{proof}

\begin{cor}\label{cor:K-series}
Assume $s \in \bG^{\star F^{\star}}$ is a semisimple element and $\sigma = \tau_{\alpha} \in \rK(\bG^F)$. Then there exists a central element $z \in \rZ(\bG^{\star})^{F^{\star}}$ such that $\omega_s\circ\alpha = \hat{z}$ and $\mathcal{E}(\bG^F,s)^{\sigma} = \mathcal{E}(\bG^F,sz)$.
\end{cor}

\begin{proof}
As $\omega_s\circ\alpha$ is a linear character with $\rO^{p'}(\bG^F)$ in its kernel it follows from (ii) of Proposition~\ref{prop:cent-chars} that $\omega_s\circ\alpha = \hat{z}$ for some $z \in \rZ(\bG^{\star})^{F^{\star}}$. If $\chi \in \mathcal{E}(\bG^F,s)$ then $\omega_{\chi} = \omega_s$ by (i) of Proposition~\ref{prop:cent-chars} and so $\chi^{\sigma} = (\omega_{\chi}\circ\alpha)\chi = \hat{z}\chi$ by (ii) of Lemma~\ref{lem:autos-triv-on-derived}. This gives us that $\mathcal{E}(\bG^F,s)^{\sigma} = \mathcal{E}(\bG^F,sz)$ by (ii) of Proposition~\ref{prop:cent-chars}.
\end{proof}

\section{Galois Automorphisms}\label{sec:gal}
\begin{pa}\label{pa:crag-setup}\label{pa:HLnotn}\label{pa:rel-stab}
Let $\bG$ be a connected reductive algebraic group and $F : \bG \to \bG$ a Frobenius endomorphism endowing $\bG$ with an $\mathbb{F}_q$-rational structure. Let $\textbf{T}\leq\textbf{B}$ be, respectively, an $F$-stable maximal torus and Borel subgroup of $\mathbf{G}$. We denote by $\bW = \rN_{\bG}(\bT)/\bT$ the corrresponding Weyl group with simple reflections $\mathbb{S} \subseteq \bW$ determined by $\bB$. Let $\Phi$ be the root system of $\bg{G}$ with respect to $\bg{T}$ and $\Delta \subseteq \Phi^+$ the set of simple and positive roots determined by $\bB$.

For each subset $J \subseteq \Delta$ of simple roots, we have a corresponding standard parabolic and Levi subgroup $\bL_J \leqslant \bP_J \leqslant \bG$. We denote by $\Cusp(\bG,F)$ the set of \emph{cuspidal pairs} $(\bL_J,\lambda)$ with $J \subseteq \Delta$ an $F$-invariant subset of simple roots and $\lambda \in \Irr(\bL_J^F)$ a cuspidal character. To each such pair $(\bL,\lambda) \in \Cusp(\bG,F)$, we have a corresponding Harish-Chandra series $\mathcal{E}(\bG^F,\bL,\lambda)$. By definition, this is the set of irreducible constituents of the corresponding Harish-Chandra induced character $R_{\bL}^{\bG}(\lambda)$. 

Thanks to \cite{geck:1993:a-note-on-harish-chandra} and \cite[Thm.~8.6]{lusztig:1984:characters-of-reductive-groups}, we know that any cuspidal character $\lambda\in\Irr(\bL^F)$ extends to its inertia group $\rN_{\bG^F}(\bL)_\la$. Hence we may fix a so-called \emph{extension map} $\Lambda$ with respect to $\bL^F\lhd \rN_{\bG^F}(\bL)$, meaning that for each cuspidal character $\la\in\Irr(\bL^F)$, $\Lambda(\la)$ is an extension of $\la$ to $\rN_{\bG^F}(\bL)_\la$. We denote by $W(\la) := \rN_{\bG^F}(\bL)_\la/\bL^F$ the so-called relative Weyl group. After choosing a square root $\sqrt{p} \in \overline{\mathbb{Q}} \subseteq \Ql$ of $p$, we have a canonical bijection $\Irr(W(\lambda)) \to \mathcal{E}(\bG^F,\bL_J,\lambda)$, see \cite[Sec.~10.6]{carter:1993:finite-groups-of-lie-type}, which we denote by $\eta \mapsto R_{\bL}^{\bG}(\la)_\eta$.

The group $W(\la)$ is a semidirect product $R(\la) \rtimes C(\la)$ of a Weyl group $R(\la)$ with a root system $\Phi_\la$ and the stabilizer $C(\la)$ in $W(\la)$ of a simple system in $\Phi_\la$, see \cite[Sec.~10.6]{carter:1993:finite-groups-of-lie-type}. Note that $W(\la)$ in general is not necessarily a reflection group but $R(\la)$ is. Given $\sigma\in\gal$ and $\eta \in \Irr(W(\lambda))$ we define $\eta^{(\sigma)}\in\Irr(W(\la^\sigma))$ as in \cite[\S3.5]{schaeffer-fry:2019:galois-automorphisms}. Let $w_0 \in \bW$ be the longest element with respect to $\mathbb{S}$. For any $w \in \bW^F$ we write $\ind(w):=|U\cap U^{w_0w}|$, where $U = \bU^F$ and $\bU \leqslant \bB$ is the unipotent radical of $\bB$.  With this, we now recall \cite[Thm.~3.8]{schaeffer-fry:2019:galois-automorphisms}, which will be key to our proofs of Theorems \ref{mainthmgalSO}(i) and \ref{mainthmgalSp}.
\end{pa}

\begin{thm}[{}{Schaeffer Fry, \cite[Thm.~3.8]{schaeffer-fry:2019:galois-automorphisms}}]\label{thm:GaloisAct}
Assume $\sigma\in\gal$ and $(\bL, \la) \in \Cusp(\bG,F)$ is a cuspidal pair. Let $\Lambda$ be an extension map with respect to $\bL^F \lhd \rN_{\bG^F}(\bL)$. Define:
%%%%
\begin{itemize}
	\item $\delta_{\la,\sigma}$ to be the linear character of $W(\la)$ such that $\Lambda(\la)^\sigma=\delta_{\la,\sigma}\Lambda({\la^\sigma})$;
	\item $\delta'_{\la,\sigma}\in\Irr(W(\la))$ to be the character such that $\delta'_{\la,\sigma}(w)=\delta_{\la,\sigma}(w)$ for $w\in C(\la)$ and $\delta'_{\la,\sigma}(w)=1$ for $w\in R(\la)$; and
	\item $\gamma_{\la,\sigma}$ to be the function on $W(\la)$ such that $\gamma_{\la,\sigma}(w)={\frac{\sqrt{\ind(w_2)}^\sigma}{\sqrt{\ind{(w_2)}}}}$, where $w=w_1w_2$ for $w_1\in R(\la)$ and $w_2\in C(\la)$.
\end{itemize}
%%%%	
Then for any $\eta\in \Irr(W(\la))$, we have
%%%%
\begin{equation*}
\big(R_{\bL}^{\bG}(\la)_\eta\big)^\sigma = R_{\bL}^{\bG}(\lambda^\sigma)_{\eta'},
\end{equation*}
%%%%
where $\eta'\in\Irr(W(\la))=\Irr(W(\la^\sigma))$ is defined by $\eta'(w)=\gamma_{\la,\sigma}(w)\delta'_{\la,\sigma}(w^{-1})\eta^{(\sigma)}(w)$ for each $w\in W(\lambda)$.
\end{thm}

\subsection*{On the group $W(\lambda)$}
\begin{pa}\label{pa:reg-embedding}
For each $w \in \bW$, we fix a representative $n_w \in \rN_{\bG}(\bT)$ such that if $w \in \bW^F$, then $n_w \in \rN_{\bG^F}(\bT)$. If $J \subseteq \Delta$ is an $F$-stable subset of simple roots, then we have an isomorphism $\rN_{\bW^F}(J) \to \rN_{\bG^F}(\bL_J)/\bL_J^F$, given by $w \mapsto n_w\bL_J^F$, where $\rN_{\bW^F}(J) = \{w \in \bW^F \mid wJ = J\}$. Given a cuspidal pair $(\bL_J,\lambda) \in \Cusp(\bG,F)$, the group $W(\lambda) \leqslant \rN_{\bG^F}(\bL_J)/\bL_J^F$ can thus be identified with the subgroup $\rN_{\bW^F}(J)_{\lambda} \leqslant \rN_{\bW^F}(J)$, and hence with a subgroup of $\bW$. We will make this identification implicitly in what follows.

Now assume $\iota : \bG \to \wt{\bG}$ is an $F$-equivariant regular embedding of $\bG$ into a group with connected centre and let us identify $\bG$ with its image $\iota(\bG)$. Given $J \subseteq \Delta$, we have a corresponding $F$-stable Levi $\wt{\bL}_J = \bL_J\cdot \rZ(\wt{\bG})$ of $\wt{\bG}$. This is standard with respect to the maximal torus $\wt{\bT} = \bT\cdot \rZ(\wt{\bG})$ and Borel $\wt{\bB} = \bB\cdot \rZ(\wt{\bG})$ of $\wt{\bG}$, see \cite[Prop.~3.4.10]{digne-michel:2020:representations-of-finite-groups-of-lie-type}. If $\lambda \in \Irr(\bL_J^F)$ is a cuspidal character, then any character $\wt{\lambda} \in \Irr(\wt{\bL}_J^F)$ covering $\lambda$ is also cuspidal. Note that the Weyl group $\wt{\bW} = \rN_{\wt{\bG}}(\wt{\bT})/\wt{\bT}$ is naturally isomorphic to $\bW$. We let $W(\wt{\lambda}) = \rN_{\bG^F}(\bL_J)_{\wt{\lambda}}/\bL_J^F$, which is a subgroup of $W(\lambda)$ by \cite[Cor.~12.5(c)]{bonnafe:2006:sln}.
\end{pa}

\begin{lem}\label{lem:stab-lam-tilde}
Assume $(\bL,\lambda) \in \Cusp(\bG,F)$ is a cuspidal pair and $\wt{\lambda} \in \Irr(\wt{\bL})$ covers $\lambda$. Then $R(\lambda) = W(\wt{\lambda})$ and there is an injective homomorphism $C(\lambda) \to H^1(F,\mathcal{Z}(\bG))$ where $\mathcal{Z}(\bG) = \rZ(\bG)/Z^{\circ}(\bG)$.
\end{lem}

\begin{proof}
Let $\wt{W}(\wt{\lambda}) = \wt{R}(\wt{\lambda}) \rtimes \wt{C}(\wt{\lambda})$ be as in \ref{pa:rel-stab} with respect to $\wt{\bG}^F$. By a result of Lusztig, $\wt{C}(\wt{\lambda}) = \{1\}$, see \cite[Thm.~8.6]{lusztig:1984:characters-of-reductive-groups}. For the first statement, it suffices to show that $\wt{R}(\wt{\lambda}) = R(\lambda)$, where we identify $\wt{W}(\wt{\lambda})$ and $W(\wt{\lambda})$. For this, let $\bT \leqslant \bP \leqslant \bG$ be an $F$-stable parabolic subgroup with Levi complement $\bT \leqslant \bM \leqslant \bP$. We assume $\bL \subsetneq \bM$ and $\bM$ is minimal with this property. As above, we set $\wt{\bM} = \bM\cdot \rZ(\wt{\bG})$, so $\wt{\bL} \subsetneq \wt{\bM}$.

If $W_{\bM^F}(\bL) := \rN_{\bM^F}(\bL)/\bL^F$ is trivial then there is nothing to consider. So suppose $W_{\bM^F}(\bL)$ is non-trivial, and hence so is $W_{\wt{\bM}^F}(\wt{\bL}) := \rN_{\wt{\bM}^F}(\wt{\bL})/\wt{\bL}^F$. By minimality, $W_{\bM^F}(\bL) \cong W_{\wt{\bM}^F}(\wt{\bL})$ is cyclic of order $2$. We will assume that $W_{\bM^F}(\bL)_{\lambda} = \rN_{\bM^F}(\bL)_{\lambda}/\bL^F$ is non-trivial. Then $R_{\bL}^{\bM}(\lambda) = \lambda_1+\lambda_2$ with $\lambda_1(1) \leqslant \lambda_2(1)$ and $\lambda_i \in \Irr(\bM^F)$. The unique non-trivial element $1 \neq w \in W_{\bM^F}(\bL)_{\lambda}$ is a generator of $R(\lambda)$ if and only if $\lambda_1(1) \neq \lambda_2(1)$.

As Harish-Chandra induction is compatible with restriction and conjugation, we have
%%%%
\begin{equation}\label{eq:res-1}
\Res_{\bM^F}^{\wt{\bM}^F}\left(R_{\wt{\bL}}^{\wt{\bM}}(\wt{\lambda})\right) = R_{\bL}^{\bM}\left(\Res^{\bM^F}_{\bL^F}(\wt{\lambda})\right) = \sum_{g \in \wt{\bL}^F/(\wt{\bL}^F)_{\lambda}} R_{\bL}^{\bM}(\lambda)^g = \sum_{g \in \wt{\bL}^F/(\wt{\bL}^F)_{\lambda}} \lambda_1^g+\lambda_2^g.
\end{equation}
%%%%
If $\lambda_1(1) \neq \lambda_2(1)$ then these characters cannot be conjugate in $\wt{\bM}^F$, so by Clifford's Theorem $R_{\wt{\bL}}^{\wt{\bM}}(\wt{\lambda}) = \wt{\lambda}_1 + \wt{\lambda}_2$ with $\wt{\lambda}_i \in \Irr(\wt{\bM}^F)$ an irreducible character with degree $[\wt{\bL}^F : (\wt{\bL}^F)_{\lambda}]\lambda_i(1)$. Hence, $\wt{\lambda}_1(1) \neq \wt{\lambda_2}(1)$ so $R(\lambda) \leqslant R(\wt{\lambda})$.

Conversely, if $W_{\bM^F}(\bL)_{\wt{\lambda}} = \rN_{\bM^F}(\bL)_{\wt{\lambda}}/\bL^F$ is non-trivial, then $R_{\wt{\bL}}^{\wt{\bM}}(\wt{\lambda}) = \wt{\lambda}_1 + \wt{\lambda}_2$ with $\wt{\lambda}_i \in \Irr(\wt{\bM}^F)$. We must have $W_{\bM^F}(\bL)_{\lambda}$ is also non-trivial, by \cite[Cor.~12.5(c)]{bonnafe:2006:sln}, so $R_{\bL}^{\bM}(\lambda) = \lambda_1 + \lambda_2$ as above. We can compare
%%%%
\begin{equation*}
\Res_{\bM^F}^{\wt{\bM}^F}\left(R_{\wt{\bL}}^{\wt{\bM}}(\wt{\lambda})\right)  = \Res_{\bM^F}^{\wt{\bM}^F}(\wt{\lambda}_1 + \wt{\lambda}_2)
\end{equation*}
%%%%
with the terms in \eqref{eq:res-1}. If $\lambda_1$ is a constituent of $\Res_{\bM^F}^{\wt{\bM}^F}(\wt{\lambda}_i)$, then $\sum_{g \in \wt{\bL}^F/(\wt{\bL}^F)_{\lambda}} \lambda_1^g$ must be also, as $\wt{\bL}^F \leqslant \wt{\bM}^F$ and the restriction is invariant under $\wt{\bM}^F$.

If $\lambda_1$ and $\lambda_2$ were both constituents of $\Res_{\bM^F}^{\wt{\bM}^F}(\wt{\lambda}_1)$, say, then we would have to have $\Res_{\bM^F}^{\wt{\bM}^F}(\wt{\lambda}_2) = 0$, which is impossible. If $\wt{\lambda}_1(1) \neq \wt{\lambda}_2(1)$ then $\Res_{\bM^F}^{\wt{\bM}^F}(\wt{\lambda}_i) = \sum_{g \in \wt{\bL}^F/(\wt{\bL}^F)_{\lambda}} \lambda_i^g$ so $\tilde{\lambda}_i(1) = [\wt{\bL}^F : (\wt{\bL}^F)_{\lambda}]\lambda_i(1)$ and $\lambda_1(1) \neq \lambda_2(1)$. This gives $R(\wt{\lambda}) \leqslant R(\lambda)$.

For the final statement, note that if $g \in \rN_{\bG^F}(\bL)_{\lambda}$, then $\wt{\lambda}^g \in \Irr(\wt{\bL}^F)$ also covers $\lambda$. Thus $\wt{\lambda}^g = \wt{\lambda}\theta$ where $\theta$ is the inflation of a linear character of $\wt{\bL}^F/\bL^F\rZ(\wt{\bL})^F \cong H^1(F,\mathcal{Z}(\bG))$. This gives a homomorphism $W(\lambda) \to H^1(F,\mathcal{Z}(\bG))$ with kernel $W(\wt{\lambda}) = R(\lambda)$.
\end{proof}

\subsection*{On the Characters $\gamma_{\la,\sigma}$ and $\delta_{\la,\sigma}$}
In the latter parts of this paper, we will need to carefully study the characters $\delta_{\lambda,\sigma}$ and $\gamma_{\lambda,\sigma}$ introduced in Theorem~\ref{thm:GaloisAct}. Here we introduce a few general statements that will be used later on. For the rest of this section, we fix a cuspidal pair $(\bL,\lambda) \in \Cusp(\bG,F)$ and a Galois automorphism $\sigma \in \gal$.

\begin{lem}\label{lem:genrsigma}
Assume $(\bL, \la)\in\Cusp(\bG,F)$ and $\sigma\in\gal$. Then
$\gamma_{\la,\sigma}$ is a character of $W(\la)$, and is moreover trivial if and only if at least one of the following holds:
\begin{itemize}
\item $q$ is a square;
\item the length $l(w_2)$, with respect to $(\bW,\mathbb{S})$, is even for each $w_2\in C(\la)$; or
\item $\sigma$ fixes $\sqrt{p}$.
\end{itemize}
Otherwise, $\gamma_{\la,\sigma}(w)=(-1)^{l(w_2)}$, where $w=w_1w_2$ with $w_1\in R(\la)$ and $w_2\in C(\la)$.
\end{lem}
\begin{proof}
By \cite[\S2.9]{carter:1993:finite-groups-of-lie-type} we have $\ind(w_2)=q^{l(w_2)}$ with $l(w_2)$ as above. Therefore, $\gamma_{\la,\sigma}$ is a character in the case of finite reductive groups.  In particular, $\gamma_{\la,\sigma}(w)=1$ when $\sqrt{q^{l(w_2)}}$ is fixed by $\sigma$ and $\gamma_{\la,\sigma}(w)=-1$ otherwise.
\end{proof}

\begin{pa}
When $\bL=\bT$, the group $\rN_{\bG^F}(\bg{T})$ is generated by $\bT^F$ and the group $\langle n_\alpha(\pm1) \mid \alpha\in\Phi \rangle$, see \cite[Lem.~3.2]{malle-spaeth:2016:characters-of-odd-degree}. Let $R_\la$ be the subgroup of the stabilizer $\rN_{\bG^F}(\bg{T})_\la$ of $\la$ in $\rN_{\bG^F}(\bT)$ generated by $\bT^F$ and $\langle n_\alpha(-1)\mid \alpha\in\Phi_\la\rangle$, so that $R_\la/\bT^F \cong R(\la)$.   Here we use the notation of the Chevalley generators as in \cite[\S 1.12]{gorenstein-lyons-solomon:1998:classification-3}. The next lemma concerns the restriction of $\Lambda(\la)$ to $R_\la$.
\end{pa}

\begin{lem}\label{lem:extRla}
Assume that $\bg{G}$ is simple of simply connected type, not of type $A_n$, and let $\la\in\Irr(\bT^F)$.  Then the following hold: 
\begin{enumerate}
\item any extension $\wh{\la}$ of $\la$ to $R_\la$ satisfies that $\wh{\la}(x)\in\{\pm1\}$ for all $x\in \langle n_\alpha(-1)\mid \alpha\in\Phi_\la\rangle$;
\item if $\la^\sigma=\la$, then $\wh{\la}^\sigma=\wh{\la}$ for any extension $\wh{\la}$ of $\la$ to $R_\la$.  In particular, $R(\la)\leq \Ker(\delta_{\la,\sigma})$.
\end{enumerate}
 \end{lem}

\begin{proof}
Let $\wh{\la}$ be an extension of $\lambda$ to $R_\la$.  By \cite[Lem.~5.1]{malle-spaeth:2016:characters-of-odd-degree}, $\Phi_\la$ consists of those $\alpha\in\Phi$ such that $\lambda(h_\alpha(t))=1$ for each $t\in\mathbb{F}_q^\times$.  In particular, this means that if $\alpha\in\Phi_\la$, then $\la(h_\alpha(-1))=1$. But since $n_\alpha(-1)^2=h_\alpha(-1)$, it follows that $\wh{\la}(n_\alpha(-1))\in\{\pm1\}$, proving (i).
For (ii), notice that
\[\wh{\la}^\sigma(tx)=\wh{\la}^\sigma(t)\wh{\la}^\sigma(x)=\la^\sigma(t)\wh{\la}(x)=\la(t)\wh{\la}(x)=\wh{\la}(tx),\]  
for $tx\in R_\la$ with $t\in \bT^F$ and $x\in \langle n_\alpha(-1)\mid \alpha\in\Phi_\la\rangle$, since $\wh{\la}(x)\in\{\pm1\}$ and $\la^\sigma=\la$.
\end{proof}

\subsection*{Square Roots}

\begin{pa}\label{pa:j-inj-hom}
Let $p$ be an odd prime.  In light of Lemma \ref{lem:genrsigma}, we will want to study the action of $\gal$ on $\sqrt{p}$.   Let us be more explicit about our choice of $\sqrt{p}$. We fix an injective group homomorphism $\jmath : \mathbb{Q/Z} \to \overline{\mathbb{Q}}^{\times}$ and let $\tilde{\jmath} : \mathbb{Q} \to \overline{\mathbb{Q}}^{\times}$ be the composition with the natural quotient map $\mathbb{Q} \to \mathbb{Q/Z}$. We have $i = \tilde{\jmath}(1/4)$ is an element such that $i^2 = -1$. We also denote this by $\sqrt{-1}$. Letting $\omega = (-1)^{\frac{p-1}{2}}$, there then exists a unique square root $\sqrt{p}$ satisfying   $\sqrt{\omega}\sqrt{p} = \sum_{n=1}^{p-1}\legendre{n}{p}\zeta_p^n$, where $\legendre{n}{p}$ denotes the Legendre symbol and $\zeta_p := \tilde{\jmath}(\frac{1}{p})$, see \cite[\S36]{bonnafe:2006:sln}. We set $\sqrt{\omega p} := \sqrt{\omega}\sqrt{p}$.

Now, let $\sigma\in\gal$ and let $k\in \mathbb{Z}$ be an integer coprime to $p$ such that $\xi^\sigma=\xi^k$ for all $p$th roots of unity $\xi\in\overline{\Q}$. Then we have

\begin{equation}\label{eq:gaussum}
\sqrt{\omega p}^\sigma=\sum_{n=1}^{p-1}{\legendre{n}{p}}\zeta_p^{k n}=\legendre{k}{p}\cdot\sum_{n=1}^{p-1}{\legendre{k n}{p}}\zeta_p^{k n}={\legendre{k}{p}}\sqrt{\omega p}.
\end{equation}
\end{pa}

\section{Galois Automorphisms Relevant for the McKay--Navarro Conjecture}\label{sec:actionH}
\begin{pa}
In this section, we turn our attention to the case of the specific Galois automorphisms that appear in the context of the McKay--Navarro conjecture, with the aim of recording several statements that will also likely be useful for future work with the conjecture for finite reductive groups. Let $\ell$ be a prime and denote by $\galh:=\galh_\ell$ the subgroup of elements $\sigma\in\gal$ such that there is some integer $r\geq0$ satisfying $\zeta^{\sigma}=\zeta^{\ell^r}$ for all roots of unity $\zeta$ of order coprime to $\ell$.  Throughout, we fix $\ell$ and $\sigma\in\galh$ and understand $r$ to be this integer corresponding to $\sigma$. 

Assume $p\neq \ell$ is an odd prime. We now discuss the action of $\galh$ on $\sqrt{p}$ in order to describe the character $\gamma_{\la,\sigma}$ in more detail in the case that $\sigma\in\galh$.   Note that ${\legendre{\ell^r n}{p}}={\legendre{\ell}{p}^r}{\legendre{n}{p}}$.  Then \eqref{eq:gaussum} implies $\sqrt{\omega p}^\sigma={\legendre{\ell}{p}}^r\sqrt{\omega p}$.
\end{pa}

\subsection*{The Case $\ell$ Odd}

\begin{pa}
First, suppose that $\ell$ is odd.  Note that if $r$ is even, then $\ell^r\equiv 1\pmod 4$, so $\sigma$ fixes $4$th roots of unity and ${\legendre{\ell}{p}}^r=1$ since $\legendre{\ell}{p}\in\{\pm1\}$.  If $r$ is odd, then $\ell^r\equiv \ell\pmod 4$.  We also remark that when $\ell\neq p$ are odd primes, we have
%%%%

\begin{equation*}\legendre{\ell}{p}=(-1)^{(p-1)(\ell-1)/4}\legendre{p}{\ell}.
\end{equation*}
%%%%
Together with \eqref{eq:gaussum} this yields the following:
\end{pa}

\begin{lem}\label{lem:omegapH}
Assume $\ell \neq p$ are both odd and $\sigma\in\galh$. Let $\omega = (-1)^{\frac{p-1}{2}}$, so  $p\equiv\omega\pmod4$.
%%%%
\begin{enumerate}
\item If $r$ is even, then $\sqrt{p}^{\sigma}=\sqrt{p}$ and $\sqrt{\omega p}^{\sigma}=\sqrt{\omega p}$.
\item If $r$ is odd, then $\sqrt{\omega p}^{\sigma} = \legendre{\ell}{p}\sqrt{\omega p}$, so $\sqrt{p}^\sigma=(-1)^{(p-1)(\ell-1)/4}\legendre{\ell}{p}\sqrt{p}=\legendre{p}{\ell}\sqrt{p}.$
\end{enumerate}
\end{lem}

We can now apply this is to the setup of \ref{pa:crag-setup}. In particular, we have $\bG^F$ is a finite reductive group defined over $\mathbb{F}_q$ with $q$ a power of $p$. Combining Lemmas~\ref{lem:genrsigma} and \ref{lem:omegapH}, we have the following.

\begin{lem}\label{lem:genrsigmaH}
Assume $\ell \neq p$ are both odd and $(\bL,\lambda) \in \Cusp(\bG,F)$ is a cuspidal pair. Let $\sigma\in\galh$ and $w=w_1w_2\in W(\lambda)$ with $w_1\in R(\la)$ and $w_2\in C(\la)$.
%%%%
\begin{enumerate}
\item If $q$ is a square or $r$ is even, then $\gamma_{\la,\sigma}=1$.
\item If $q$ is not a square and $r$ is odd, then
%%%%
\begin{equation*}
\gamma_{\la,\sigma}(w) = \legendre{p}{\ell}^{l(w_2)} = \begin{cases}
(-1)^{l(w_2)}\legendre{\ell}{p}^{l(w_2)} &\text{if $p\equiv \ell\equiv 3\pmod 4$}\\
\legendre{\ell}{p}^{l(w_2)} &\text{otherwise},
\end{cases}
\end{equation*}
%%%%
where $l(w_2)$ is the length of $w_2$ in the Weyl group $(\bW,\mathbb{S})$ of $\bg{G}$.
\end{enumerate}
\end{lem}

\begin{pa}
The particular case that $\ell\mid(q-1)$ has been studied in various contexts related to (refinements of) the McKay conjecture, see, e.g., \cite{malle-spaeth:2016:characters-of-odd-degree, CSS21}, due to the relatively nice nature of height-zero characters in this case.  
We remark that if $\ell$ is a  prime dividing $(q-1)$, we have $q\equiv 1\pmod \ell$ is a square modulo $\ell$, and hence either $q$ is a square or $p$ is a square modulo $\ell$.  Then Lemma~\ref{lem:genrsigmaH}  yields:  
\end{pa}

\begin{cor}\label{cor:rsigmalinearprime}
Assume $\ell$ and $p$ are odd and $\ell \mid (q-1)$. Then for any cuspidal pair $(\bL,\la)\in\Cusp(\bG, F)$ and $\sigma\in\galh$, the character $\gamma_{\la,\sigma}$ is trivial.
\end{cor}

\subsection*{The Case $\ell=2$}
\begin{pa}
We now focus on the situation when $\ell=2$ and analyse the character $\gamma_{\la, \sigma}$ in this case. We begin by noting that
%%%%
\begin{equation*} 
\legendre{2}{p} = \begin{cases}
1 &\text{if $p\equiv\pm1\pmod 8$},\\
-1 &\text{if $p\equiv \pm3\pmod 8$},
\end{cases}
\end{equation*}
%%%%
so comparing with \eqref{eq:gaussum} yields:
\end{pa}

\begin{lem}\label{lem:omegapH2}
Let $\sigma\in\galh$ with $\ell=2$ and let $p$ be an odd prime. Let $\omega = (-1)^{\frac{p-1}{2}}$, so $p\equiv \omega\pmod 4$.

\begin{enumerate}
\item If $p\equiv \pm 1\pmod 8$, then $\sqrt{\omega p}^{\sigma} = \sqrt{\omega p}$.
\item If $p\equiv \pm 3\pmod 8$, then $\sqrt{\omega p}^{\sigma} = (-1)^r\sqrt{\omega p}$.
\end{enumerate}
\end{lem}

In \cite[Lem.~4.10]{schaeffer-fry:2019:galois-automorphisms}, the character $\gamma_{\la,\sigma}$ is described explicitly in the case of a specific choice of $\sigma$. We  now record an extension of that statement to all of $\galh$.

\begin{lem}\label{lem:rw}
Assume $\ell=2$, $p$ is odd, and $(\bL, \la)\in\Cusp(\bG, F)$. Let $\sigma\in\galh$ and $w=w_1w_2\in W(\lambda)$ with $w_1\in R(\la)$ and $w_2\in C(\la)$.
%%%%
\begin{enumerate}
\item If $q\equiv 1\pmod 8$ then $\gamma_{\la,\sigma}(w)=1$.
\item If $q\equiv -1\pmod 8$ then
$\gamma_{\la,\sigma}(w)=
\begin{cases}
1 &\text{if $i^\sigma=i$},\\
(-1)^{l(w_2)} &\text{if $i^\sigma=-i$}.
\end{cases}$

\item If $q\equiv 3\pmod 8$ then 
$\gamma_{\la,\sigma}(w) =
\begin{cases}
1 &\text{if $i^\sigma=(-1)^r i$},\\
(-1)^{l(w_2)} &\text{if $i^\sigma=(-1)^{r-1}i$}.
\end{cases}$

\item If $q\equiv -3\pmod 8$ then $\gamma_{\la,\sigma}(w)=(-1)^{r \cdot l(w_2)}$.
\end{enumerate}
%%%%
Here $l(w_2)$ denotes the length of $w_2$ in the Weyl group $(\bW,\mathbb{S})$ of $\bg{G}$.
\end{lem}

\begin{proof}
This follows from Lemma~\ref{lem:genrsigma} and Lemma~\ref{lem:omegapH2}.
\end{proof}

\begin{pa}
We conclude with one further remark.  When $\ell=2$, the power $r$ defined by $\sigma\in\galh$ is even if and only if $\sigma$ fixes third roots of unity.  Hence in this case, the action of $\sigma$ on $\sqrt{p}$ and descriptions of $\gamma_{\la,\sigma}$  obtained in Lemmas~\ref{lem:omegapH2} and \ref{lem:rw} can alternatively be described entirely in terms of $p\pmod 8$ and the action of $\sigma$ on third and fourth roots of unity.  The same will be true in Corollary~\ref{cor:rdelta} below, where we describe $\gamma_{\la,\sigma}\delta_{\la,\sigma}$ in the case of $\Sp(V)$.
\end{pa}

\section{Power Maps and Regular Unipotent Elements}\label{sec:power-maps-1}
\begin{pa}
In this section, we return to the setup of \ref{pa:crag-setup}, so that $\bG$ is a connected reductive algebraic group. We fix an integer $k \in \mathbb{Z}$ coprime to $p$.  Then the power map $g \mapsto g^k$ on $\bG$ induces a permutation of the set of unipotent conjugacy classes $\Clu(\bG)$. In fact, this map is known to be the identity. In other words, any unipotent element is \emph{rational}. This was shown, for instance, by Tiep--Zalesski \cite[Thm.~4.3]{tiep-zalesski:2004:unipotent-elements}, with other proofs given by Lusztig \cite[Prop.~2.5(a)]{lusztig:2009:remarks-on-springers-reps} and Liebeck--Seitz \cite[Cor.~3]{liebeck-seitz:2012:unipotent-nilpotent-classes}.

We now consider the analogous question over $\mathbb{F}_q$. Let $\mathcal{O} \in \Clu(\bG)^F$ be an $F$-stable unipotent class. Then the power map induces a permutation $\pi_k : \Cl_{\bG^F}(\mathcal{O}^F) \to \Cl_{\bG^F}(\mathcal{O}^F)$ of the set of $\bG^F$-conjugacy classes contained in the fixed point set $\mathcal{O}^F$. If $\bG$ is simple and simply connected, then \cite[Thm.~1.7]{tiep-zalesski:2004:unipotent-elements} describes exactly when every unipotent element of $\bG^F$ is rational, which means $\pi_k$ is the identity for any class $\mathcal{O} \in \Clu(\bG)^F$.

Here we wish to describe the actual permutation.  Unlike \cite{tiep-zalesski:2004:unipotent-elements}, we will focus just on the case of good characteristic. Note that after \cite[Thm.~1.5]{tiep-zalesski:2004:unipotent-elements}, we know this permutation is a product of transpositions. In this section, we consider the case of regular unipotent elements. This will be used to settle the general case for symplectic and special orthogonal groups in Section~\ref{sec:power-maps-2}. If $\bG=\SO(V)$, then every unipotent element of $\bG^F$ is rational by \cite[Thm.~1.9]{tiep-zalesski:2004:unipotent-elements}, even strongly rational, so the main case of interest for us is when $\bG = \Sp(V)$. However, our argument below also covers the case of $\SO(V)$.

For any $x \in \bG^F$ we denote by $H^1(F,A_{\bG}(x))$ the $F$-conjugacy classes of $A_{\bG}(x)$. The $F$-class of $g\rC_{\bG}^{\circ}(x) \in A_{\bG}(x)$ is denoted by $[g\rC_{\bG}^{\circ}(x)]$. If $\mathscr{L} : \bG \to \bG$ is the Lang map, defined by $\mathscr{L}(g) = g^{-1}F(g)$, then we have a bijection $\Cl_{\bG^F}(\mathcal{O}^F) \to H^1(F,A_{\bG}(x))$ given by ${}^gx \mapsto [\mathscr{L}(g)\rC_{\bG}^{\circ}(x)]$, where $\mathcal{O}$ is the $\bG$-class of $x$. We leave the proof of the following to the reader.
\end{pa}

\begin{lem}\label{lem:perm-conj}
Assume $\mathcal{O} \in \Clu(\bG)^F$ is an $F$-stable unipotent class and $u\in \mathcal{O}^F$. If $g \in \bG$ is such that $u^k = {}^gu$, then we have a commutative diagram
%%%%
\begin{center}
\begin{tikzcd}
H^1(F,A_{\bG}(u)) \arrow[r]\arrow[d, "\tau"] & \Cl_{\bG^F}(\mathcal{O}^F) \arrow[d, "\pi_k"]\\
H^1(F,A_{\bG}(u)) \arrow[r] & \Cl_{\bG^F}(\mathcal{O}^F)
\end{tikzcd}
\end{center}
%%%%
where $\tau([a\rC_{\bG}^{\circ}(u)]) = [g^{-1}aF(g)\rC_{\bG}^{\circ}(u) ]$. Moreover, if there is an $F$-stable torus $\bS \leqslant \bG$ such that $\rC_{\bG}(u) = \bS\rC_{\bG}^{\circ}(u)$ then $\tau([a\rC_{\bG}^{\circ}(u)]) = [a\mathscr{L}(g)\rC_{\bG}^{\circ}(u)]$.
\end{lem}

\begin{pa}\label{pa:torus-borel}
Let $X = X(\bT)$ and $\wc{X} = \wc{X}(\bT)$ be the character and cocharacter groups of $\bT$ with the usual perfect pairing $\langle -,-\rangle : X \times \wc{X} \to \mathbb{Z}$. If $\mathbb{Z}_{(p)}$ is the localisation of $\mathbb{Z}$ at the prime ideal $(p)$, then we may choose an isomorphism $\imath : \mathbb{Z}_{(p)}/\mathbb{Z} \to \mathbb{F}^{\times}$. With respect to this choice, we have a surjective homomorphism of abelian groups $\tilde{\imath} : \mathbb{Q}\wc{X} \to \bT$, as in \cite[\S3]{bonnafe:2005:quasi-isolated}, where $\mathbb{Q}\widecheck{X} = \mathbb{Q} \otimes_{\mathbb{Z}} \wc{X}$ is a $\mathbb{Q}$-vector space. If $\mathbb{Q}X = \mathbb{Q}\otimes_{\mathbb{Z}} X$ then $\langle-,-\rangle$ extends naturally to a non-degenerate bilinear form $\mathbb{Q}X \times \mathbb{Q}\wc{X} \to \mathbb{Q}$.

For each root $\alpha \in \Phi$ we have a corresponding coroot $\wc{\alpha} \in \wc{X}$. Let $\mathbb{Q}\widecheck{\Phi} \subseteq \mathbb{Q}\wc{X}$ be the subspace spanned by $\widecheck{\Phi} = \{\wc{\alpha} \mid \alpha \in \Phi\}$. Then we have a set of fundamental dominant coweights $\wc{\Omega} = \{\wc{\omega}_{\alpha} \mid \alpha \in \Delta\} \subseteq \mathbb{Q}\wc{\Phi}$ defined such that $\langle \alpha,\widecheck{\omega}_{\beta}\rangle = \delta_{\alpha,\beta}$. By \cite[\S1, no.~10]{bourbaki:2002:lie-groups-chap-4-6} we have a vector
%%%%
\begin{equation*}
\wc{\rho}_{\Delta} := \sum_{\alpha \in \Delta} \wc{\omega}_{\alpha} = \frac{1}{2}\sum_{\alpha \in \Phi^+} \widecheck{\alpha} \in \mathbb{Z}\wc{\Omega} \subseteq \mathbb{Q}\wc{\Phi}.
\end{equation*}
%%%%
In Table~\ref{tab:half-sum}, we describe the image of $\wc{\rho}_{\Delta}$ in the fundamental group $\wc{\Pi} := \mathbb{Z}\widecheck{\Omega}/\mathbb{Z}\widecheck{\Phi}$ when $\Phi$ is irreducible. Here we use the notation in the plates of \cite{bourbaki:2002:lie-groups-chap-4-6}. For the cases $\G_2$, $\F_4$, $\E_6$, and, $\E_8$ we have $\wc{\rho}_{\Delta}$ is always contained in $\mathbb{Z}\wc{\Phi}$.
\end{pa}

\begin{table}[t]
\centering
\begin{tabular}{*{8}{c}}
\toprule
$\Phi$ & $\A_{2m+\epsilon}$ & $\B_{2m+\epsilon}$ & $\C_n$ & $\D_{2m+\epsilon}$ & $\E_7$\\
\midrule
$\widecheck{\rho}_{\Delta}$ & $\epsilon\widecheck{\omega}_m$ & $(m+\epsilon)\widecheck{\omega}_1$ & $\widecheck{\omega}_n$ & $(m+\epsilon)\widecheck{\omega}_1$ & $\widecheck{\omega}_7$\\
\bottomrule
\end{tabular}
\caption{Representative of the image of $\widecheck{\rho}_{\Delta}$ in the fundamental group $\mathbb{Z}\widecheck{\Omega}/\mathbb{Z}\widecheck{\Phi}$.}
\label{tab:half-sum}
\end{table}

\begin{prop}\label{prop:coprime-perm-reg}
Assume $p$ is good for $\bG$ and $u \in \bB^F$ is a regular unipotent element with $\mathcal{O}$ the $\bG$-class containing $u$. If $r \in \mathbb{Z}$ is such that $\imath(r/(q-1)+\mathbb{Z}) = k\pmod{p}$, then there is an element $g \in \bB$ such that ${}^gu = u^k$ and $\mathscr{L}(g)\rC_{\bG}^{\circ}(u) = z\rC_{\bG}^{\circ}(u)$, where $z := \tilde{\imath}(r\widecheck{\rho}_{\Delta})$ is contained in $\rZ(\bG)$. If $\tau$ is as in Lemma~\ref{lem:perm-conj}, then $\tau([a\rC_{\bG}^{\circ}(u)]) = [az\rC_{\bG}^{\circ}(u)]$.
\end{prop}

\begin{proof}
Let $\bU \leqslant \bB$ be the unipotent radical of $\bB$. For each $\alpha \in \Delta$ we pick a closed embedding $x_{\alpha} : \mathbb{F} \to \bU$ onto the corresponding root subgroup. Denote by $\bU_{\der} \lhd \bU$ the derived subgroup of $\bU$ and let $\overline{u} = \prod_{\alpha \in \Delta} x_{\alpha}(c_{\alpha})\bU_{\der}$, with $c_{\alpha} \in \mathbb{F}$, be the image of $u$ in $\bU/\bU_{\der}$. Clearly $\overline{u}^k = \prod_{\alpha \in \Delta} x_{\alpha}(kc_{\alpha})\bU_{\der}$. If $t = \tilde{\imath}(r/(q-1) \cdot \widecheck{\rho}_{\Delta})$, then we have $\alpha(t) = k$ for all $\alpha \in \Delta$. Note that $\wc{\rho}_{\Delta}$ is invariant under all graph automorphisms and
%%%%
\begin{equation*}
\mathscr{L}(t) = \tilde{\imath}\left(\frac{(q-1)r}{(q-1)} \cdot \widecheck{\rho}_{\Delta}\right) = \tilde{\imath}(r\widecheck{\rho}_{\Delta}) = z.
\end{equation*}
%%%%
Moreover, $z \in \rZ(\bG)$ because  $\alpha(z) = \imath(\langle \alpha,r\widecheck{\rho}_{\Delta}\rangle + \mathbb{Z}) = 1$ for any root $\alpha \in \Phi$, as $\langle \alpha,r\widecheck{\rho}_{\Delta}\rangle = r\langle \alpha,\widecheck{\rho}_{\Delta}\rangle \in \mathbb{Z}$.

It is clear that ${}^t\overline{u} = \overline{u}^k$, so $t^{-1}(u^k)t \in u\bU_{\der}$. It is stated in the proof of \cite[Prop.~12.2.2]{digne-michel:2020:representations-of-finite-groups-of-lie-type} that $u\bU_{\der}$ is the $\bU$-conjugacy class of $u$, see also \cite[Prop.~12.2.5]{digne-michel:2020:representations-of-finite-groups-of-lie-type}. Hence, there exists $v \in \bU$ such that $t^{-1}(u^k)t = {}^vu$ so ${}^{tv}u = u^k$. As $\mathscr{L}(t) \in \rZ(\bG)$, we have $\mathscr{L}(tv) = \mathscr{L}(v)\mathscr{L}(t) \in \rC_{\bG}(u)$, which implies $\mathscr{L}(v) \in \rC_{\bU}(u) \leqslant \rC_{\bG}^{\circ}(u)$ by \cite[Prop.~12.2.7]{digne-michel:2020:representations-of-finite-groups-of-lie-type} because $p$ is good. Hence $\rC_{\bG}^{\circ}(u)\mathscr{L}(tv) = \rC_{\bG}^{\circ}(u)z$.
\end{proof}

\begin{pa}\label{rem:kmodp}
Let us maintain the assumptions of Proposition~\ref{prop:coprime-perm-reg}. As was already used above, the natural map $\rZ(\bG) \to \rC_{\bG}(u)$ factors through an $F$-equivariant isomorphism $\mathcal{Z}(\bG) \to A_{\bG}(u)$, where $\mathcal{Z}(\bG) = \rZ(\bG)/\rZ^{\circ}(\bG)$. This defines a bijection $H^1(F,\mathcal{Z}(\bG)) \to H^1(F,A_{\bG}(u))$. As $\mathcal{Z}(\bG)$ is abelian, we have $H^1(F,\mathcal{Z}(\bG)) = \mathcal{Z}(\bG)/\mathscr{L}(\mathcal{Z}(\bG))$. Hence, we have $u$ and $u^k$ are $\bG^F$-conjugate if and only if $z \in \mathscr{L}(\mathcal{Z}(\bG))$, where $z = \tilde{\imath}(r\wc{\rho}_{\Delta})$ is as in Proposition~\ref{prop:coprime-perm-reg}.

Recall that $\mathbb{Q}\wc{X} = \mathbb{Q}\wc{\Phi} \oplus \mathbb{Q}\Phi^{\perp}$ where $\mathbb{Q}\Phi^{\perp}$ is the $\mathbb{Q}$-subspace spanned by $\Phi^{\perp} = \{y \in \wc{X} \mid \langle x,y\rangle = 0$ for all $x \in \Phi\}$. We thus have a natural projection map $\mathbb{Q}\wc{X} \to \mathbb{Q}\wc{\Phi}$ and we denote by $\wc{X}_{\mathrm{ss}} \subseteq \mathbb{Z}\wc{\Omega}$ the image of $\wc{X}$, which contains $\mathbb{Z}\wc{\Phi}$. If $\bG$ is semisimple, then $\wc{X}_{\mathrm{ss}} = \wc{X}$. The map $\tilde{\imath}$ defines an isomorphism $\Tor_{p'}(\mathbb{Z}\wc{\Omega}/\wc{X}_{\mathrm{ss}}) \to \mathcal{Z}(\bG)$ where $\Tor_{p'}$ denotes the $p'$-torsion subgroup. As $2\wc{\rho}_{\Delta} \in \mathbb{Z}\wc{\Phi}$ we see that $z^2 = 1$. The following characterises exactly when $z$ is the identity.
\end{pa}

\begin{lem}\label{lem:kmodp}
Assume $r \in \mathbb{Z}$ is such that $\imath(r/(q-1)+\mathbb{Z}) = k\pmod{p}$ and let $z = \tilde{\imath}(r\wc{\rho}_{\Delta})$. Then $z = 1$ if and only if at least one of the following holds:
%%%%
\begin{enumerate}
	\item $p=2$,
	\item $\wc{\rho}_{\Delta} \in \wc{X}_{\mathrm{ss}}$,
	\item $\log_p(q)$ is even,
	\item $k \pmod{p}$ is a square in $\mathbb{F}_p$.
\end{enumerate}
\end{lem}

\begin{proof}
Clearly $z=1$ if either (i) or (ii) holds. So assume neither of these two conditions hold. Then $z=1$ if and only if $2 \mid r$. As $k \pmod{p} \in \mathbb{F}_p$, we have $\tilde{\imath}(s/(p-1)) = k \pmod{p}$ for some $1 \leqslant s < p$. If $a = \log_p(q)$ then
%%%%
\begin{equation*}
r = s \cdot \frac{q-1}{p-1} = s(p^{a-1}+\cdots+p+1).
\end{equation*}
%%%%
Hence, as $p \neq 2$ we have $2 \mid r$ if and only if either $a$ is even or $2 \mid s$. This last condition is equivalent to $k\pmod{p} \in \mathbb{F}_p$ being a square in $\mathbb{F}_p$.
\end{proof}

\begin{pa}
For an example, suppose $\bG$ is a simple group of type $\A_{n-1}$. Then $\wc{\Pi}$ is a cyclic group of order $n$ and $\wc{\Pi}_{\bG} := \mathbb{Z}\wc{\Omega}/\wc{X}_{\mathrm{ss}}$ is a cyclic group of order $m \mid n$, which we view additively. The Lang map on $\mathcal{Z}(\bG) \cong \rZ(\bG)$ is identified with the map $\mathscr{L} : \wc{\Pi}_{\bG} \to \wc{\Pi}_{\bG}$ given by $\mathscr{L}(c) = (q-1)c$ if $F$ is split and $\mathscr{L}(c) = -(q+1)c$ if $F$ is twisted. Identifying $\mathscr{L}(\wc{\Pi}_{\bG})$ with $\wc{\Pi}_{\bG}/\Ker(\mathscr{L})$, we see easily that $|\mathscr{L}(\wc{\Pi}_{\bG})| = m/\gcd(m,q\pm 1)$. Hence, if $z \neq 1$, in which case $2 \mid m$ and $p \neq 2$, we have $u$ is $\bG^F$-conjugate to $u^k$ if and only if $|\mathscr{L}(\wc{\Pi}_{\bG})|$ is even.

Given this discussion and the data in Table~\ref{tab:half-sum}, the other simple groups are easily dealt with in this way. This seems to simplify the discussion in \cite[\S5]{tiep-zalesski:2004:unipotent-elements} and the proof of \cite[Thm.~1.7]{tiep-zalesski:2004:unipotent-elements} in the case of good characteristic, especially in the case of twisted groups. We will need the following, whose proof is easy and left to the reader. We simply note that when $\bG = \SO(V)$, we have $\widecheck{\omega}_1 \in \wc{X}_{\mathrm{ss}} = \wc{X}$.
\end{pa}

\begin{lem}\label{cor:reg-class-conj}
Assume $p$ is good for $\bG$ and $u \in \bG^F$ is a regular unipotent element. If $\bG = \GL(V)$ or $\SO(V)$, then $u$ is rational in $\bG^F$. If $\bG = \Sp(V)$, then $u$ and $u^k$ are $\bG^F$-conjugate if and only if $k \pmod{p} \in \mathbb{F}_q$ is a square.
\end{lem}

\section{A Variation of GGGCs for Disconnected Groups}\label{sec:GGGC-ext}
\begin{pa}
In this section, we assume $\bG$ is a (possibly disconnected) reductive algebraic group with Frobenius endomorphism $F : \bG \to \bG$. If $u \in \bG^{\circ F}$ is a unipotent element, then one has a corresponding GGGC $\Gamma_u = \Gamma_u^{\bG^{\circ F}}$ of the finite group $\bG^{\circ F}$. We would like to construct such characters for the group $\bG^F$. One option is the \emph{induced GGGC} $\Ind_{\bG^{\circ F}}^{\bG^F}(\Gamma_u)$. However, we would like to consider a different construction, which instead involves extending $\Gamma_u$ and then inducing to $\bG^F$.

If $r = p^a$ with $a \geqslant 0$, then we denote by $F_r : \mathbb{F} \to \mathbb{F}$ the map defined by $F_r(k) = k^r$. Let $\Frob_1(\bG)$ denote the automorphism group of $\bG$ as an algebraic group. If $r > 1$, then let $\Frob_r(\bG)$ denote the set of Frobenius endomorphisms on $\bG$ endowing $\bG$ with an $\mathbb{F}_r$-structure. We let $\Frob(\bG) = \bigcup_{a \geqslant 0} \Frob_{p^a}(\bG)$. Moreover, $\Frob(\bG,F) \subseteq \Frob(\bG)$ denotes those $\sigma \in \Frob(\bG)$ commuting with $F$, and $\Frob_r(\bG,F) = \Frob(\bG,F) \cap \Frob_r(\bG)$. If $\sigma \in \Frob_r(\bG)$, then we have a corresponding $F_r$-semilinear endomorphism $\sigma : \lie{g} \to \lie{g}$ of the Lie algebra $\lie{g}$ of $\bG^{\circ}$, which we also denote by $\sigma$, see \cite[11.2, 11.3]{taylor:2018:action-of-automorphisms-symplectic}.
\end{pa}

\begin{pa}\label{pa:springer-prop}\label{pa:notation-1}
In what follows, we will freely use the notation and terminology of \cite{brunat-dudas-taylor:2020:unitriangular-shape-of-decomposition-matrices} particularly that of \cite[\S6]{brunat-dudas-taylor:2020:unitriangular-shape-of-decomposition-matrices}. In particular, we assume throughout that $\bG^{\circ}$ is proximate and that $\mathcal{K} = (\phi_{\spr},\kappa,\chi_p)$ is a Kawanaka datum with respect to which our GGGCs will be defined. After \cite[Prop.~11.5]{taylor:2018:action-of-automorphisms-symplectic}, we may, and will, assume that for any $\sigma \in \Frob_r(\bG)$ with $r = p^a \geqslant 1$, the following hold: $\phi_{\spr}\circ \sigma = \sigma\circ\phi_{\spr}$, and $\kappa(\sigma(X),\sigma(Y)) = F_r(\kappa(X,Y))$.

We fix a unipotent element $u \in \bG^{\circ F}$ and a corresponding Dynkin cocharacter $\lambda \in \mathcal{D}_u(\bG^{\circ})^F$, see \cite[\S4.2]{brunat-dudas-taylor:2020:unitriangular-shape-of-decomposition-matrices} for the precise definition of $\mathcal{D}_u(\bG^{\circ}) \subseteq \wc{X}(\bG)$. Let $\bL = \rC_{\bG}(\lambda)$ be the centraliser of the cocharacter. Then $\bL^{\circ} = \rC_{\bG^{\circ}}(\lambda)$ is a Levi subgroup of $\bG^{\circ}$. Recall that for any integer $i < 0$ we have a closed connected $F$-stable unipotent subgroup $\bU(\lambda,i) \leqslant \bG$. This group is defined such that the Lie algebra $\mathrm{Lie}(\bU(\lambda,i)) = \bigoplus_{j \leqslant i} \lie{g}(\lambda,j)$ is a direct sum of weight spaces with respect to $\lambda$, see \cite[\S3]{brunat-dudas-taylor:2020:unitriangular-shape-of-decomposition-matrices} for more details.

If $\mathcal{U}(\bG^{\circ}) \subseteq \bG^{\circ}$ is the closed set of unipotent elements, then we have a function $\eta_u : \mathcal{U}(\bG^{\circ})^F \to \Ql$ whose restriction to $\bU(\lambda,-2)^F$, denoted by $\eta_{u,\lambda}$, is a linear character. This character satisfies ${}^x\eta_{u,\lambda} = \eta_{{}^xu,{}^x\lambda}$ for all $x \in \bG^F$. We have a corresponding character $\zeta_{u,\lambda} \in \Irr(\bU(\lambda,-1)^F)$ defined uniquely by the condition that
%%%%
\begin{equation*}
\Ind_{\bU(\lambda,-2)^F}^{\bU(\lambda,-1)^F}(\eta_{u,\lambda}) = q^{\dim \lie{g}(\lambda,-1)/2}\zeta_{u,\lambda}.
\end{equation*}
%%%%
We then have $\Gamma_u = \Ind_{\bU(\lambda,-1)^F}^{\bG^{\circ F}}(\zeta_{u,\lambda})$.
\end{pa}

\begin{pa}\label{pa:notation-2}
Recall that $\ell \neq p$ is a prime. We denote by $\rC_{\bG}(\lambda)_{\zeta_{u,\lambda}}^F$ the stabiliser of $\zeta_{u,\lambda}$ in $\rC_{\bG}(\lambda)^F$. We will assume that $S_{u,\lambda} \leqslant \rC_{\bG}(\lambda)_{\zeta_{u,\lambda}}^F$ is a fixed Sylow $\ell$-subgroup. The group $\rC_{\bG}(\lambda)$ normalises $\bU(\lambda,-1)$, so $\bU(\lambda,-1)^F\cdot S_{u,\lambda}$ is a group.  In fact, this group is a semidirect product, since $\bU(\lambda,-1)^F$ is a $p$-group. Moreover, $\zeta_{u,\lambda}$ has a unique extension $\widehat{\zeta}_{u,\lambda} \in \Irr(\bU(\lambda,-1)^FS_{u,\lambda})$ such that $o(\widehat{\zeta}_{u,\lambda}) = o(\zeta_{u,\lambda})$ by \cite[Thm.~1.5]{isaacs:2018:characters-of-solvable}. Here, if $G$ is a finite group and $\chi \in \Irr(G)$ is an irreducible character, then $o(\chi)$ denotes the order of the corresponding determinant character. We now define a character
%%%%
\begin{equation*}
\widehat{\Gamma}_{u,\ell} = \Ind_{\bU(\lambda,-1)^FS_{u,\lambda}}^{\bG^F}(\widehat{\zeta}_{u,\lambda}).
\end{equation*}
%%%%
Note that in the special case $\bG = \bG^{\circ}$, we have $\widehat{\Gamma}_{u,\ell}$ is a summand of $\Gamma_u$, which is part of the ideas considered in \cite{brunat-dudas-taylor:2020:unitriangular-shape-of-decomposition-matrices}. More generally, we have the following.
\end{pa}

\begin{lem}\label{lem:restrict-gamhat}
For any unipotent element $u \in \bG^{\circ F}$ we have
%%%%
\begin{equation*}
\Res_{\bG^{\circ F}}^{\bG^F}(\widehat{\Gamma}_{u,\ell}) = \sum_{g \in \bG^F/\bG^{\circ F}S_{u,\lambda}} \Gamma_{gug^{-1}}.
\end{equation*}
%%%%
Moreover, if $\rC_{\bG}(u)^F_{\ell} \leqslant \rC_{\bG}(u)^F$ is a Sylow $\ell$-subgroup, then $\rC_{\bG}(u)^F_{\ell} \leqslant \bG^{\circ F}S_{u,\lambda}$.
\end{lem}

\begin{proof}
The first statement follows from the usual Mackey formula for finite groups, together with \cite[Prop.~11.10]{taylor:2018:action-of-automorphisms-symplectic}. If $x \in \rC_{\bG}(u)^F$ then ${}^x\lambda \in \mathcal{D}_u(\bG^{\circ})^F$. By \cite[Lem.~3.6]{brunat-dudas-taylor:2020:unitriangular-shape-of-decomposition-matrices}, and a standard application of the Lang--Steinberg Theorem, there exists $h \in \rC_{\bG^{\circ}}^{\circ}(u)^F$ such that $\lambda = {}^{hx}\lambda$. Therefore, if $\bL = \rC_{\bG}(\lambda)$ then $g = hx \in \rC_{\bL}(u)^F \leqslant \bL_{\zeta_{u,\lambda}}^F$ which implies $\rC_{\bG}(u)^F \leqslant \rC_{\bG^{\circ}}^{\circ}(u)^F\cdot \bL_{\zeta_{u,\lambda}}^F$. Using conjugacy of Sylow subgroups in the quotient $\rC_{\bG}(u)^F/\rC_{\bG^{\circ}}^{\circ}(u)^F$ gives the second statement.
\end{proof}

\begin{rem}
We will be interested in applying this construction in the special case where $\bG = \GO(V)$ is an orthogonal group and $p \neq 2$ is odd. In this case, $\bG/\bG^{\circ}$ has order $2$. If the centraliser $\rC_{\bG}(u)$ of $u \in \bG^{\circ F}$ is not contained in $\bG^{\circ}$ then $\bG^F = \bG^{\circ F}\rC_{\bG}(u)_2^F$ and Lemma~\ref{lem:restrict-gamhat} shows that $\widehat{\Gamma}_{u,2}$ is an extension of the GGGC $\Gamma_u$.
\end{rem}

The following extends \cite[Prop.~11.10]{taylor:2018:action-of-automorphisms-symplectic} to give the action of the automorphism group on these characters.

\begin{thm}\label{prop:GGGCext-auts}
Recall our assumption that $\bG^{\circ}$ is proximate, $p$ is good for $\bG$, and $\ell \neq p$ is a prime. For any unipotent element $u \in \bG^{\circ F}$ and $\sigma \in \Frob(\bG,F)$, we have ${}^{\sigma}\widehat{\Gamma}_{u,\ell} = \widehat{\Gamma}_{\sigma(u),\ell}$.
\end{thm}

\begin{proof}
Recall that the restriction of $\sigma$ to $\bG^F$ is an automorphism. It is easy to see that $\sigma(\bU(\lambda,-1)) = \bU(\sigma\cdot\lambda,-1)$, where $\sigma\cdot\lambda$ is as in \cite[11.6]{taylor:2018:action-of-automorphisms-symplectic}, and we get that ${}^{\sigma}\zeta_{u,\lambda} = \zeta_{\sigma(u),\sigma\cdot\lambda}$ by arguing as in \cite[Prop.~11.10]{taylor:2018:action-of-automorphisms-symplectic}. Moreover, $\sigma(\rC_{\bG}(\lambda)_{\zeta_{u,\lambda}}^F) = \rC_{\bG}(\sigma\cdot \lambda)_{\zeta_{\sigma(u),\sigma\cdot\lambda}}^F$ and there exists an element $g \in \rC_{\bG}(\sigma\cdot \lambda)_{\zeta_{\sigma(u),\sigma\cdot\lambda}}^F$ such that ${}^g\sigma(S_{u,\lambda}) = S_{\sigma(u),\sigma\cdot\lambda}$. Now certainly
%%%%
\begin{equation*}
{}^{\sigma}\widehat{\Gamma}_{u,\ell} = {}^{g\sigma}\widehat{\Gamma}_{u,\ell} = \Ind_{\bU(\sigma\cdot \lambda,-1)^F S_{\sigma(u),\sigma\cdot\lambda}}^{\bG^F}({}^{g\sigma}\widehat{\zeta}_{u,\lambda})
\end{equation*}
%%%%
and ${}^{g\sigma}\widehat{\zeta}_{u,\lambda}$ is an extension of $\zeta_{\sigma(u),\sigma\cdot\lambda}$. By the unicity of the extension, ${}^{g\sigma}\widehat{\zeta}_{u,\lambda} = \widehat{\zeta}_{\sigma(u),\sigma\cdot\lambda}$ and the result follows.
\end{proof}

We now want to extend \cite[Prop.~4.10]{schaeffer-fry-taylor:2018:on-self-normalising-Sylow-2-subgroups}, which gives the action of $\gal$ on GGGCs. Understanding the effect of $\gal$ on $\widehat{\Gamma}_u$ seems to be more complicated and appears to be related to stronger notions of conjugacy, such as being strongly real. The following, though weaker than the result in \cite{schaeffer-fry-taylor:2018:on-self-normalising-Sylow-2-subgroups},  will be sufficient for our purposes.

\begin{prop}\label{prop:GGGCext-gal}
Recall our assumption that $\bG^{\circ}$ is proximate, $p$ is good for $\bG$, and $\ell \neq p$ is a prime. Assume $\sigma \in \gal$ is a Galois automorphism and $k \in \mathbb{Z}$ is an integer coprime to $p$ such that $\xi^{\sigma} = \xi^k$ for all $p$th roots of unity $\xi \in \Ql^{\times}$. If $u \in \bG^{\circ F}$ is a unipotent element and $u^k$ is $\bG^{\circ F}$-conjugate to $u$, then $\widehat{\Gamma}_{u,\ell}^{\sigma} = \widehat{\Gamma}_{u,\ell}$.
\end{prop}

\begin{proof}
Let $\bL = \rC_{\bG}(\lambda)$. From the proof of \cite[Prop.~4.10]{schaeffer-fry-taylor:2018:on-self-normalising-Sylow-2-subgroups}, there exists an element $x \in \bL^{\circ}$ such that $\phi_{\spr}({}^xu) = k\phi_{\spr}(u) \in \lie{g}(\lambda,2)_{\reg}$ and ${}^xu$ is $\bG^{\circ F}$-conjugate to $u^k$. By assumption, ${}^{gx}u = u$ for some $g \in \bG^{\circ F}$. The natural map $\rC_{\bL^{\circ}}(u)/\rC_{\bL^{\circ}}^{\circ}(u) \to \rC_{\bG^{\circ}}(u)/\rC_{\bG^{\circ}}(u)$ is an isomorphism. Hence, a quick application of the Lang--Steinberg Theorem in the group $\rC_{\bL^{\circ}}^{\circ}(u)$ implies that, after possibly replacing $x$ by $xa$ with $a \in \rC_{\bL^{\circ}}(u)$, we may assume that $x \in \bL^{\circ F}$.

Arguing as in \cite[Prop.~4.10]{schaeffer-fry-taylor:2018:on-self-normalising-Sylow-2-subgroups}, we see that $\eta_u^{\sigma} = \eta_{{}^xu} = {}^x\eta_u$. It follows that $\zeta_{u,\lambda}^{\sigma} = \zeta_{{}^xu,\lambda} = {}^x\zeta_{u,\lambda}$. Being a Galois conjugate, the character $\zeta_{u,\lambda}^{\sigma} = {}^x\zeta_{u,\lambda}$ has the same stabiliser as $\zeta_{u,\lambda}$ in $\rC_{\bG}(\lambda)^F$. This implies $x$ normalises the stabiliser of $\zeta_{u,\lambda}$ in $\rC_{\bG}(\lambda)^F$, so $S_{u,\lambda}^{xg} = S_{u,\lambda}$ for some $g \in \rC_{\bG}(\lambda)_{\zeta_{u,\lambda}}^F$. By the unicity of the extension, we have $\widehat{\zeta}_{u,\lambda}^{\sigma xg} = \widehat{\zeta}_{u,\lambda}$ so
%%%%
\begin{equation*}
\widehat{\Gamma}_{u,\ell}^{\sigma} = \widehat{\Gamma}_{u,\ell}^{\sigma xg} = \Ind_{\bU(\lambda,-1)^F\cdot S_{u,\lambda}^{xg}}^{\bG^F}(\widehat{\zeta}_{u,\lambda}^{\sigma xg}) = \widehat{\Gamma}_{u,\ell}.
\end{equation*}
%%%%
\end{proof}

\section{Classical Groups}\label{sec:classical-groups}
\begin{assumption}
From now on, we assume that $p \neq 2$.
\end{assumption}

\begin{pa}\label{pa:bilinear-forms}
We fix a finite dimensional $\mathbb{F}$-vector space $V$, with $\dim(V)>1$, and let $\mathcal{B} : V \times V \to \mathbb{F}$ be a non-degenerate symmetric or alternating bilinear form. We denote by $\GL(V\mid\mathcal{B})$, resp., $\SL(V\mid\mathcal{B})$, the subgroup of isometries of $\mathcal{B}$ in the general linear group, resp., special linear group. If $\mathcal{B}$ is symmetric, then $\GO(V) := \GL(V\mid \mathcal{B})$ is the orthogonal group and $\SO(V) := \SL(V \mid \mathcal{B}) = \GO(V)^{\circ}$ is the special orthogonal group. If $\mathcal{B}$ is alternating, then $\Sp(V) := \GL(V\mid \mathcal{B}) = \SL(V \mid \mathcal{B})$ is the symplectic group.

For calculations, it will be convenient to make concrete choices for the form $\mathcal{B}$. We define a totally ordered set
%%%%
\begin{equation}\label{eq:basis}
\mathcal{I} = \begin{cases}
\{1 \prec \cdots \prec n \prec 0 \prec -n \prec \cdots \prec-1\} &\text{if $\dim(V)=2n+1$ is odd}\\
\{1 \prec \cdots \prec n \prec -n \prec \cdots \prec -1\} &\text{if $\dim(V)=2n$ is even}
\end{cases}
\end{equation}
%%%%
and choose a corresponding basis $\mathcal{V} = (v_i \mid i \in \mathcal{I}) \subset V$ ordered by $\prec$. For any integer $m \in \mathbb{Z}$, let $\sgn(m) = 1$ if $m \geqslant 0$ and $\sgn(m) = -1$ if $m < 0$. Let $\epsilon \in \{0,1\}$ be such that $\mathcal{B}(v,w) = (-1)^{\epsilon}\mathcal{B}(w,v)$ for all $v,w \in V$. Then we assume that $\mathcal{B}(v_i,v_j) = \sgn(i)^{\epsilon}\cdot\delta_{i,-j}$ for all $i,j \in \mathcal{I}$
\end{pa}

\subsection*{Tori and Normalisers}
\begin{pa}\label{pa:notnT}\label{pa:Wgens}
Let $\bG = \GL(V \mid \mathcal{B})$. Via our choice of basis $\mathcal{V}$, we will interchangeably describe elements of $\GL(V)$ either as linear maps or matrices. We will denote by $\bT \leqslant \bB \leqslant \bG^{\circ}$ the subgroups of diagonal and upper triangular matrices. For $i \in \mathcal{I}$ we let $\wc{\varepsilon}_i : \mathbb{F}^{\times} \to \GL(V)$ be the homomorphism defined such that $\wc{\varepsilon}_i(\zeta)v_j = \zeta^{\delta_{i,j}}v_j$ for any $\zeta \in \mathbb{F}^{\times}$ and $j \in \mathcal{I}$. We have $\bT = \bT_1\times\cdots\times \bT_n$, where $\bT_i$ is the image of $\mathbf{d}_i : \mathbb{F}^{\times} \to \bT$ defined by setting $\mathbf{d}_i(\zeta) = \wc{\varepsilon}_i(\zeta)\wc{\varepsilon}_{-i}(\zeta^{-1})$.

Denote by $\mathfrak{S}_{\mathcal{I}}$ the symmetric group on $\mathcal{I}$. Let $W_n \leqslant \mathfrak{S}_{\mathcal{I}}$ be the subgroup generated by $\mathbb{S} = \{s_1,\dots,s_{n-1},s_n\}$, where $s_i = (i,i+1)(-i,-i-1)$ for $1\leqslant i\leqslant n-1$ and $s_n = (n,-n)$. The pair $(W_n,\mathbb{S})$ is a Coxeter system of type $\B_n$. We define elements:
%%%%
\begin{itemize}
	\item $t_m = (m,-m) = s_m\cdots s_{n-1}s_n s_{n-1}\cdots s_m$ if $1 \leqslant m \leqslant n$,
	\item $u_m = (m,-m)(n,-n) = s_m\cdots s_{n-2}(s_ns_{n-1}s_n s_{n-1})s_{n-2}\cdots s_m$ if $1 \leqslant m < n$.
\end{itemize}
%%%%
The subgroup $W_n' \leqslant W_n$ generated by $\mathbb{S}' = \{s_1,\dots,s_{n-1},s_ns_{n-1}s_n\}$ is a Coxeter group of type $\D_n$. We have a further subgroup $W_n'' \leqslant W_n'$, of type $\B_{n-1}$, generated by $\{s_1,\dots,s_{n-2},u_{n-1}\}$. When $n = 1$, we have $W_n'' = W_n'$ is trivial. Let $\wt{\mathbb{S}} = \mathbb{S} \cup \{s_0\}$, where $s_0 = t_1$. Assume $0 \leqslant a \leqslant n$ and let $b = n-a$. Then we have a reflection subgroup  $W_a\times W_b \leqslant W_n$ of type $\B_a\times\B_b$ with Coxeter generators $\wt{\mathbb{S}}\setminus \{s_a\}$.
\end{pa}

\begin{pa}\label{pa:lifts-reflections}
Recall that for any $\sigma \in \mathfrak{S}_{\mathcal{I}}$, we have a corresponding transformation $p_{\sigma} \in \GL(V)$ such that $p_{\sigma}(v_i) = v_{\sigma(i)}$ for all $i \in \mathcal{I}$. The groups $\bW = \rN_{\bG}(\bT)/\bT$ and $\bW^{\circ} = \rN_{\bG^{\circ}}(\bT)/\bT$ denote the Weyl groups of $\bG$ and $\bG^{\circ}$. For each $s \in \mathbb{S}$, we define an element $n_{s} \in \rN_{\bG}(\bT)$ as follows. We let $n_{s_i} := p_{s_i}$ for any $1 \leqslant i < n$ and
%%%%
\begin{equation*}
n_{s_n} := \begin{cases}
p_{s_n}\wc{\varepsilon}_n(-1) &\text{if }\bG = \Sp(V),\\
p_{s_n}\wc{\varepsilon}_0(-1) &\text{if }\bG = \GO(V)\text{ and }\dim(V)\text{ is odd}\\
p_{s_n} &\text{if }\bG = \GO(V)\text{ and }\dim(V)\text{ is even}.
\end{cases}
\end{equation*}
%%%%
The map defined by $n_s\bT \mapsto s$, for $s \in \mathbb{S}$, extends to a unique isomorphism $\bW \to W_n$. We will implicitly identify $\bW$ and $W_n$ in this way.

By Matsumoto's Theorem \cite[Thm.~1.2.2]{geck-pfeiffer:2000:characters-of-finite-coxeter-groups}, we obtain a unique element $n_w = n_{s_1}\cdots n_{s_k} \in \rN_{\bG}(\bT)$, where $w = s_1\cdots s_k$, with $s_i \in \mathbb{S}$, is a reduced expression for $w$. For $\bG = \GO(V)$, the map $w \mapsto n_w$ is a homomorphism $\bW \to \rN_{\bG}(\bT)$, i.e., the Weyl group lifts. When $\dim(V)$ is even and $\bG = \GO(V)$, the group $\bW^{\circ}$ is mapped onto $W_n'$.
\end{pa}

\subsection*{Subspace Subgroups and Conformal Groups}
\begin{pa}\label{pa:subspace-subgrp}
Assume $0 \leqslant m \leqslant n$ is an integer and let $\mathcal{I}_1 = \{1,\dots,m,-m,\dots,-1\}$ and $\mathcal{I}_0 = \mathcal{I}\setminus\mathcal{I}_1$. We have corresponding subspaces $V_i \subseteq V$ spanned by $\{v_k \mid k\in \mathcal{I}_i\}$ and tori $\bS_i = \prod_{k \in \mathcal{I}_i\setminus\{0\}} \bT_k$ with $i \in \{0,1\}$. The decomposition $V = V_1 \oplus V_0$ is an orthogonal decomposition of $V$ and we have a corresponding subgroup $\bM_1 \times \bM_0 \leqslant \bG$, where $\bM_i = \GL(V_i \mid \mathcal{B})$. Note that the group $\bL_J = \bS_1\times\bM_0^{\circ}$ is a standard Levi subgroup of $\bG^{\circ}$ with $J \subseteq \Delta$. Identifying $\bW_i = \rN_{\bM_i}(\bS_i)/\bS_i$ with a subgroup of $\bW = W_n$, we have $\bW_1\times\bW_0 = W_m\times W_{n-m}$ is the subgroup of type $\B_m\times\B_{n-m}$ in \ref{pa:Wgens}. We recall the following well-known result.
\end{pa}

\begin{lem}\label{lem:normaliser}
Assume $\bL_J = \bS_1\times\bM_0^{\circ}$ is a standard Levi subgroup as in \ref{pa:subspace-subgrp} such that $1 \leqslant m \leqslant n-2$. Then we have $\rN_{\bW^{\circ}}(J) \cong W_m$ is a Weyl group of type $\B_m$ with generators
%%%%
\begin{equation*}
\begin{cases}
\{s_1,\dots,s_{m-1},u_m\} &\text{when $\bG = \GO(V)$ and $\dim(V)$ is even},\\
\{s_1,\dots,s_{m-1},t_m\} &\text{otherwise}.
\end{cases}
\end{equation*}
\end{lem}

\begin{pa}\label{pa:conformal}
Now assume $\dim(V)$ is even. In this case, we let $\wt{\bG} = \bG\cdot \rZ(\GL(V))$ be the corresponding conformal group and similarly $\wt{\bG}^{\circ} = \bG^{\circ}\cdot \rZ(\GL(V))$. The inclusion $\bG^{\circ} \to \wt{\bG}^{\circ}$ is a regular embedding, as in \ref{pa:reg-embedding}, and the group $\wt{\bT} = \bT\cdot \rZ(\GL(V))$ is a maximal torus of $\wt{\bG}^{\circ}$. We define a closed embedding $\mathbf{z} : \mathbb{F}^{\times} \to \wt{\bT}$ by setting $\mathbf{z}(\zeta) = \wc{\varepsilon}_1(\zeta)\cdots \wc{\varepsilon}_n(\zeta)$. Writing $\bZ$ for the image of $\mathbf{z}$,  we have $\wt{\bT} = \bT \times \bZ$. Moreover, we have $\wt{\bM}_0 = \bM_0\cdot \bZ$ and $\wt{\bM}_0^{\circ} = \bM_0^{\circ}\cdot \bZ$ are also isomorphic to conformal groups, and $\wt{\bL}_J = \bS_1\times \wt{\bM}_0^{\circ}$. We note that for any $1 \leqslant k \leqslant n$, we have
%%%%
\begin{equation*}
{}^{t_k}\mathbf{z}(\zeta) = \mathbf{d}_k(\zeta^{-1})\mathbf{z}(\zeta) \quad\text{and}\quad {}^{u_k}\mathbf{z}(\zeta) = \mathbf{d}_k(\zeta^{-1})\mathbf{d}_n(\zeta^{-1})\mathbf{z}(\zeta).
\end{equation*}
%%%%
Furthermore, for $\prod_{i=1}^n \mathbf{d}_i(\zeta_i)\in\bT$ and $1 \leqslant k \leqslant n$, we have
%%%%
\begin{equation*}
{}^{t_k}\left(\left(\prod_{i=1}^n \mathbf{d}_i(\zeta_i)\right)\mathbf{z}(\zeta)\right)=\mathbf{d}_k(\zeta_k^{-1}\zeta^{-1})\left(\prod_{i\neq k} \mathbf{d}_i(\zeta_i)\right)\mathbf{z}(\zeta).
\end{equation*}
\end{pa}

\subsection*{Automorphisms and Frobenius Endomorphisms}
\begin{pa}
We want to define an analogue for algebraic groups of the group $\rK(G) \leqslant \Aut(G)$ defined in \ref{pa:aso}. For this, we first consider the central automorphisms of \ref{pa:central-auts} in the context of algebraic groups.

Let $\bH$ be an affine algebraic group. If $\alpha \in \Hom(\bH,\rZ(\bH))$ is a homomorphism of algebraic groups then $\tau_{\alpha} : \bH \to \bH$ defined by $\tau_{\alpha}(g) = g\cdot\alpha(g)$ is also a homomorphism of algebraic groups. Hence $\rC_{\Aso(\bH)}(\Inn(\bH))$ consists of those $\tau_{\alpha}$ that are invertible.

Denote by $\bH_{\mathrm{u}} \subseteq \bH$ the set of unipotent elements and define $\rO^{p'}(\bH) = \langle (\bH^{\circ})_{\mathrm{u}}\rangle \leqslant \bH^{\circ}$ to be the subgroup generated by the unipotent elements of $\bH^{\circ}$. We then define a submonoid
%%%%
\begin{equation*}
\rK(\bH) = \{\varphi \in \rC_{\Aso(\bH)}(\Inn(\bH)) \mid \varphi|_{\rO^{p'}(\bH)} = \mathrm{Id}_{\rO^{p'}(\bH)}\} \leqslant \Aso(\bH).
\end{equation*}
\end{pa}

\begin{rem}\label{rem:op-prime-red-grp}
The group $\rO^{p'}(\bH)$ is an $\Aso(\bH)$-stable closed connected subgroup of $\bH$. If $\bH^{\circ}$ is reductive then $\rO^{p'}(\bH) = [\bH^{\circ},\bH^{\circ}]$ is the commutator subgroup of $\bH^{\circ}$.
\end{rem}

\begin{pa}\label{pa:ker-res-aso}
We now wish to describe the monoids $\Aso(\bG)$ and $\Aso(\bG^{\circ})$ where $\bG = \GL(V \mid \mathcal{B})$ as above. The description of $\Aso(\bG^{\circ})$ is well known and can be found in \cite[Thm.~1.15.7]{gorenstein-lyons-solomon:1998:classification-3}. We will describe a precise generating set for $\Aso(\bG^{\circ})$ below.

Here we focus on $\Aso(\bG)$ by considering the natural restriction map $\Aut_1(\bG) \to \Aut_1(\bG^{\circ})$. Suppose $\phi \in \Aut_1(\bG)$ is in the kernel of this restriction map. Following the proof of Proposition~\ref{prop:aut-grp-quasisimp} we see that ${}^gx =\phi({}^gx) = {}^{\phi(g)}x$ for all $x \in \bG^{\circ}$ and $g \in \bG$ because $\bG^{\circ}$ is normal in $\bG$. Hence $g^{-1}\phi(g) \in \rC_{\bG}(\bG^{\circ})$. The following shows that $\phi$ must be a central automorphism.
\end{pa}

\begin{lem}\label{lem:centraliser-orthogonal}
We have $\rC_{\bG}(\bG^{\circ}) = \rZ(\bG)$ and $\rZ(\bG) \leqslant \bG^{\circ}$ if and only if $\dim(V)$ is even.
\end{lem}

\begin{proof}
It is well known that $\rZ(\bG) = \{\pm\ID_V\}$ which implies the second statement and this yields the first statement when $\dim(V)$ is odd because $\bG \cong \bG^{\circ} \times \rZ(\bG)$ in this case. Hence we can assume $\dim(V)$ is even and $\bG = \GO(V)$.

As is well known, the natural homomorphism $\bG/\bG^{\circ} \to \Out(\bG^{\circ})$ to the outer automorphism group, given by $g\bG^{\circ} \mapsto \iota_g|_{\bG^{\circ}}\Inn(\bG^{\circ})$, is injective. This implies that $\rC_{\bG}(\bG^{\circ}) = \rC_{\bG^{\circ}}(\bG^{\circ}) = \rZ(\bG^{\circ}) = \rZ(\bG)$.
\end{proof}

\begin{rem}
This implies the natural map $\Inn(\bG) \to \Aut(\bG^{\circ})$ given by $\iota_g \mapsto \iota_g|_{\bG^{\circ}}$ is injective. Hence, the normal subgroup $\{\iota_g \mid g \in \bG^{\circ}\} \leqslant \Inn(\bG)$ is isomorphic to $\Inn(\bG^{\circ})$ so we will again denote this, unambiguously, by $\Inn(\bG^{\circ})$.
\end{rem}

\begin{pa}\label{pa:auts-ortho}\label{pa:gamma-auts}
If $q = p^a$ is a power of $p$ with $a \geqslant 1$ an integer, then we denote by $F_q : V \to V$ the unique semi-linear map fixing the basis $\mathcal{V}$ and satisfying $F_q(av) = a^qv$ for all $a \in \mathbb{F}$. Precomposing with $F_q$ gives a Frobenius endomorphism $F_q : \GL(V) \to \GL(V)$. We also denote by $F_q$ the restriction of this map to any $F_q$-stable subgroup.

Assume $\bG = \GO(V)$ and $\dim(V)$ is even. We define $\gamma_1 = \iota_{n_w} \in \Inn(\bG)$ to be the inner automorphism induced by the representative $n_w \in \rN_{\bG}(\bT)$ of $w = s_n \in \bW$ chosen in \ref{pa:lifts-reflections}. By Lemma~\ref{lem:centraliser-orthogonal} we have $\Inn(\bG) = \Inn(\bG^{\circ}) \rtimes \langle \gamma_1\rangle$. We will also denote the restriction $\gamma_1|_{\bG^{\circ}} \in \Aso(\bG^{\circ})$ by $\gamma_1$.

Recall that $[\bG,\bG] = [\bG^{\circ},\bG^{\circ}] = \bG^{\circ}$. We define $\gamma_2 = \tau_{\alpha} \in \Aso(\bG)$, where $\alpha : \bG \to \rZ(\bG)$ is the unique non-trivial homomorphism, so that $\gamma_2(g) = g\cdot\alpha(g)$ for all $g \in \bG$. As $\dim(V)$ is even we have $\rZ(\bG) \leqslant \bG^{\circ} = \Ker(\alpha)$ so $\rK(\bG) = \langle \gamma_2\rangle_{\mon} \cong C_2$ is a group. 
To have a consistent notation we let $\gamma_1 = \gamma_2$ be the identity if either $\bG = \Sp(V)$ or $\dim(V)$ is odd. In these cases we have $\rK(\bG) = \langle \gamma_2\rangle$ is trivial. Now let $\Gamma(\bG) = \langle F_p\rangle_{\mon} \times \langle \gamma_1\rangle \leqslant \Aso(\bG)$ be the submonoid generated be these elements. The image of the map $\Gamma(\bG) \to \Aso(\bG^{\circ})$ given by restriction will be denoted by $\Gamma(\bG^{\circ})$.
\end{pa}

\begin{lem}\label{lem:asogeny-orthogonal}
If $\bG = \GL(V \mid \mathcal{B})$ is a classical group then the following hold:
%%%%
\begin{enumerate}
	\item $\Aso(\bG^{\circ}) = \Inn(\bG^{\circ}) \rtimes \Gamma(\bG^{\circ})$,
	\item $\Aso(\bG) = K(\bG) \times (\Inn(\bG^{\circ}) \rtimes \Gamma(\bG))$.
\end{enumerate}
\end{lem}

\begin{proof}
(i). By \cite[Thm.~1.15.7]{gorenstein-lyons-solomon:1998:classification-3} we have $\Aut_1(\bG^{\circ}) = \Inn(\bG^{\circ})\rtimes (\langle F_p\rangle_{\grp} \times \langle \gamma_1 \rangle_{\grp})$ and from this the statement is clear.

(ii). The restriction map $\Aut_1(\bG) \to \Aut_1(\bG^{\circ})$ is visibly surjective and from the discussion in \ref{pa:ker-res-aso} we see that the kernel consists of central automorphisms. Hence $\Aut_1(\bG) = K(\bG)(\Inn(\bG)\langle F_p\rangle_{\grp})$ and $K(\bG) \cap \Inn(\bG)\langle F_p\rangle_{\grp}$ is trivial because the restriction map gives an isomorphism $\Inn(\bG)\langle F_p\rangle_{\grp} \to \Aut_1(\bG^{\circ})$ as is easily seen. The statement is an easy consequence of this.
\end{proof}

\begin{rem}
Here the semidirect products are semidirect products of monoids. Those not familiar with this setup can simply replace $\Aso(\bG)$ with $\Aut_1(\bG)$ and $\Gamma(\bG)$ by the subgroup generated by $\Gamma(\bG)$ to have a semidirect product of groups.
\end{rem}

\begin{rem}
Note that we have
%%%%
\begin{equation*}
\lvert \rK(\bG) \rvert = \lvert \Hom(\bG/\bG^{\circ}\rZ(\bG),\rZ(\bG))\rvert = \lvert \bG/\bG^{\circ}\rZ(\bG) \rvert
\end{equation*}
%%%%
and this can be calculated from Lemma~\ref{lem:centraliser-orthogonal}.
\end{rem}

\begin{pa}
For $q = p^a$ with $a \geqslant 1$, we let $\Gamma_q(\bG) = F_q\langle \gamma_1\rangle \subseteq \Gamma(\bG)$ and $\Gamma_q(\bG^{\circ}) = F_q\langle \gamma_1\rangle \subseteq \Gamma(\bG^{\circ})$ be the coset defined by $F_q$. We say $F \in \Gamma_q(\bG)$ is \emph{split} if $F = F_q$ and \emph{twisted} if $F = F_q\gamma_1$.

It follows from Lemma~\ref{lem:asogeny-orthogonal} that, up to conjugation by elements of $\bG^{\circ}$, any Frobenius endomorphism endowing $\bG$ with an $\mathbb{F}_q$-rational structure is contained in $F_q\Gamma(\bG)\rK(\bG)$. To our knowledge, the Frobenius endomorphisms $F\gamma_2$ of $\bG$, with $F \in \Gamma_q(\bG)$, are not particularly well studied or discussed in the literature.

It seems that the groups $\bG^F$ and $\bG^{F\gamma_2}$ can be subtly different. For instance, whilst $\bG^F$ is a split extension of $\bG^{\circ F} = \bG^{\circ F\gamma_2}$ this does not seem to always be the case for $\bG^{F\gamma_2}$. Discussing this would lead us too far away from our main goal, so we will abstain from further discussion of the endomorphisms in $\Gamma_q(\bG)\gamma_2$. For now we need the following.
\end{pa}

\begin{lem}\label{lem:cent-fin-con-comp}
If $F \in \Gamma_q(\bG)$ then $\rC_{\bG}(\rO^{p'}(\bG^F)) = \rZ(\bG)$.
\end{lem}

\begin{proof}
It suffices to show that $\rC_{\bG}(\rO^{p'}(\bG^F)) \leqslant \rC_{\bG^{\circ}}(\rO^{p'}(\bG^F))\cdot \rZ(\bG)$ because, as in the proof of Proposition~\ref{prop:aut-grp-quasisimp}, the argument in \cite[Lem.~6.1]{bonnafe:2006:sln} shows that $\rC_{\bG^{\circ}}(\rO^{p'}(\bG^F)) = \rZ(\bG^{\circ}) = \rZ(\bG)$. Thus we can assume $\bG = \GO(V)$. If $\dim(V)$ is odd then this is clear as $\bG = \bG^{\circ}\rZ(\bG)$.

So suppose $\dim(V)$ is even. As in the proof of Lemma~\ref{lem:centraliser-orthogonal} it is enough to show that $\gamma_1|_{\rO^{p'}(\bG^F)}$ is not inner. But this follows from the explicit description of $\Out(\rO^{p'}(\bG^F))$ given in \cite[Thm.~2.5.12]{gorenstein-lyons-solomon:1998:classification-3}, see also the detailed discussion in \cite[\S2.7]{gorenstein-lyons-solomon:1998:classification-3}. Note that when $F = F_q\gamma_1$ is twisted then $\gamma_1|_{\rO^{p'}(\bG^F)} = F_q|_{\rO^{p'}(\bG^F)}$.
\end{proof}

\begin{rem}
An alternative argument is to show that $\bG^{\gamma_1} \cap \bG^F$ does not contain all $p$-elements for which one can use a count due to Steinberg \cite[Cor.~7.4.6]{digne-michel:2020:representations-of-finite-groups-of-lie-type}.
\end{rem}

\begin{pa}
Let us recall our notation $\Aso(\bG,F)$ and $\Inn(\bG,F)$ introduced in \ref{pa:aso}. The natural map $\Inn(\bG, F) \to \Aut(\bG^{\circ F})$, given by restriction, is injective because by Lemma~\ref{lem:cent-fin-con-comp} $\rC_{\bG}(\bG^{\circ F}) \leqslant \rC_{\bG}(\rO^{p'}(\bG^F)) = \rZ(\bG^F)$ by \cite[Lem.~6.1]{bonnafe:2006:sln}. So we can identify the image of $\rC_{\Inn(\bG^{\circ})}(F)$ in $\Aut(\bG^F)$ with $\Inn(\bG^{\circ},F)$.

Now let $F \in \Gamma_q(\bG)$ and denote by $\Gamma(\bG,F)$ the image of $\rC_{\Gamma(\bG)}(F) = \Gamma(\bG)$ in $\Aut(\bG^{F})$, similarly by $\Gamma(\bG^{\circ},F)$ the image of $\rC_{\Gamma(\bG^{\circ})}(F) = \Gamma(\bG^{\circ})$ in $\Aut(\bG^{\circ F})$ under the natural restriction maps. We also denote by $\rK(\bG,F) \leqslant \Aut(\bG^F)$ the image of $\rC_{\rK(\bG)}(F) = \rK(\bG)$ under this restriction map.
\end{pa}

\begin{lem}\label{lem:asogeny-orthogonal-fin}
If $F \in \Gamma_q(\bG)$ then the following hold:
%%%%
\begin{enumerate}
	\item $\Aso(\bG^{\circ},F) = \Inn(\bG^{\circ},F) \rtimes \Gamma(\bG^{\circ},F)$,
	\item $\Aso(\bG,F) = \rK(\bG,F) \times (\Inn(\bG^{\circ},F) \rtimes \Gamma(\bG,F))$,
	\item $\Aut(\bG^F)_{\bG^{\circ F}} = \rK(\bG^F) \times (\Inn(\bG^{\circ},F) \rtimes \Gamma(\bG,F))$,
\end{enumerate}
%%%%
where $\Aut(\bG^F)_{\bG^{\circ F}}$ is the stabiliser of $\bG^{\circ F}$.
\end{lem}

\begin{proof}
(i) and (ii) are simple consequences of Lemma~\ref{lem:asogeny-orthogonal}. We prove (iii) by considering the natural restriction map $\Aut(\bG^F)_{\bG^{\circ F}} \to \Aut(\bG^{\circ F})$. By Lemma~\ref{lem:cent-fin-con-comp} we see that the kernel of this map is given by central automorphisms. Hence, the kernel is contained in $\rK(\bG^F)$. After (i), (ii), and Proposition~\ref{prop:aut-grp-quasisimp}, to show that the map is surjective it suffices to show that the image of $\rK(\bG^F)$ is $\rK(\bG^{\circ F})$.

If $\bG = \bG^{\circ} = \Sp(V)$ then $\bG^F$ is perfect so $\rK(\bG^F) = \rK(\bG^{\circ F})$ is trivial. Hence, we can assume that $\bG = \GO(V)$. To obtain the statement we follow the discussion in Remark~\ref{rem:cent-in-der}. Note that $\rE^2(\bG^F) = \rO^{p'}(\bG^F) = [\bG^F,\bG^F]$ because $V = \bG^{\circ F}/\rO^{p'}(\bG^F)$ has a complement $U \leqslant \bG^F/\rO^{p'}(\bG^F)$. Indeed, if $\dim(V)$ is even then we may take $U = \langle n_w\rO^{p'}(\bG^F)\rangle/\rO^{p'}(\bG^F)$, where $n_w \in \rN_{\bG^F}(\bT)$ is the element of order $2$ defined in \ref{pa:auts-ortho}, and if $\dim(V)$ is odd then we may take $U = \rZ(\bG^F)\rO^{p'}(\bG^F)/\rO^{p'}(\bG^F)$.

If $\alpha \in \Hom(V,\rZ(\bG^F))$ then inflating we obtain an element of $\Hom(\bG^F,\rZ(\bG^F))$. From this it follows that $\rK(\bG^F)$ maps onto $\rK(\bG^{\circ F})$ via the restriction map.

The above argument shows that the restriction map $\Aut(\bG^F)_{\bG^{\circ F}} \to \Aut(\bG^{\circ F})$ has a kernel complemented in $\rK(\bG^F)$. From this we easily obtain a complement to the kernel in $\Aut(\bG^F)_{\bG^{\circ F}}$ and (iii) follows.
\end{proof}

\begin{rem}
With some effort, by studying the distribution of involutions across the cosets of $\rO^{p'}(\bG^F)$ in $\bG^F$, it seems possible to show that $\bG^{\circ F}$ is characteristic in $\bG^F$. In other words, $\Aut(\bG^F)_{\bG^{\circ F}} = \Aut(\bG^F)$. However, as we will not need this here we do not attempt a proof.
\end{rem}

\begin{rem}
From the discussion in Remark~\ref{rem:cent-in-der} the size of $\rK(\bG^F)$ may be expressed as
%%%%
\begin{equation*}
\lvert \rK(\bG^F)\rvert
=
\lvert{\Hom(\bG^F/\rO^{p'}(\bG^F)\rZ(\bG^F),\rZ(\bG^F))}\rvert
=
\lvert \bG^F/\rO^{p'}(\bG^F)\rZ(\bG^F)\rvert
\end{equation*}
%%%%
which can be calculated explicitly. If $\bG = \Sp(V)$ this is $1$ and if $\dim(V)$ is odd then this is $2$. If $\bG = \GO(V)$ and $\dim(V)$ is even then this is either $2$ or $4$ depending on whether $\rZ(\bG^F) \leqslant \rO^{p'}(\bG^F)$. An explicit criterion for this to occur is given in Lemma~\ref{lem:omega-inv}.
\end{rem}

\subsection*{Dual groups and Dual automorphisms}
\begin{pa}\label{pa:centraliser-semisimple}
Let us now assume that $\bG = \SO(V)$ is a special orthogonal group. We let $V^{\star} \subseteq V$ be the subspace spanned by $\mathcal{V}\setminus \{v_0\}$, with $\mathcal{V}$ as in \eqref{eq:basis}. The dual of $\bG$ is then the group $\SL(V^{\star}\mid \mathcal{B}^{\star})$ where $\mathcal{B}^{\star}$ is the form from \ref{pa:bilinear-forms} satisfying $\mathcal{B}^{\star}(v,w) = (-1)^{\dim(V)}\mathcal{B}^{\star}(w,v)$ for all $v,w \in V^{\star}$. In particular, $\bG = \bG^{\star}$ if $\dim(V)$ is even. Recalling from \cite[5.3]{taylor:2018:action-of-automorphisms-symplectic} the notion of duality between isogenies, we have the following.
\end{pa}

\begin{lem}\label{lem:dual-asos}
Assume $\bG = \SO(V)$ is a special orthogonal group and let ${}^{\star} : \Gamma(\bG) \to \Gamma(\bG^{\star})$ be the unique monoid isomorphism satisfying $F_p^{\star} = F_p$ and $\gamma_1^{\star} = \gamma_1$. Then for any $\sigma \in \Gamma(\bG)$, we have $\sigma^{\star} \in \Gamma(\bG^{\star})$ is an asogeny dual to $\sigma$.
\end{lem}

\section{Action of Automorphisms on Unipotent Characters}\label{sec:aut-uni-char}
\begin{assumption}
From now on, we assume that $\bG$ is either $\GL(V\mid\mathcal{B})$ or $\SL(V\mid\mathcal{B})$ and $F \in \Gamma_q(\bG)$ is a Frobenius endomorphism, with $q$ a power of $p\neq 2$, as in Section~\ref{sec:classical-groups}.
\end{assumption}

We recall that the set $\mathcal{E}(\bG^F,1) \subseteq \Irr(\bG^F)$ of unipotent characters consists of those irreducible characters $\chi \in \Irr(\bG^F)$ whose restriction to $\bG^{\circ F}$ contains a unipotent character.

\begin{lem}\label{lem:ortho-unip}
Assume $\bG = \GO(V)$ is an orthogonal group with $\dim(V) \geqslant 5$. If $\chi \in \mathcal{E}(\bG^F,1)$ is a unipotent character, then $\chi^{\sigma} = \chi$ for any $\sigma \in \Aut(\bG^F)$.
\end{lem}

\begin{proof}
Let $\chi \in \mathcal{E}(\bG^F,1)$ be a unipotent character and let $\psi_1 \in \mathcal{E}(\bG^{\circ F},1)$ be an irreducible constituent of the restriction $\Res_{\bG^{\circ F}}^{\bG^F}(\chi)$. As $\rZ(\bG)^F = \rZ(\bG^{\circ})^F \leqslant \bG^{\circ F}$, it follows that $\omega_{\chi} = \omega_{\psi_1}$, which is trivial by Proposition~\ref{prop:cent-chars}(i). Thus if $\sigma \in \langle \gamma_2\rangle\rK(\bG^F)$, then the statement holds by (ii) of Lemma~\ref{lem:autos-triv-on-derived}.

As $\gamma_1$ is inner, it suffices to consider the case where $\sigma \in \Inn(\bG^{\circ},F)\langle F_p\rangle$. We have $\Res_{\bG^{\circ F}}^{\bG^F}(\chi)$ is either irreducible or the sum of two distinct irreducible constituents. If $\Res_{\bG^{\circ F}}^{\bG^F}(\chi) = \psi_1+\psi_2$ with $\psi_1 \neq \psi_2$ irreducible, then $\chi = \Ind_{\bG^{\circ F}}^{\bG^F}(\psi_1)$. It is well known that $\psi_1^{\sigma} = \psi_1$, so $\chi^{\sigma} = \chi$ in this case. Now consider the case where $\Res_{\bG^{\circ F}}^{\bG^F}(\chi) = \psi_1$ is irreducible. Then $\Ind_{\bG^{\circ F}}^{\bG^F}(\psi_1) = \chi + \theta\chi$, where $\theta$ is the inflation of the non-trivial character of $\bG^F/\bG^{\circ F}$.  Since $\sigma$ fixes $\theta$, it suffices to show that one character lying over $\psi:=\psi_1$ is fixed by $\sigma$.

Denote by $S$ the quotient group $\rO^{p'}(\bG^F)/\rZ(\rO^{p'}(\bG^F))$, which is a finite simple group. Recall that the unipotent character $\psi$ restricts irreducibly to $\rO^{p'}(\bG^F)$ and has $\rZ(\rO^{p'}(\bG^F))$ in its kernel. This follows from \cite[Prop.~11.3.8]{digne-michel:2020:representations-of-finite-groups-of-lie-type} because $\rO^{p'}(\bG^F)$ is the image of $\bG_{\simc}^F$ under a simply connected covering $\bG_{\simc} \to \bG$ and $\rZ(\rO^{p'}(\bG^F)) \leqslant \rZ(\bG^F)$. Let $\eta \in \Irr(S)$ be the deflation of $\Res_{\rO^{p'}(\bG^F)}^{\bG^{\circ F}}(\psi)$ to $S$. If $H \leqslant \Aut(S)$ is the stabiliser of $\eta$, then by \cite[Thm.~2.4]{malle:2008:extensions-of-unipotent-characters} there exists a character $\tilde{\eta} \in \Irr(H)$ extending $\eta$. If $\alpha : \bG^F \to \Aut(S)$ is the natural map, then $\langle \alpha(\bG^F),\sigma\rangle \leqslant H$.  Then inflating the restriction $\Res_{\alpha(\bG^F)}^{H}(\tilde{\eta})$ yields a $\sigma$-invariant irreducible character of $\bG^F$ extending $\psi$.
\end{proof}

\begin{pa}\label{pa:klein-four-sit}
We record here an elementary observation that will be used in the next proof and later sections. For $i\in\{1, 2\}$, let  $G_i$ be a finite group with a subgroup $N_i \leqslant G_i$ of index two. Let $G = G_1\times G_2$ and $N = N_1\times N_2$ and assume we have a subgroup $N \leqslant G^{\circ} \leqslant G$. For any subgroup $H \leqslant G$, we set $\bar{H} = HN/N$, so $\bar{G} = \bar{G}_1\times\bar{G}_2 \cong \rC_2\times \rC_2$ is a Klein four group and $\bar{G}^{\circ} \leqslant \bar{G}$.

Let $\chi \in \Irr(G^{\circ})$ be an irreducible character and let $\chi_i \in \Irr(N_i)$ be such that $\chi_1\boxtimes\chi_2 \in \Irr(N)$ is a constituent of the restriction $\Res_N^{G^{\circ}}(\chi)$. The group $\bar{G}$ acts naturally by conjugation on $\Irr(N)$ and
%%%%
\begin{equation*}
\bar{G}_{\chi_1\boxtimes\chi_2} = (\bar{G}_1)_{\chi_1} \times (\bar{G}_2)_{\chi_2}.
\end{equation*}
%%%%
If $\bar{G}^{\circ}$ is the diagonal embedding of $\rC_2$ into $\bar{G}$ then
%%%%
\begin{itemize}
	\item either $\bar{G}^{\circ}_{\chi_1\boxtimes\chi_2} \neq 1$ and $\chi = \Res_{G^{\circ}}^{G}(\widetilde{\chi}_1\boxtimes\widetilde{\chi}_2)$ with $\widetilde{\chi}_i \in \Irr(G_i)$ extending $\chi_i$
	\item or $\bar{G}^{\circ}_{\chi_1\boxtimes\chi_2} = 1$ and $\chi = \Ind_{N}^{G^{\circ}}(\chi_1\boxtimes\chi_2)$.
\end{itemize}

For the following we recall also that the group $(G/G^{\circ})^{\vee}$ of linear characters of $G/G^{\circ}$ also acts naturally on $\Irr(G)$ by tensoring with the inflation.
\end{pa}

\begin{prop}\label{prop:unip-qi}
Assume $\bG = \GO(V)$ and $s \in \bG^{\circ F}$ is a semisimple element such that $s^2 = 1$. We identify $A_{\bG^{\circ}}(s)^F$ with $\rC_{\bG^{\circ}}(s)^F/\rC_{\bG}^{\circ}(s)^F$. Let $\gamma_1 \in \Aso(\bG,F)$ be as in \ref{pa:gamma-auts} and assume $\gamma_1(s) = s$. If $\chi \in \mathcal{E}(\rC_{\bG^{\circ}}(s)^F,1)$ is a unipotent character, then the following hold:
%%%%
\begin{enumerate}
	\item $\chi^{F_p} = \chi$,
	\item if $s \neq 1$ then $\chi^{\gamma_1} = \chi$ if and only if $(A_{\bG^{\circ}}(s)^{F})^{\vee}_{\chi} \neq (A_{\bG^{\circ}}(s)^{F})^{\vee}$.
\end{enumerate}
\end{prop}

\begin{proof}
We have $V = V_1\oplus V_{-1}$, where $V_{\zeta}$ is the $\zeta$-eigenspace of $s$. Moreover, $\rC_{\bG}(s) = \GO(V_1) \times \GO(V_{-1})$. The quotient $A_{\bG}(s)^F \cong \rC_{\bG}(s)^F/\rC_{\bG}^{\circ}(s)^F$ is a Klein four group and $A_{\bG^{\circ}}(s)^F \cong \rC_{\bG^{\circ}}(s)^F/\rC_{\bG}^{\circ}(s)^F$ is the diagonally embedded subgroup. Let $\chi_1\boxtimes \chi_2 \in \mathcal{E}(\rC_{\bG}^{\circ}(s)^F,1)$ be an irreducible constituent of $\Res_{\rC_{\bG}^{\circ}(s)^F}^{\rC_{\bG^{\circ}}(s)^F}(\chi)$. By the discussion in \ref{pa:klein-four-sit}, either $\chi$ is the restriction of a character of $\rC_{\bG}(s)^F$, which is $F_p$-fixed by Lemma~\ref{lem:ortho-unip}, or is the induction of $\chi_1\boxtimes\chi_2$, which is known to be $F_p$-fixed. Hence, $\chi$ is $F_p$-fixed, giving (i). Statement (ii) follows from the description above.
\end{proof}

\section{Power Maps and Classical Groups}\label{sec:power-maps-2}
\begin{pa}\label{pa:uni-class-parts}
Recall that a partition of an integer $N \geqslant 0$ is a weakly decreasing sequence of positive integers $\mu = (\mu_1,\mu_2,\dots)$ such that $\sum_{i \in \mathbb{N}} \mu_i = N$. For any $m \in \mathbb{N}$, we set $r_m(\mu) := |\{i \in \mathbb{N} \mid \mu_i = m\}|$. If $\epsilon \in \{0,1\}$, we denote by $\mathcal{P}_{\epsilon}(N)$ the set of all partitions $\mu$ of $N$ such that $r_m(\mu) \equiv 0 \pmod{2}$ whenever $m \equiv \epsilon \pmod{2}$. For $\mu \in \mathcal{P}(N)$, we set
%%%%
\begin{align*}
a(\mu,\epsilon) &:= |\{m \in \mathbb{N} \mid r_m(\mu) \neq 0\text{ and }m\equiv 1+\epsilon \pmod{2}\}|,\\
\delta(\mu,\epsilon) &:= \begin{cases}
1 &\text{if $N$ is even and $r_m(\mu)$ is odd for some }m \equiv 1+\epsilon \pmod{2},\\
0 &\text{otherwise}.
\end{cases}
\end{align*}
%%%%

Let $\bG=\GL(V\mid \mathcal{B})$ with $N = \dim(V) \geqslant 0$ and $\epsilon \in \{0,1\}$  such that $\mathcal{B}(v,w) = (-1)^{\epsilon}\mathcal{B}(w,v)$ for all $v,w \in V$. We define $c(\bG) \in \{0,1\}$ so that $|\bG/\bG^{\circ}| = 2^{c(\bG)}$. We have a bijection $\mathcal{P}_{\epsilon}(N) \to \Clu(\bG)$, denoted by $\mu \mapsto \mathcal{O}_{\mu}$, such that $\mu$ gives the sizes of the Jordan blocks in the Jordan normal form of any $u \in \mathcal{O}_{\mu}$ in its natural action on $V$. Note that any unipotent class $\mathcal{O}_{\mu} \in \Clu(\bG)$ is contained in $\bG^{\circ}$ but may not be a single $\bG^{\circ}$-conjugacy class. Recall that a unipotent element $u \in \bG^{\circ}$ is said to be \emph{$\bG^{\circ}$-distinguished} if any torus of $\rC_{\bG}^{\circ}(u) = \rC_{\bG^{\circ}}^{\circ}(u)$ is contained in $Z^{\circ}(\bG^{\circ}) = Z^{\circ}(\bG)$. Equivalently $\rC_{\bG}^{\circ}(u)$ is a unipotent group as $\rZ^{\circ}(\bG)$ is trivial for our chosen $\bG$.
\end{pa}

\begin{lem}\label{lem:comp-grp-orders}
Let $\pi : \bG^{\circ} \to \bG_{\ad}$ be an adjoint quotient of $\bG^{\circ}$.  Then $\pi$ defines a bijection $\Clu(\bG^{\circ}) \to \Clu(\bG_{\ad})$ between the sets of unipotent classes. If $u \in \mathcal{O}_{\mu}$, then $A_{\bG}(u)$, $A_{\bG^{\circ}}(u)$, and $A_{\bG_{\ad}}(\pi(u))$ are all elementary abelian $2$-groups whose respective orders are the maximum of $1$ and: $2^{a(\mu,\epsilon)}$, $2^{a(\mu,\epsilon) - c(\bG)}$, and $2^{a(\mu,\epsilon) -c(\bG) - \delta(\mu,\epsilon)}$.
\end{lem}

\begin{proof}
This follows easily from \cite[Thm.~3.1]{liebeck-seitz:2012:unipotent-nilpotent-classes}, see also \cite[\S3.3.5]{liebeck-seitz:2012:unipotent-nilpotent-classes} and \cite[\S13.1]{carter:1993:finite-groups-of-lie-type}.
\end{proof}

\begin{lem}\label{lem:unique-complement}
If $u \in \bG$ is a $\bG^{\circ}$-distinguished unipotent element, then the subset $A \subseteq \rC_{\bG}(u)$ of all semisimple elements is a subgroup and $\rC_{\bG}(u) = \rC_{\bG}^{\circ}(u) \times A$. Moreover, the group $\bR = \rC_{\bG}(A)$ is reductive, $u \in \bR$ is regular unipotent, and $A = \rZ(\bR)$.
\end{lem}

\begin{proof}
Let $\epsilon \in \{0,1\}$ be such that $\mathcal{B}(v,w) = (-1)^{\epsilon}\mathcal{B}(w,v)$ for all $v,w \in V$. Let $\mu := (n_1,\dots,n_t) \in \mathcal{P}_{\epsilon}(N)$ be a partition such that $n_1 > \cdots > n_t$ and $n_i \equiv 1+\epsilon \pmod{2}$ for all $1 \leqslant i \leqslant t$. Correspondingly, we can choose an orthogonal decomposition $V_1\oplus\cdots\oplus V_t$ of $V$ such that $\dim(V_i) = n_i$. We also have a subgroup $\bR = \bR_1\times\cdots\times \bR_t$ of $\bG$ where $\bR_i = \GL(V_i\mid\mathcal{B})$ and we identify $\mathcal{B}$ with its restriction to $V_i$.

Fix a regular unipotent element $u \in \bR$.  Then $u \in \mathcal{O}_{\mu}$ is $\bG^{\circ}$-distinguished and all distinguished unipotent elements arise in this way, see \cite[Prop.~3.5]{liebeck-seitz:2012:unipotent-nilpotent-classes}. Recall that $\bR_i = \GO(V_i)$, not $\SO(V_i)$, when $\epsilon = 0$. We have $\rZ(\bR)$ is an elementary abelian $2$-group and $\rC_{\bG}(u) = \rC_{\bG}^{\circ}(u) \times \rZ(\bR)$. Indeed, $\rC_{\bG}^{\circ}(u) \cap \rZ(\bR)$ is trivial, as $\rC_{\bG}^{\circ}(u)$ is unipotent and $|\rZ(\bR)| = 2^t = 2^{a(\mu,\epsilon)} = |A_{\bG}(u)|$. Finally, it is clear that $\bR = \rC_{\bG}(A)$ and $A = \rZ(\bR)$.
\end{proof}

\begin{pa}
We now assume that $u \in \bG^{\circ F}$ is any unipotent element. Fix an $F$-stable maximal torus $\bS \leqslant \rC_{\bG}^{\circ}(u)$ and consider the subgroup $\bM = \rC_{\bG}(\bS)$, which is necessarily $F$-stable. We have an isomorphism
%%%%
\begin{equation}\label{eq:struc-M}
\bM \cong \GL(U_1) \times \cdots \times \GL(U_d) \times \GL(V_0 \mid \mathcal{B})
\end{equation}
%%%%
for some decompositions $V = V_1\oplus V_0$, as in \ref{pa:subspace-subgrp}, and $V_1 = (U_1\oplus U_1') \oplus \cdots \oplus (U_d\oplus U_d')$, where $U_i$ and $U_i'$ are totally isotropic subspaces. The element $u$ is $\bM^{\circ}$-distinguished, where $\bM^{\circ} = \rC_{\bG^{\circ}}(\bS)$  is an $F$-stable Levi subgroup of $\bG^{\circ}$. We now show that $u$ is regular in a reductive subgroup of $\bG$. This idea also appears in work of Testerman \cite[\S3]{testerman:1995:A1-type-overgroups}.
\end{pa}

\begin{prop}\label{prop:inj-R-cmp}
Assume $u \in \bG^F$ is a unipotent element and $\bS \leqslant \rC_{\bG}^{\circ}(u)$ is an $F$-stable maximal torus. Let $\bM = \rC_{\bG}(\bS)$ and denote by $A \subseteq \rC_{\bM}(u)$ the subset of semisimple elements. Then $u$ is regular unipotent in the $F$-stable reductive subgroup $\bR = \rC_{\bG}(A)$ and $A = \rZ(\bR)$. Moreover, the natural map $\rZ(\bR) \to A_{\bR}(u)$ is surjective and the natural map $A_{\bR}(u) \to A_{\bG}(u)$ is injective.
\end{prop}

\begin{proof}
Note that if $\bH = \GL(V)$, then a unipotent element $v \in \bH$ is distinguished if and only if it is regular. For such an element, $\rZ(\bH) \leqslant \rC_{\bH}(v)$ gives the subset of semisimple elements, which is clearly a subgroup. After Lemma~\ref{lem:unique-complement} and \eqref{eq:struc-M}, it follows immediately that $A = \rZ(\bR)$ is a subgroup, $\bR = \rC_{\bG}(A)$ is reductive, and $u \in \bR$ is regular.

The natural maps $\rZ(\bR)/Z^{\circ}(\bR) \to A_{\bR}(u) \to A_{\bM}(u)$ are isomorphisms, which follows from Lemma~\ref{lem:unique-complement}. Now, as $\bM = \rC_{\bG}(\bS)$ is the centralizer of a torus, a standard argument shows that the map $A_{\bM}(u) \to A_{\bG}(u)$ is injective, see \cite[1.4]{spaltenstein:1985:on-the-generalized-springer-correspondence}. Hence $A_{\bR}(u) \to A_{\bG}(u)$ is injective.
\end{proof}

\begin{cor}\label{cor:so-conj}
If $\bG = \GO(V)$, then any unipotent element $u \in \bG^{\circ F}$ is rational in $\bG^{\circ F}$.
\end{cor}

\begin{proof}
The group $\bR^{\circ}$ is a direct product of general linear and special orthogonal groups. As $u \in \bR^{\circ}$ is a regular element, it follows from Lemma~\ref{cor:reg-class-conj} that $u$ and $u^k$ are conjugate in $\bR^{\circ F}$ and hence also in $\bG^{\circ F}$.
\end{proof}

\begin{cor}\label{cor:sp-conj}
Assume $\bG = \Sp(V)$ with $\dim(V) = 2n$. We have $\bR = \bR^{\circ}$ is a direct product of symplectic and general linear groups. If $g \in \bR$ is such that ${}^gu = u^k$, then the map $\tau$ of Lemma~\ref{lem:perm-conj} is given by $\tau([\rC_{\bG}^{\circ}(u)a]) = [\rC_{\bG}^{\circ}(u)a\mathscr{L}(g)]$. In particular, if $u \in \mathcal{O}_{\mu}$ with $\mu \in \mathcal{P}_1(2n)$, then $u$ is $\bG^F$-conjugate to $u^k$ if and only if one of the following hold:
%%%%
\begin{enumerate}
	\item $r_{2m}(\mu)$ is even for all $1 \leqslant m \leqslant n$,
	\item $k \pmod{p} \in \mathbb{F}_q$ is a square.
\end{enumerate}
%%%%
\end{cor}

\begin{proof}
It is known that there exists an $F$-stable torus $\bS \leqslant \bG$ such that $\rC_{\bG}(u) = \bS\rC_{\bG}^{\circ}(u)$, see \cite[Prop.~9.3]{brunat-dudas-taylor:2020:unitriangular-shape-of-decomposition-matrices}. So the first statement follows from Lemma~\ref{lem:perm-conj}. Deconstructing the isomorphisms in \eqref{eq:struc-M} and Lemma~\ref{lem:unique-complement} we see that $\bR = \bR_1\times\cdots\times \bR_d\times\bR_{d+1}\times\cdots\times \bR_{d+t}$ is a direct product where $\bR_i \cong \GL(U_i)$, for $1 \leqslant i \leqslant d$, and $\bR_{d+i} \cong \Sp(V_i)$, for $1 \leqslant i \leqslant t$. In particular, we have $\dim(V_1)>\cdots > \dim(V_t)$ are all even. As these spaces have distinct dimensions it follows that the subgroups $\bR_{d+1},\dots,\bR_{d+t}$ are $F$-stable.

We write $u = u_1u_2\cdots u_{d+t}$ with $u_j \in \bR_j$ a regular unipotent element. It is clear that $u$ acts on $V$ as a direct sum of Jordan blocks with sizes
%%%%
\begin{equation*}
(\dim(U_1),\dim(U_1),\dots,\dim(U_d),\dim(U_d),\dim(V_1),\dots,\dim(V_t)).
\end{equation*}
%%%%
From this description we see that $A_{\bR}(u)$ is an elementary abelian $2$-group of order $2^t$, where $t\geqslant 0$ is the number of even numbers occurring an odd number of times in $\mu$. If (i) holds then $t = 0$ so $A_{\bR}(u)$ is trivial and it follows that $u$ and $u^k$ are conjugate in $\bR^F$ by Proposition~\ref{prop:coprime-perm-reg}.

Now suppose (i) does not hold, or equivalently $t > 0$. By \cite[Lem.~2.2]{taylor:2013:on-unipotent-supports} and the proof of \cite[Prop.~2.4]{taylor:2013:on-unipotent-supports} we see that $F$ acts trivially on $A_{\bG}(u)$. Thus the natural map $H^1(F,A_{\bR}(u)) \to H^1(F,A_{\bG}(u))$ is simply identified with the natural inclusion map $A_{\bR}(u) \to A_{\bG}(u)$. So $u$ and $u^k$ are conjugate in $\bG^F$ if and only if they are conjugate in $\bR^F$. This happens if and only if (ii) holds by Lemma~\ref{cor:reg-class-conj}.
\end{proof}

\section{Cuspidal Characters of Orthogonal Groups}\label{sec:cusp-chars}
\begin{pa}
In this section, we will give analogues of the calculations in \cite[\S14]{taylor:2018:action-of-automorphisms-symplectic} for the special orthogonal groups. For this, we introduce a little notation.  We refer the reader to \cite{taylor:2018:action-of-automorphisms-symplectic} for more details. For any integers $m,\delta \in \mathbb{N}_0$, we let
%%%%
\begin{equation*}
\Delta(\delta,m) = \begin{cases}
1 &\text{if }\delta = 0\text{ and }m \neq 0,\\
0 &\text{otherwise}.
\end{cases}
\end{equation*}
%%%%
Given integers $r,n \in \mathbb{N}_0$ and $d \in \{0,1\}$, we will denote by $\mathbb{X}_{n,d}^r := \mathbb{X}_{n,d}^{r,0}$ the sets of symbols defined in \cite[\S1]{lusztig-spaltenstein:1985:on-the-generalised-springer-correspondence}, see \cite[13.4]{taylor:2018:action-of-automorphisms-symplectic}.

For any $n$, we let ${}^{\dagger}\mathbb{X}_{n,0}^0\subseteq \mathbb{X}_{n,0}^0$ be the set of symbols $\ssymb{a_1&\cdots & a_m}{b_1 & \cdots & b_m}$ such that $\sum_{i=1}^m a_i > \sum_{i=1}^m b_i$. For compatibility, we also let ${}^{\dagger}\mathbb{X}_{n,1}^0 = \mathbb{X}_{n,1}^0$. If $e \in \mathbb{N}$, then we will need the following \emph{special symbols}
%%%%
\begin{align*}
\mathcal{S}_{e,1} &= \symb{0 & 1 & \cdots & e-1 & e}{1 & 2 & \cdots & e} \in {}^{\dagger}\mathbb{X}_{e(e+1),1}^0,\\
\mathcal{S}_{e,0} &= \symb{1 & 2 & \cdots & e}{0 & 1 & \cdots & e-1} \in {}^{\dagger}\mathbb{X}_{e^2,0}^0.
\end{align*}
%%%%
Recall we have a bijection ${}^{\dagger}\mathbb{X}_{n,1}^0 \to \Irr(W_n)$ and an injection ${}^{\dagger}\mathbb{X}_{n,0}^0 \to \Irr(W_n')$, defined as in \cite[4.5]{lusztig:1984:characters-of-reductive-groups}, where $W_n' \leqslant W_n$ are as in \ref{pa:Wgens}.

Let $\bG = \SO(V)$ be a special orthogonal group with $\dim(V) = 2n+\delta$ and $\delta \in \{0,1\}$. We have a regular embedding $\iota : \bG \to \wt{\bG}$, where $\wt{\bG} = \bG\cdot \rZ(\GL(V))$ is the corresponding conformal group, as in \ref{pa:conformal}. We denote by $\wt{\bG}^{\star}$ a dual group of $\wt{\bG}$ and by $\iota^{\star} : \wt{\bG}^{\star} \to \bG^{\star}$ a surjective homomorphism dual to the embedding. The dual group $\bG^{\star}$ is as in \ref{pa:centraliser-semisimple}. By the classification of quasi-isolated semisimple elements, a semisimple element $s \in \bG^{\star}$ is quasi-isolated if and only if $s^2=1$, see \cite[Prop.~4.11, Exmp.~4.10]{bonnafe:2005:quasi-isolated}. If $V_{\pm1}^{\star}(s)$ is the $(\pm1)$-eigenspace for $s$, then $\dim(V_{\pm1}^{\star}(s))$ is even and $\rC_{\bG^{\star}}^{\circ}(s) = \SL(V_1^{\star}(s) \mid \mathcal{B}^{\star}) \times \SL(V_{-1}^{\star}(s) \mid \mathcal{B}^{\star})$.

Now, if $s \in \bG^{\star F^{\star}}$ is $F^{\star}$-fixed and $\wt{s} \in \wt{\bG}^{\star F^{\star}}$ is such that $\iota^{\star}(\wt{s}) = s$, then $\iota^{\star}$ restricts to a surjective homomorphism $\rC_{\wt{\bG}^{\star}}(\wt{s})^{F^{\star}} \to \rC_{\bG^{\star}}^{\circ}(s)^{F^{\star}}$. Inflating through $\iota^{\star}$ we identify the unipotent characters of these groups. Conjugating $\wt{s}$ into $\wt{\bT}^{\star}$, we identify the Weyl group $\wt{\bW}^{\star}(s)$ of $C_{\wt{\bG}^{\star}}(\wt{s})$ with a subgroup of $\wt{\bW}^{\star}$. The following may be extracted from the work of Lusztig in a way entirely analogous to that of \cite[Lem.~14.3]{taylor:2018:action-of-automorphisms-symplectic}, so we omit the details. Here $n_{\wt{\chi}}$, for $\wt{\chi} \in \Irr(\wt{\bG}^F)$, is defined as in \cite[14.1]{taylor:2018:action-of-automorphisms-symplectic}
\end{pa}

\begin{lem}[Lusztig]\label{lem:lusztig-cusp-special}
Assume $\bG = \SO(V)$ and let $\dim(V) = 2n+\delta$ with $\delta \in \{0,1\}$. Let $\widetilde{s} \in \widetilde{\bG}^{\star F^{\star}}$ be a semisimple element whose image $s = \iota^{\star}(\widetilde{s}) \in \bG^{\star F^{\star}}$ is quasi-isolated with $2a = \dim(V_1^{\star}(s))$ and $2b = \dim(V_{-1}^{\star}(s))$. If $\wt{\chi} \in \mathcal{E}(\wt{\bG}^F,\wt{s})$ is a cuspidal character, then the following hold:
%%%%
\begin{enumerate}
	\item $\wt{\chi}$ is the unique cuspidal character in $\mathcal{E}(\wt{\bG}^F,\wt{s})$ and $a = e(e+\delta)$ and $b = f(f+\delta)$, for some non-negative integers $e$ and $f$, and either $e+\delta \geqslant 2$ or $f+\delta \geqslant 2$,
	\item $\wt{\chi}$ is contained in the family $\mathcal{E}(\widetilde{\bG}^F,\widetilde{s},\widetilde{\mathfrak{C}}) \subseteq \mathcal{E}(\widetilde{\bG}^F,\widetilde{s})$, where $\Irr(\wt{\bW}^{\star}(\wt{s}) \mid \widetilde{\mathfrak{C}})$ contains the special character $\mathcal{S}_{e,\delta} \boxtimes \mathcal{S}_{f,\delta} \in {}^{\dagger}\mathbb{X}_{a,\delta}^0 \times {}^{\dagger}\mathbb{X}_{b,\delta}^0$,
	\item $n_{\widetilde{\chi}} = 2^{e+f-\Delta(\delta,e)-\Delta(\delta,f)}$.
\end{enumerate}
%%%%
\end{lem}

\begin{pa}\label{pa:ident-uni-elts}
Let $\bG_{\ad}$ be the adjoint group of the same type as $\bG$, and let $\pi : \wt{\bG} \to \bG_{\ad}$ be an adjoint quotient. Then $\pi\circ\iota : \bG \to \bG_{\ad}$ is also an adjoint quotient of $\bG$. The maps $\pi\circ\iota$ and $\pi$ define bijections $\Clu(\bG) \rightarrow \Clu(\bG_{\ad}) \leftarrow \Clu(\wt{\bG})$ between the unipotent classes of these groups, see \cite[Prop.~5.1.1]{carter:1993:finite-groups-of-lie-type}. Hence we can label the unipotent classes as in Section~\ref{sec:power-maps-2}. Moreover, we have $|A_{\wt{\bG}}(u)| = |A_{\bG_{\ad}}(\pi(u))|$ by Lemma~\ref{lem:comp-grp-orders}.
\end{pa}

\begin{prop}\label{prop:mult1statement}
Assume $\bG = \SO(V)$ and $\wt{s} \in \wt{\bG}^{\star F}$ is a semisimple element whose image $s = \iota^{\star}(\wt{s}) \in \bG^{\star}$ is quasi-isolated. If $\wt{\chi} \in \mathcal{E}(\wt{\bG}^F,\wt{s})$ is a cuspidal irreducible character and $u \in \mathcal{O}_{\wt{\chi}}^*$ is a representative of the wave front set, then we have $n_{\wt{\chi}} = |A_{\wt{\bG}}(u)|$.
\end{prop}

\begin{proof}
Let $\dim(V) = 2n+\delta$ with $\delta \in \{0,1\}$. We may identify $E = \mathcal{S}_{e,\delta} \boxtimes \mathcal{S}_{f,\delta}$ with an irreducible character of $\wt{\bW}^{\star}(s)$. According to \cite[4.5(a), 6.3(b)]{lusztig:2009:unipotent-classes-and-special-Weyl}, we have $j_{\wt{\bW}^{\star}(s)}^{\wt{\bW}^{\star}}(E) = \mathcal{S}_{e,\delta} \oplus \mathcal{S}_{f,\delta} \in \mathbb{X}_{n,\delta}^0$ with the $j$-induction defined in \cite[\S5.2]{geck-pfeiffer:2000:characters-of-finite-coxeter-groups} and the addition of symbols defined as in \cite{lusztig-spaltenstein:1985:on-the-generalised-springer-correspondence}, see \cite[13.6]{taylor:2018:action-of-automorphisms-symplectic}. If $\Lambda_{0,\delta}^2$ is the unique element of $\mathbb{X}_{0,\delta}^2$ then we get a symbol $\Spr(j_{\wt{\bW}^{\star}(s)}^{\wt{\bW}^{\star}}(E)) = \mathcal{S}_{e,\delta} \oplus \mathcal{S}_{f,\delta}\oplus \Lambda_{0,\delta}^2 \in \mathbb{X}_{n,\delta}^2$. From this symbol, we may extract the partition of a unipotent class, as in \cite[\S2]{geck-malle:2000:existence-of-a-unipotent-support}, which will be exactly the class $\mathcal{O}_{\wt{\chi}}^*$. By symmetry, it suffices to treat the case where $0 \leqslant e \leqslant f$. We set $k = f-e \geqslant 0$.

Assume first that $\delta = 1$. Then the symbol $\Spr(j_{\wt{\bW}^{\star}(s)}^{\wt{\bW}^{\star}}(E))$ is exactly the same as that written in the first case occurring in the proof of \cite[Prop.~14.4]{taylor:2018:action-of-automorphisms-symplectic}. The unipotent class corresponding to this symbol is parameterised by the partition $\lambda = 2\mu+1$, where
%%%%
\begin{equation*}
\mu = (k+2e,\dots,k+1,k,k-1,k-1,\dots,1,1,0,0).
\end{equation*}
%%%%
As $\lambda$ contains $k + (2e+1) = e+f+1$ distinct odd numbers, we have $n_{\wt{\chi}} = 2^{e+f} = |A_{\wt{\bG}}(u)|$ by Lemmas~\ref{lem:lusztig-cusp-special} and \ref{lem:comp-grp-orders}.

Now assume $\delta = 0$. We have $\Spr(j_{\wt{\bW}^{\star}(s)}^{\wt{\bW}^{\star}}(E))$ is the symbol
%%%%
\begin{equation*}
\kbordermatrix{
& 0 & 1 & \cdots & k-1  & k    & k+1  & \cdots & k+e-1  \\
& 1 & 4 & \cdots & 3k-2 & 3k+2 & 3k+6 & \cdots & 3k+4e-2\\
& 0 & 3 & \cdots & 3k-3 & 3k   & 3k+4 & \cdots & 3k+4e-4
}.
\end{equation*}
%%%%
The unipotent class corresponding to this symbol is parameterised by the partition $\lambda = 2\mu+1$, where
%%%%
\begin{equation*}
\mu = (2k+4e-1,\dots,k+1,k,k-1,k-1,\dots,1,1,0,0).
\end{equation*}
%%%%
As $\lambda$ contains $k + 2e = e+f$ distinct odd numbers and at least one odd number occurs an odd number times, we have $n_{\wt{\chi}} = 2^{e+f-\Delta(0,e)-\Delta(0,f)} = |A_{\wt{\bG}}(u)|$.
\end{proof}

\begin{pa}
Recall from \ref{pa:gamma-auts} that we have defined a subgroup $\Gamma(\bG,F) \leqslant \Aut(\bG^F)$ generated by graph, field, and central automorphisms. We denote by $\mathscr{G}(\bG,F)$ the direct product $\Gamma(\bG,F) \times \gal$, which acts on $\Irr(\bG^F)$ through the actions of $\Aut(\bG^F)$ and $\gal$.
\end{pa}

\begin{thm}\label{prop:CquasiisolatedGGGR}
Assume $\bG = \SO(V)$ and $s \in \bG^{\star F^{\star}}$ is quasi-isolated. Then for any $\sigma \in \Aso(\bG,F)$, we have $\mathcal{E}(\bG^F,s)^{\sigma} = \mathcal{E}(\bG^F,s)$. Moreover, any cuspidal character $\chi \in \mathcal{E}(\bG^F,s)$ satisfies $\mathscr{G}(\bG,F)_{\chi} = \mathscr{G}(\bG,F)$.
\end{thm}

\begin{proof}
That the series $\mathcal{E}(\bG^F,s)$ is $\sigma$-stable is shown as in \cite[Thm.~14.6]{taylor:2018:action-of-automorphisms-symplectic}, but with the additional consideration of \cite[Lem.~3.4]{schaeffer-fry-taylor:2018:on-self-normalising-Sylow-2-subgroups} in the case that $\sigma \in \gal$. Let $\wt{s} \in \wt{\bG}^{\star}$ be a semisimple element such that $\iota(\wt{s}) = s$ and let $\wt{\chi} \in \mathcal{E}(\wt{\bG}^F,\wt{s})$ be a character covering $\chi$. Recall from Lemma \ref{lem:lusztig-cusp-special} that there is a unique cuspidal character in $\mathcal{E}(\wt{\bG}^F,\wt{s})$. The constituents of $\Res_{\bG^F}^{\wt{\bG}^F}(\wt{\chi})$ are the cuspidal characters in $\mathcal{E}(\bG^F,s)$, so $\Res_{\bG^F}^{\wt{\bG}^F}(\wt{\chi})$ is $\sigma$-stable. If $\chi$ extends to $\widetilde{\bG}^F$, then it is the unique cuspidal character contained in $\mathcal{E}(\bG^F,s)$ so must be $\mathscr{G}(\bG,F)$-invariant.

Assume now that $\chi$ does not extend to $\wt{\bG}^F$.   Then we must have $\dim(V)$ is even and $\Res_{\bG^F}^{\wt{\bG}^F}(\wt{\chi}) = \chi_1+\chi_2$ with $\chi_1 \neq \chi_2$ distinct irreducible characters. (Recall that $\wt{\bG}^F/\rZ(\wt{\bG}^F)\bG^F$ has size 2.) As $\chi \in \{\chi_1,\chi_2\}$, we must show that $\chi_i^{\sigma} = \chi_i$.

Let $u \in \mathcal{O}_{\chi}^*$. Using Lemma~\ref{lem:comp-grp-orders}, and the remarks in \ref{pa:ident-uni-elts}, we can compute $|A_{\wt{\bG}}(u)|$ and $|A_{\bG}(u)|$. It follows from the exact sequence in \cite[(2.1)]{taylor:2013:on-unipotent-supports} that $\rZ(\bG)$ embeds in $A_{\bG}(u)$. Arguing exactly as in the proof of \cite[Prop.~5.4]{taylor:2013:on-unipotent-supports}, we may find two unipotent elements $u_1,u_2 \in \mathcal{O}_{\chi}^{*F}$ such that $\langle \Gamma_{u_i},\chi_j\rangle = \delta_{i,j}$. We claim that the character $\Gamma_{u_i}$ is $\mathscr{G}(\bG,F)$-invariant. From this it follows that $\chi_i$ is $\mathscr{G}(\bG,F)$-invariant.

Firstly, we see that $\Gamma_{u_i}$ is $\gal$-invariant by Corollary~\ref{cor:so-conj} and \cite[Prop.~4.10]{schaeffer-fry-taylor:2018:on-self-normalising-Sylow-2-subgroups}. The action of $\Aso(\bG,F)$ on $\Gamma_{u_i}$ is described by \cite[Prop.~11.10]{taylor:2018:action-of-automorphisms-symplectic}. The same argument as used in \cite[Lem.~13.2]{taylor:2018:action-of-automorphisms-symplectic} shows that $\Gamma_{u_i}$ is $F_p$-invariant. We now just need to show that $\Gamma_{u_i}$ is $\GO(V)^F$-invariant. By Lemma~\ref{lem:comp-grp-orders}, and the description of $\mathcal{O}_{\chi}^*=\mathcal{O}_{\wt{\chi}}^*$ given in the proof of Proposition~\ref{prop:mult1statement}, we have $\rC_{\GO(V)^F}(u_i)$ is not contained in $\rC_{\bG^F}(u_i)$. As $\bG^F$ has index $2$ in $\GO(V)^F$, we have the $\bG^F$-class of $u_i$ is $\GO(V)^F$-invariant.
\end{proof}

\begin{pa}
Let $\bG = \GO(V)$ be an orthogonal group. We let $\mathscr{G}(\bG,F) = \Gamma(\bG,F) \times \gal$ as above. Recall that a character $\chi \in \Irr(\bG^F)$ is said to be cuspidal if the restriction $\Res_{\bG^{\circ F}}^{\bG^F}(\chi)$ contains a cuspidal irreducible character of $\bG^{\circ F}$. In the following sections, we will need to have some control over cuspidal characters of $\GO(V)$ to construct certain extensions. Using the characters constructed in Section~\ref{sec:GGGC-ext}, we get the following.
\end{pa}

\begin{thm}\label{prop:cuspidalisol}
Assume $\bG = \GO(V)$. Let $s \in \bG^{\circ \star F}$ be a quasi-isolated semisimple element. If $\chi \in \Irr(\bG^F)$ is a cuspidal character lying over $\chi^{\circ} \in \mathcal{E}({\bG^\circ}^F,s)$ with $s^2 = 1$, then $\chi$ is $\mathscr{G}(\bG,F)$-invariant.
\end{thm}

\begin{proof}
By Proposition~\ref{prop:CquasiisolatedGGGR}, the character $\chi^{\circ}$ is $\mathscr{G}(\bG^{\circ},F)$-invariant so extends to $\bG^F$.  Hence $\chi$ is one such extension. All other extensions are of the form $\theta\chi$, where $\theta$ is the inflation of an irreducible character of $\bG^F/\bG^{\circ F}$. Hence, it suffices to show that just one of the two extensions of $\chi^{\circ}$ is $\mathscr{G}(\bG,F)$-invariant.

Now, we may find a GGGC $\Gamma_u$ of $\bG^{\circ F}$, with $u \in \mathcal{O}_{\chi}^{*F}$ as in the proof of Proposition~\ref{prop:CquasiisolatedGGGR}, such that $\langle \Gamma_u,\chi^{\circ}\rangle = 1$. If $\rC_{\bG}(u)_2^F \leqslant \rC_{\bG}(u)^F$ is a Sylow $2$-subgroup, then $\rC_{\bG}(u)_2^F$ is not contained in $\bG^{\circ F}$, so by Lemma~\ref{lem:restrict-gamhat}, the character $\widehat{\Gamma}_{u,2}$ of $\bG^F$ extends $\Gamma_u$. It must contain only one of the two possible extensions of $\chi^{\circ}$. Moreover, arguing as in the proof of Proposition~\ref{prop:CquasiisolatedGGGR}, it follows from Theorem~\ref{prop:GGGCext-auts} and Proposition~\ref{prop:GGGCext-gal} that $\widehat{\Gamma}_{u,2}$ is $\mathscr{G}(\bG,F)$-invariant. Hence, the unique extension of $\chi^{\circ}$ that it contains must also be $\mathscr{G}(\bG,F)$-invariant.
\end{proof}

\section{Automorphisms and Quasi-Isolated Series}\label{sec:quasiisol}
\begin{pa}\label{pa:ortho-decomp}
Recall the notation for cuspidal pairs and Harish-Chandra series from Section~\ref{sec:gal}. If $s \in \bG^{\circ \star F^{\star}}$ is a semisimple element, then we let $\Cusp_s(\bG^{\circ},F) \subseteq \Cusp(\bG^{\circ},F)$ be the set of cuspidal pairs $(\bL,\lambda)$ such that $\mathcal{E}(\bG^{\circ F},\bL,\lambda) \subseteq \mathcal{E}(\bG^{\circ F},s)$, see \cite[Thm.~11.10]{bonnafe:2006:sln}. Then we have
%%%%
\begin{equation*}
\mathcal{E}(\bG^{\circ F},s) = \bigsqcup_{(\bL,\lambda) \in \Cusp_s(\bG^{\circ},F)/\bW^F} \mathcal{E}(\bG^{\circ F},\bL,\lambda),
\end{equation*}
%%%%
where the pairs are taken up to the natural action of $\bW^F$ by conjugation. In the following sections, as in \cite[\S15]{taylor:2018:action-of-automorphisms-symplectic}, we will study the actions of $\Aut(\bG^F)$ and $\gal$ on the Harish-Chandra series above. The proof of the following is identical to that of \cite[Lem.~15.4]{taylor:2018:action-of-automorphisms-symplectic}.
\end{pa}

\begin{lem}\label{lem:cusp-pair}
Let $\bG=\GL(V\mid \mathcal{B})$.  Assume $s \in \bG^{\circ \star F^{\star}}$ is a semisimple element with $s^2 = 1$ and let $(\bL,\lambda) \in \Cusp_s(\bG^{\circ},F)$ be a cuspidal pair. Then $\bL = Z^{\circ}(\bL) \times \bL_{\der}$ and there exists an orthogonal decomposition $V = V_1 \oplus V_0$ as in \ref{pa:subspace-subgrp} such that the following hold:
%%%%
\begin{enumerate}
	\item $\bL_{\der} = \bM_0^{\circ}$ and $Z^{\circ}(\bL) = \bS_1 \leqslant \bM_1^{\circ}$ is a maximally split torus of $\bM_1^{\circ}$,
	\item $\lambda = \lambda_1\boxtimes\psi$ with $\la_1\in\Irr(\bS_1^F)$, $\psi \in \mathcal{E}({\bM_0^\circ}^F,s_0)$ a cuspidal character, $s_0^2 = 1$, and $\la_1^2 = 1$.
\end{enumerate}
\end{lem}

\begin{pa}\label{strN}
Now assume $(\bL,\lambda) \in \Cusp_s(\bG^{\circ},F)$ is as in Lemma~\ref{lem:cusp-pair}. Using the notation of Section~\ref{sec:classical-groups}, we let $\tau = p_{{t}_1}p_{{t}_n} \in \rN_{\bG^{\circ}}(\bT)$. An argument with Weyl groups along the lines of that in \cite[15.6]{taylor:2018:action-of-automorphisms-symplectic} yields that
%%%%
\begin{equation*}
\rN_{\bG^{\circ}}(\bL) = (\rN_{\bM_1^{\circ}}(\bS_1) \times \bM_0^{\circ})\langle \tau\rangle \leqslant \rN_{\bG}(\bL) = \rN_{\bM_1}(\bS_1) \times \bM_0.
\end{equation*}
%%%%
With this in place, we have the following analogue of \cite[15.6]{taylor:2018:action-of-automorphisms-symplectic}.
\end{pa}

\begin{lem}\label{extGO}
Assume $\bG$, $s$, and $(\bL,\lambda)$, are as in Lemma~\ref{lem:cusp-pair}. Then for any extension $\Lambda(\lambda) \in \Irr(\rN_{\bG^{\circ F}}(\bL)_{\lambda})$ of $\la$, the following hold:
%%%%
\begin{enumerate}
	\item $\Lambda(\lambda)$ extends to a character of the form $\Lambda_1(\lambda_1) \boxtimes \Lambda_0(\psi)$ in $\rN_{\bG^F}(\bL)_{\lambda} = \rN_{\bM_1^F}(\bS_1)_{\lambda_1} \times \bM_0^F$, where $\Lambda_1(\la_1)$ and $\Lambda_0(\psi)$ are extensions of $\la_1$ and $\psi$ to $ \rN_{\bM_1^F}(\bS_1)_{\la_1}$ and $\bM_0^F$, respectively;
	\item $\Lambda(\lambda)$ is $\Gamma(\bG^{\circ},F)$-invariant.
\end{enumerate}
\end{lem}

\begin{proof}
By \cite[Lem.~15.7]{taylor:2018:action-of-automorphisms-symplectic}, and its proof, we can assume $\bG = \GO(V)$. To prove part (i), it suffices to know that $\la_1$ and $\psi$ extend,  using \ref{pa:klein-four-sit}. The group $\bS_1$ has a $\Gamma(\bG^{\circ},F)$-stable complement in $\rN_{\bM_1}(\bS_1)$ given by the permutation matrices. As $\lambda_1$ is linear, it extends trivially to $\rN_{\bM_1^F}(\bS_1)_{\lambda_1}$ and this must be $\Gamma(\bG^{\circ},F)$-invariant. That $\psi$ has an invariant extension is Theorem~\ref{prop:cuspidalisol}.
\end{proof}

\begin{proof}[of Theorem~\ref{mainthmgalSO}.(ii)]
Let $\bG = \SO(V)$ and assume $(\bL,\lambda) \in \Cusp_s(\bG,F)$ is as in Lemma~\ref{lem:cusp-pair}. After Lemma~\ref{extGO}, we see that $\lambda$ has an $F_p$-invariant extension $\Lambda(\lambda)$. Since $F_p$ is the identity on $W(\lambda)$, the statement follows using \cite[Thm.~4.6]{malle-spaeth:2016:characters-of-odd-degree}, as in \cite[Thm.~15.9]{taylor:2018:action-of-automorphisms-symplectic}.
\end{proof}

\section{Galois Automorphisms: Quasi-Isolated Series of \texorpdfstring{$\SO(V)$}{SO(V)}}\label{sec:galoisSO}

\begin{pa}\label{pa:not-cusp-pair}
Let $\bG=\GL(V\mid\mathcal{B})$ and let $s \in \bG^{\circ \star F^{\star}}$ be a quasi-isolated semisimple element. We will adopt the notation and assumptions of Lemma~\ref{lem:cusp-pair} throughout the whole of this section.

We now introduce some further notation to be used in this and the next section. Slightly abusing notation, we let $W^{\circ}(\lambda) = \rN_{\bG^{\circ F}}(\bL)_{\lambda}/\bL^F$ and $W(\lambda) = \rN_{\bG^F}(\bL)_{\lambda}/\bL^F$. Let $\Lambda$, $\Lambda_1$, and $\Lambda_0$ be extension maps with respect to $\bL^F\lhd \rN_{{\bG^\circ}^F}(\bL)$, $\bS_1^F\lhd \rN_{\bM_1^F}(\bS_1)$, and ${\bM_0^\circ}^F\lhd \bM_0^F$, respectively, such that
%%%%
\begin{equation*}
\Lambda(\lambda) = \Res_{\rN_{{\bG^\circ}^F}(\bL)_{\lambda}}^{\rN_{{\bG}^F}(\bL)_{\lambda}}\left(\Lambda_1(\lambda_1) \boxtimes \Lambda_0(\psi)\right)
\end{equation*}
%%%%
as in Lemma~\ref{extGO}. Then we similarly let $W(\lambda_1) = \rN_{\bM_1^F}(\bS_1)_{\lambda_1}/\bS_1^F$ and $W^\circ(\lambda_1) = \rN_{{\bM_1^\circ}^F}(\bS_1)_{\lambda_1}/\bS_1^F$.  We also have $W(\la)=W(\la_1)\times \bM_0^F/\bM_0^{\circ F}$, since the cuspidal character $\psi$ is invariant under $\bM_0^F$. Let $W^\circ(\la)= R(\la)\rtimes C(\la)$ and $W^\circ(\la_1)= R(\la_1) \rtimes C(\la_1)$ be decompositions as in \ref{pa:rel-stab}. 

For any $\sigma \in \gal$, note that $\Lambda_1(\lambda_1^\sigma)$ and $\Lambda_1(\lambda_1)^\sigma$ are extension of $\lambda_1^\sigma$ and that $\Lambda_0(\psi)^{\sigma}$ and $\Lambda_0(\psi^{\sigma})$ are extensions of $\psi^\sigma$.  Then Gallagher's Theorem \cite[Cor.~6.17]{isaacs:2006:character-theory-of-finite-groups} implies that there exist linear characters $\delta_{\lambda_1,\sigma} \in \Irr(W(\lambda_1))$ and $\delta_{\psi,\sigma} \in \Irr(\bM_0^F)$ such that $\Lambda_1(\lambda_1)^{\sigma} = \delta_{\lambda_1,\sigma}\Lambda_1(\lambda_1^{\sigma})$ and similarly $\Lambda_0(\psi)^{\sigma} = \delta_{\psi,\sigma}\Lambda_0(\psi^{\sigma})$. Moreover, for $\delta_{\la, \sigma}$ such that $\Lambda(\la)^\sigma=\delta_{\la, \sigma}\Lambda(\la^\sigma)$ as in Theorem~\ref{thm:GaloisAct}, we have
%%%%
\begin{equation*}
\delta_{\lambda,\sigma} = \Res_{W^{\circ}(\lambda)}^{W(\la)} \left(\delta_{\lambda_1,\sigma} \boxtimes \delta_{\psi,\sigma}\right).
\end{equation*}
%%%%
Further, recall that the characters $\gamma_{\la,\sigma}$ and $\delta_{\la,\sigma}'$ of Theorem~\ref{thm:GaloisAct} are characters of $C(\la)$.

With this, we have mostly reduced ourselves to the case of principal series characters and cuspidal characters. The rest of this section is devoted to completing the proof of Theorem~\ref{mainthmgalSO}. We first consider the case of principal series characters.
\end{pa}

\begin{prop}\label{prop:isolprincipalSO}
Assume $\bG = \GO(V)$, $F \in \Gamma_q(\bG)$, and $\lambda\in\Irr(\bT^F)$ is a character such that $\lambda^2=1$. Then each member of $\mathcal{E}(\bG^{\circ F}, \bT, \la)$ is invariant under $\gal$.
\end{prop}

\begin{proof}
Note that in this case, $\la=\la_1$.  As in \ref{pa:notnT}, we have $\bT = \bT_1\times\cdots \times \bT_n$. Each $\bT_i$ is $F$-stable, so $\bT^F = \bT_1^F \times \cdots \times \bT_n^F$. Acting with $\bW^F$, we can assume that $\lambda = 1_{\bT_1^F}\boxtimes\cdots \boxtimes 1_{\bT_a^F}\boxtimes \varepsilon_{\bT_{a+1}^F}\boxtimes \cdots \boxtimes \varepsilon_{\bT_{n-1}^F}\boxtimes \lambda_n$ for some $0 \leqslant a < n$, with $\lambda_n \in \{1_{\bT_n^F},\varepsilon_{\bT_n^F}\}$. If $F$ is split, then we will also assume that $\lambda_n = \varepsilon_{\bT_n^F}$. Let $b = n-a$.

We will assume that $\dim(V)$ is even, as the statement in the case that $\dim(V)$ is odd follows from a simplified version of the argument. Note that $\bW^F\leqslant \bW = W_n$ is the subgroup $W_n''$ of type $\B_{n-1}$, as in \ref{pa:Wgens}, when $F$ is twisted. In the notation of \ref{pa:subspace-subgrp}, we have $W(\lambda) \leqslant \bW_1\times \bW_0$ is $W_a\times W_b$ if $F$ is split and $W_a\times W_b''$, of type $\B_a\times\B_{b-1}$, if $F$ is twisted. Moreover, using Lemma~\ref{lem:stab-lam-tilde} and \ref{pa:conformal}, $R(\la) \leqslant \bW_1\times\bW_0$ is the subgroup $W_a'\times W_b'$ of type $\D_a\times \D_b$ if $F$ is split and $W_a'\times W_b''$, of type $\D_a\times \B_{b-1}$, if $F$ is twisted. In all cases, the group $C(\la)$ is generated by $u_1$, in the notation of \ref{pa:Wgens}, which has order two and even length in $\bW$. Then for each $\sigma\in\gal$, Lemma~\ref{lem:genrsigma} implies $\gamma_{\la,\sigma}=1$, and we further have $\Lambda_\la(n_{u_1})^2=\la(1)=1$ for any extension $\Lambda_\la:=\Lambda(\la)$ of $\la$ to $\rN_{\bG^{\circ F}}(\bg{T})_\la$. Then $\Lambda_\la(n_{u_1})\in\{\pm1\}$ is fixed by $\sigma$, and hence $\delta_{\la, \sigma}'=1$.

It remains to consider the values of characters of $W^{\circ}(\la)$ and its corresponding Hecke algebra. If $F$ is split, then the sequence of subgroups $R(\lambda) \leqslant W^{\circ}(\lambda) \leqslant W(\lambda)$ fits the setup of \ref{pa:klein-four-sit}. In particular, each irreducible character of $W^{\circ}(\lambda)$ either extends to $W(\lambda)$ or is induced from $R(\lambda)$. Moreover, there is a corresponding sequence of subalgebras of the Hecke algebra, and the same statement must be true for the representations of these algebras, by \cite[\S10.4.2]{geck-pfeiffer:2000:characters-of-finite-coxeter-groups}.
Using this and the fact that $R(\la)$ and $W(\la)$ are products of type $\B$ and $\D$ groups, it follows that all irreducible representations of the algebras corresponding to $R(\lambda)$ and $W(\lambda)$ are defined over $\mathbb{Q}$.  This implies that all irreducible representations of the Hecke algebra $\End_G(R_\bT^\bG(\la))$ are also defined over $\mathbb{Q}$. When $F$ is twisted, a similar argument applies. In particular, we have $\eta=\eta^\sigma=\eta^{(\sigma)}$ for all $\eta\in\Irr(W(\la))$ and all $\sigma\in\gal$.  This completes the proof, using Theorem~\ref{thm:GaloisAct}.
\end{proof}

\begin{proof}[of Theorem~\ref{mainthmgalSO}.(i)]
Let $\bG = \GO(V)$ and $\sigma \in \gal$. As in the proof of Proposition~\ref{prop:isolprincipalSO}, we will just treat the case $\dim(V)$ even. After Proposition~\ref{prop:isolprincipalSO}, it suffices to prove the statement for all members of a Harish-Chandra series $\mathcal{E}(\bG^{\circ F},\bL,\lambda)$ with $(\bL,\lambda) \in \Cusp_s(\bG,F)$ a cuspidal pair such that $\bL \neq \bT$. Recall the notation introduced in \ref{pa:not-cusp-pair}. By Theorem~\ref{prop:cuspidalisol}, the character $\delta_{\psi,\sigma}$ is trivial. Moreover, it follows from the proof of Proposition~\ref{prop:isolprincipalSO} that $\delta_{\lambda_1,\sigma}$ is trivial on $C(\la_1)$.  We let $m = \dim(\bS_1)$. As $\bL \neq \bT$, we have $0 \leqslant m \leqslant n-2$.  Then we similarly see that $\delta_{\lambda_1,\sigma}$ is trivial on $n_{t_m}$, which generates $W(\la_1)/W^\circ(\la_1)$, implying $\delta_{\lambda,\sigma}'$ must be trivial.

 As in the case of Proposition~\ref{prop:isolprincipalSO}, we may arrange that $\lambda_1 = 1_{\bT_1^F}\boxtimes\cdots\boxtimes 1_{\bT_a^F}\boxtimes \varepsilon_{\bT_{a+1}^F}\boxtimes \cdots \boxtimes \varepsilon_{\bT_m^F}$ for some $0 \leqslant a \leqslant m$. Let $b = m-a$. Identifying $\rN_{\bG^{\circ F}}(\bL)/\bL^F$ with the Coxeter group $W_m$, as in Lemma~\ref{lem:normaliser}, we have the group $W^{\circ}(\lambda)$ is the subgroup $W_a\times W_b$ of \ref{pa:Wgens}. Using the notation of \ref{pa:conformal}, we let $\wt{\lambda} = \lambda_1\boxtimes\wt{\psi} \in \Irr(\wt{\bL}_J^F)$ with $\wt{\psi} \in \Irr(\wt{\bM}_0^{\circ F})$ a character covering $\psi$.

By Lemma~\ref{lem:stab-lam-tilde}, and the description of the action of $u_m$ in \ref{pa:conformal}, we see that $W^{\circ}(\wt{\lambda}) = R(\lambda) = W_a\times W_b' \leqslant W^{\circ}(\lambda)$ is of type $\B_a\times\D_b$. Hence $C(\lambda)$ is generated by $u_m$ as in Lemma~\ref{lem:normaliser}. The element $u_m$ has even length in $\bW$, so $\gamma_{\lambda,\sigma}$ is trivial by Lemma~\ref{lem:genrsigma}. Finally, as $W(\lambda)$ is a product of type $\B$ Weyl groups, we have $\eta^{(\sigma)} = \eta$ for all $\eta \in \Irr(W(\lambda))$, so the result follows from Theorem~\ref{thm:GaloisAct}.
\end{proof}

\begin{proof}[of Theorem~\ref{mainthmgalSO}.(iii)]
We may again assume $\dim(V)$ is even and $F = F_q \in \Gamma_q(\bG)$ is split and $\gamma = \gamma_1$ is as in \ref{pa:auts-ortho}. Let $(\bL,\lambda) \in \Cusp(\bG,F)$ be a cuspidal pair and set $\wt{\bL} = \bL\cdot \rZ(\wt{\bG}) \leqslant \wt{\bG}$. Suppose $\wt{\lambda} \in \Irr(\wt{\bL}^F)$ covers $\lambda$. Identifying $\wt{W}(\wt{\lambda})$ with $W(\wt{\lambda}) \leqslant W(\lambda)$, we have by \cite[Thm.~13.9]{bonnafe:2006:sln} that
%%%%
\begin{equation*}
\Res_{\bG^F}^{\wt{\bG}^F}\left(R_{\wt{\bL}}^{\wt{\bG}}(\wt{\lambda})_{\wt{\eta}}\right) = R_{\bL}^{\bG}(\lambda)_{\Ind^{W(\lambda)}_{W(\wt{\la})}\wt{\eta}}
\end{equation*}
%%%%
for any $\wt{\eta} \in \Irr(W(\wt{\lambda}))$. Here we extend the map $\Irr(W(\lambda)) \to \Irr(\bG^F)$, given by $\eta\mapsto R_{\bL}^{\bG}(\lambda)_{\eta}$, linearly to all class functions.

Consider first the case where $\bL = \bT$. As above we arrange that $\lambda = 1_{\bT_1^F}\boxtimes\cdots \boxtimes 1_{\bT_a^F}\boxtimes \varepsilon_{\bT_{a+1}^F}\boxtimes \cdots \boxtimes \varepsilon_{\bT_n^F}$ for some $1 < a < n$. Note we assume $s \neq \pm1$. Using the above identity, as well as \cite[Thm. 4.6]{malle-spaeth:2016:characters-of-odd-degree}, the statement holds in this case. Note, by Lusztig's classification this series contains all the characters $\chi \in \mathcal{E}(\bG^F,s)$ satisfying the condition $\Inn(\bG,F)_{\chi} = \Inn(\bG,F)$. Now assume $\bL \neq \bT$. The graph automorphism $\gamma_1$ induces the identity on $W(\lambda)$ so the statement again follows from \cite[Thm. 4.6]{malle-spaeth:2016:characters-of-odd-degree}.
\end{proof}

\section{Galois Automorphisms: Quasi-Isolated Series of \texorpdfstring{$\Sp(V)$}{Sp(V)}}\label{sec:gal-qiso-sp}
\begin{pa}\label{strNSp2}
In this section, we consider the action of $\gal$ on quasi-isolated series when $\bG = \Sp(V)$. As in \ref{pa:not-cusp-pair}, we need to understand the action of $\gal$ on principal series characters and cuspidal characters. The following treats the cuspidal case and extends \cite[Thm.~14.6]{taylor:2018:action-of-automorphisms-symplectic} to the setting of Galois automorphisms.
\end{pa}

\begin{thm}\label{thm:cuspidalsSp}
Assume $\bG = \Sp(V)$ and $s \in \bG^{\star F^{\star}}$ is a quasi-isolated semisimple element. Let $\sigma \in \gal$ be a Galois automorphism, let $k \in \mathbb{Z}$ be an integer coprime to $p$ such that $\xi^{\sigma} = \xi^k$ for all $p$th roots of unity $\xi \in \overline{\Q}^{\times}$, and let $\omega=(-1)^{(p-1)/2}$. If $\chi \in \mathcal{E}(\bG^F,s)$ is a cuspidal character, then we have $\chi^{\sigma} = \chi$ if and only if one of the following holds:
%%%%
\begin{enumerate}
	\item $s = 1$, so that $\chi$ is unipotent,
	\item $k \pmod{p} \in \mathbb{F}_q$ is a square.
\end{enumerate}
In particular, we have $\Q(\chi)=\Q$ if $s=1$ or $q$ is square, and $\Q(\chi)=\Q(\sqrt{\omega p})$ otherwise. 
\end{thm}

\begin{proof}
Unipotent characters of classical groups are rational-valued, so we may assume $s \neq 1$. Assume $\mathcal{O}_{\chi}^*$ is the wave-front set of $\chi$, as in \cite[14.1]{taylor:2018:action-of-automorphisms-symplectic}. Let $\dim(V)= 2n$ and let $\mu \in \mathcal{P}_1(2n)$ be the partition parameterising the unipotent class $\mathcal{O}_{\chi}^*$, under the bijection in \ref{pa:uni-class-parts}. From the proof of \cite[Prop.~14.4]{taylor:2018:action-of-automorphisms-symplectic}, we see that $r_{2m}(\mu) = 1$ for some $0 < m \leqslant n$. Then arguments along the lines of those used in Proposition~\ref{prop:CquasiisolatedGGGR}, instead appealing to Corollary~\ref{cor:sp-conj} and \cite[Prop.~4.10]{schaeffer-fry-taylor:2018:on-self-normalising-Sylow-2-subgroups}, complete the proof of the first statement.
The statement about character fields follows from this, combined with \eqref{eq:gaussum} and  Lemma~\ref{lem:kmodp}, since if $s\neq 1$ and $q$ is not square, $\chi$ is fixed by $\sigma$ if and only if $\sqrt{\omega p}$ is.
\end{proof}

\begin{rem}
This implies that the Harish-Chandra series in \ref{pa:ortho-decomp} can be permuted non-trivially by the action of $\gal$, unlike the case of the action of $\Aut(\bG^F)$.
\end{rem}

\begin{lem}\label{lem:deltatypeC} 
Assume $\bG = \Sp(V)$. Let $\la\in\Irr(\bT^F)$ be a nontrivial character such that $\la^2=1$ and let $\sigma\in\gal$. Then
%%%%
\begin{enumerate}
\item $\eta=\eta^\sigma=\eta^{(\sigma)}$ for every $\eta\in\Irr(W(\la))$;
\item $l(w_2)$ is odd for $1\neq w_2\in C(\la)$;
\item if $q\equiv 1\pmod 4$, then any extension of $\lambda$ to $\rN_{\bG^F}(\bg{T})_\la$ is fixed by $\sigma$; in particular, $\delta_{\la, \sigma}=\delta'_{\la,\sigma}$ is trivial;
\item if $q\equiv 3\pmod 4$, then there is a $\sigma$-invariant extension of $\lambda$ to $\rN_{\bG^F}(\bg{T})_\la$ if and only if $\sigma$ fixes $4$th roots of unity.  In particular, $\delta_{\la, \sigma}=\delta'_{\la,\sigma}$ is trivial if and only if $\sigma$ fixes fourth roots of unity.
\end{enumerate}
\end{lem} 

\begin{proof}
As in the proof of Proposition~\ref{prop:isolprincipalSO}, we arrange that $\lambda = 1_{\bT_1^F}\boxtimes\cdots \boxtimes 1_{\bT_a^F}\boxtimes \varepsilon_{\bT_{a+1}^F}\boxtimes \cdots \boxtimes \varepsilon_{\bT_n^F}$ for some $0 \leqslant a < n$. Note we have $b = n-a > 0$ because $\la$ is assumed to be nontrivial. Taking $m=a$ in \ref{pa:subspace-subgrp}, we have $W(\la) = \bW_1 \times \bW_0$ is a Weyl group of type $\B_a \times \B_b$. This gives (i). Letting $\wt{\lambda} = \lambda \boxtimes 1_{\bZ^F} \in \Irr(\wt{\bT}^F)$, we see using Lemma~\ref{lem:stab-lam-tilde} that $R(\la) = W(\wt{\lambda}) = W_a\times W_b'$ and $C(\la) = \langle s_n \rangle$, which proves (ii).

Let $c := n_{s_n} \in \rN_{\bG^F}(\bT)_{\lambda}$ represent $s_n \in C(\la)$, and let $\xi \in \bT_n^F \cong \mathbb{F}_q^{\times}$ be a generator of the cyclic group. Since $c^2=\mathbf{d}_n(-1)$ and $-1 = \xi^{(q-1)/2}$, we see that for any extension $\Lambda_\la:=\Lambda(\la)$ of the linear character $\la$ to $\rN_{\bG^F}(\bg{T})_\la$, we have
%%%%
\begin{equation}\label{valonc}
\Lambda_\la(c)^2=\Lambda_\la(c^2)=\varepsilon_{\bT_n^F}(\xi^{(q-1)/2})=\varepsilon_{\bT_n^F}(\xi)^{(q-1)/2}=(-1)^{(q-1)/2}.
\end{equation}
%%%%
By combining this with Lemma~\ref{lem:extRla}, we see that if $q\equiv 1\pmod 4$, then all values of $\Lambda_\la$ are in $\{\pm1\}$, giving (iii).
On the other hand, if $q\equiv 3\pmod 4$, then \eqref{valonc} implies $\Lambda_\la$ must take primitive fourth roots of unity as values on $c$, proving (iv).
\end{proof}

\begin{pa}
When $\lambda\in\Irr(\bT^F)$ is a nontrivial character satisfying $\lambda^2=1$, combining Lemmas~\ref{lem:deltatypeC} and \ref{lem:genrsigma} yields a complete description of $\gamma_{\la,\sigma}$ and $\delta_{\la,\sigma}$, which can be viewed as characters of $C(\la)\cong \rC_2$.  In particular, note that these depend only on $q$ and $\sigma$, and we obtain the following.
\end{pa}

\begin{cor}\label{cor:rdelta}
Assume $\bG = \Sp(V)$.  Let $\omega \in \{\pm1\}$ be such that $p\equiv\omega\pmod 4$. Let $\la\in\Irr(\bT^F)$ be a nontrivial character such that $\la^2=1$ and let $\sigma\in\gal$.  Identifying $\Irr(C(\lambda))$ with $\{\pm1\}$, the following hold:
\begin{enumerate}
\item if $q$ is a square, then $\gamma_{\la,\sigma}\delta_{\la,\sigma}=1$;
\item if $q$ is not a square, then $\gamma_{\la,\sigma}\delta_{\la,\sigma}=\alpha$, where $\alpha\in\{\pm1\}$ is such that$\sqrt{\omega p}^\sigma=\alpha\sqrt{\omega p}$.  
\end{enumerate}
\end{cor}

With this in place, we may now give a proof of Theorem~\ref{mainthmgalSp}.

\begin{proof}[of Theorem~\ref{mainthmgalSp}]
This follows immediately from Theorem~\ref{thm:cuspidalsSp}, Lemma~\ref{lem:deltatypeC}, and Corollary~\ref{cor:rdelta}, using Theorem~\ref{thm:GaloisAct} and Lemma~\ref{lem:cusp-pair}.
\end{proof}

\section{Galois Automorphisms Relevant for the McKay--Navarro Conjecture II}\label{sec:actionH-Sp}
\begin{pa}
We now return to the setting of Section~\ref{sec:actionH} and consider those Galois automorphisms relevant to the McKay--Navarro conjecture in the context of symplectic groups. Such explicit statements will be particularly useful in future work. As in Section~\ref{sec:actionH}, we assume $\sigma \in \mathcal{H}$ is fixed and $r \geqslant 0$ is the integer such that $\zeta^{\sigma} = \zeta^{\ell^r}$ for all roots of unity $\zeta \in \overline{\Q}^{\times}$ of order coprime to $\ell$.
\end{pa}

\begin{lem}\label{cor:rdeltaodd}
Assume $\bG = \Sp(V)$ and $\ell \neq p$ are both odd. Let $\la\in\Irr(\bT^F)$ be a nontrivial character such that $\la^2=1$. Identifying $\Irr(C(\la))$ with the group $\{\pm1\}$, we have:
%%%%
\begin{enumerate}
\item If $q$ is square or $r$ is even, then $\gamma_{\la,\sigma}\delta_{\la,\sigma}=1$.  
\item If $q$ is not square and $r$ is odd, then $\gamma_{\la,\sigma}\delta_{\la,\sigma}=\legendre{\ell}{p}$.
\end{enumerate}
%%%%
In particular, note that if $\ell\mid(q-1)$, then $\gamma_{\la,\sigma}\delta_{\la,\sigma}=1$ unless $q\equiv\ell\equiv3\pmod4$ and $r$ is odd, in which case $\gamma_{\la,\sigma}\delta_{\la,\sigma} = \delta_{\la,\sigma}=-1$.
\end{lem}

\begin{proof}
This follows from Lemma~\ref{lem:omegapH} and Corollary~\ref{cor:rdelta}.
\end{proof}

\begin{lem}\label{cor:rdelta2}
Assume $\bG = \Sp(V)$ and $\ell=2$. Let $\la\in\Irr(\bT^F)$ be a nontrivial character such that $\la^2=1$.  Then $\gamma_{\la,\sigma}\delta_{\la,\sigma}\in\Irr(C(\la))\cong\{\pm1\}$ satisfies:
%%%%
\begin{itemize}
\item If $q\equiv \pm1\pmod 8$, then $\gamma_{\la,\sigma}\delta_{\la,\sigma}=1$.
\item If $q\equiv \pm 3\pmod 8$, then $\gamma_{\la,\sigma}\delta_{\la,\sigma}=(-1)^r$.
\end{itemize}
%%%%
\end{lem}

\begin{proof}
This follows from Lemma~\ref{lem:omegapH2} and Corollary~\ref{cor:rdelta}.
\end{proof}

\begin{rem} We end this section with some useful remarks on specific situations, which follow from Theorem~\ref{thm:cuspidalsSp}, Lemma~\ref{lem:deltatypeC}, and Lemma~\ref{cor:rdeltaodd}.
Let $\bG = \Sp(V)$ with $\ell \neq p$ and $p$ odd.
%%%%
\begin{enumerate}
\item For $\ell\mid (q-1)$ odd, every member of $\mathcal{E}(\bG^F,s)$ with $s\neq 1$ quasi-isolated is fixed by all of $\galh$ if and only if at least one of $q$ or $\ell$ is congruent to $1$ modulo $4$.

\item There has been much recent interest in the particular case that $\sigma:=\sigma_1$ is the Galois automorphism fixing $\ell'$-roots of unity and mapping $\ell$-power roots of unity to their $\ell+1$ power, see, e.g., \cite{navarrotiepexponent,RSV20-1,NRSV21}. In this situation, every member of any $\mathcal{E}(\bG^F,s)$ with $s^2=1$ is fixed by $\sigma_1$.  Hence, the action of $\sigma_1$ on $\Irr_{\ell'}(\bG^F)$ in this case is determined by the action on irreducible characters of general linear groups, which is well understood, see, e.g. \cite{SrinivasanVinroot-1, schaeffer-fry-taylor:2018:on-self-normalising-Sylow-2-subgroups}.
\end{enumerate}
\end{rem}

\section{Automorphisms Not Lifting to Asogenies}\label{sec:notaso}
\begin{pa}\label{pa:}
In this section, we assume that $\bG = \SO(V)$. We wish to consider the action of the group of automorphisms $\rK(\bG^F)$, defined in Proposition~\ref{prop:aut-grp-quasisimp}, on the irreducible characters $\Irr(\bG^F)$. By \cite[Prop.~24.21]{malle-testerman:2011:linear-algebraic-groups}, and a theorem of Steinberg \cite[Thm.~12.4]{steinberg:1968:endomorphisms-of-linear-algebraic-groups}, we have $\bG^F/\rO^{p'}(\bG^F) \cong \rZ(\bG^F)$ has order $2$. Thus, by Remark~\ref{rem:cent-in-der} we have $\rK(\bG^F)$ is non-trivial only when $\rZ(\bG^F) \leqslant \rO^{p'}(\bG^F)$.

The classification of when $\rZ(\bG^F) \leqslant \rO^{p'}(\bG^F)$ is well known and is typically calculated using a result of Zassenhaus \cite[\S2, Thm.]{zassenhaus:1962:on-the-spinor-norm} on the spinor norm. The following slightly more general statement is also well known, but we include a proof that is in keeping with the theory of algebraic groups.
\end{pa}

\begin{lem}\label{lem:omega-inv}
Assume $\bG = \SO(V)$ and $n = \lfloor\dim(V)/2\rfloor$. We let $\epsilon$ be $1$ if $F$ is split and $-1$ if $F$ is twisted. Let $0 \leqslant a < n$ be an integer and set $b=n-a$. Denote by $\mathcal{O} \subseteq \bG$ the $\bG$-conjugacy class of $s = \mathbf{d}_{a+1}(-1)\cdots \mathbf{d}_n(-1)$ so that $\dim(V_{-1}(s)) = 2b$, where $V_{-1}(s)$ is the $(-1)$-eigenspace of $s$.
%%%%
\begin{enumerate}
	\item If $s \in \rZ(\bG)^F$, then $s$ is contained in $\rO^{p'}(\bG^F)$ if and only if $q^n \equiv \epsilon \pmod{4}$,
	\item If $s \not\in \rZ(\bG)^F$, then $\mathcal{O}' = \mathcal{O} \cap \rO^{p'}(\bG^F)$ is a single $\rO^{p'}(\bG^F)$-conjugacy class and it contains $s$ if and only if $q^b \equiv \epsilon \pmod{4}$.
\end{enumerate}
\end{lem}

\begin{proof}
Let $\pi : \bG_{\simc} \to \bG$ be a simply connected cover of $\bG$. We denote again by $F$ a Frobenius endomorphism of $\bG_{\simc}$ such that $F\circ\pi = \pi\circ F$. Let $\bT_{\simc} = \pi^{-1}(\bT)$, an $F$-stable maximal torus of $\bT$. The image $\pi(\bG_{\simc}^F)$ is precisely $\rO^{p'}(\bG^F)$. If $t \in \bG^F$ then $t \in \rO^{p'}(\bG^F)$ if and only if $\pi^{-1}(t) \cap \bG_{\simc}^F \neq \emptyset$. Note that if it is non-empty then $\pi^{-1}(t) \cap \bG_{\simc}^F = \pi^{-1}(t)$ because $F$ fixes $\Ker(\pi)$.

First let us note that as the centraliser of a semisimple element contains a maximal torus, it follows from \cite[Lem. 11.4.11]{digne-michel:2020:representations-of-finite-groups-of-lie-type} that any $\bG^F$-conjugacy class of semisimple elements contained in $\rO^{p'}(\bG^F)$ is a single $\rO^{p'}(\bG^F)$-class. Now, we claim that at most one $\bG^F$-class in $\mathcal{O}^F$ is contained in $\rO^{p'}(\bG^F)$. Fix $\tilde{s} \in \pi^{-1}(s)$ then as $F(s) = s$ we have $\tilde{s}^{-1}F(\tilde{s}) \in \Ker(\pi) \leqslant \rZ(\bG_{\simc})$. Recall we have a homomorphism $\gamma_s : \rC_{\bG}(s) \to \Ker(\pi)$ given by $\gamma_s(x) = [\tilde{x},\tilde{s}]$ where $\tilde{x} \in \pi^{-1}(x)$. This factors through an injective homomorphism $\rC_{\bG}(s)/\rC_{\bG}^{\circ}(s) \to \Ker(\pi)$, see \cite[Cor.~2.8]{bonnafe:2005:quasi-isolated}. Now assume $s' = {}^gs \in \bG^F$ for some $g \in \bG$ and $\tilde{s}' \in \pi^{-1}(s')$ then $\tilde{s}'^{-1}F(\tilde{s}') = \tilde{s}^{-1}F(\tilde{s})\gamma_s(g^{-1}F(g))$. From this the claim follows easily.

Let $\tilde{\imath}_{\simc} : \mathbb{Q}\wc{\Phi} \to \bT_{\simc}$ be a surjective group homomorphism as in \ref{pa:torus-borel}. We use the notation for roots as in \cite{bourbaki:2002:lie-groups-chap-4-6}, which is consistent with our setup. We may take $\tilde{s} = \tilde{\imath}_{\simc}(\frac{1}{2}\wc{\varpi}_a) \in \pi^{-1}(s)$ as a preimage of $s$, where
%%%%
\begin{equation*}
\wc{\varpi}_a = \begin{cases}
\wc{\alpha}_{a+1}+2\wc{\alpha}_{a+2}+\cdots+(b-2)\wc{\alpha}_{n-2} + \frac{1}{2}\left((b-2)\wc{\alpha}_{n-1}+b\wc{\alpha}_n\right) &\dim(V)\text{ even,}\\
\wc{\alpha}_{a+1}+2\wc{\alpha}_{a+2}+\cdots+(b-2)\wc{\alpha}_{n-2} +(b-1)\wc{\alpha}_{n-1} + \frac{b}{2}\wc{\alpha}_n &\dim(V)\text{ odd.}
\end{cases}
\end{equation*}
%%%%
A straightforward calculation shows that if $F$ is split then $F(\tilde{s}) = \tilde{s}$ if and only if $\frac{(q-1)b}{4} \in \mathbb{Z}$. Similarly if $F$ is twisted then $F(\tilde{s}) = \tilde{s}$ if and only if $\frac{(q-1)b}{4} \in \frac{1}{2} + \mathbb{Z}$. One readily checks the equivalence of these conditions with those given in (i) and (ii).
\end{proof}

\begin{pa}\label{pa:K-inv-char}
Assume now that $\rZ(\bG^F) \leqslant \rO^{p'}(\bG^F)$. Let $s \in \bG^{\star F^{\star}}$ be a quasi-isolated semisimple element. If $s\in \rO^{p'}(\bG^{\star F^{\star}})$ then by (iii) of Proposition~\ref{prop:cent-chars} we have $\omega_s = 1_{\rZ(\bG^F)}$ is trivial. In particular, by (i) of Proposition~\ref{prop:cent-chars} we have for any $\chi \in \mathcal{E}(\bG^F,s)$ that $\omega_{\chi} = 1_{\rZ(\bG^F)}$ is trivial so $\rK(\bG^F)_{\chi} = \rK(\bG^F)$ by (ii) of Lemma~\ref{lem:autos-triv-on-derived}.

Now assume $s \not\in \rO^{p'}(\bG^{\star F^{\star}})$ and $1 \neq z \in \rZ(\bG^{\star})$ is the unique nontrivial element, which is necessarily $F^{\star}$-fixed. The central character $\omega_{\chi} = \omega_s$ is non-trivial by (iii) of Proposition~\ref{prop:cent-chars}. It follows from Corollary~\ref{cor:K-series} that if $sz$ is not $\bG^{\star F^{\star}}$-conjugate to $s$ then no character in $\mathcal{E}(\bG^F,s)$ is $\rK(\bG^F)$-invariant. It is straightforward to see that $s$ and $sz$ are not even $\bG^{\star}$-conjugate unless $\dim(V) = 4m$ and $\rC_{\bG^{\star}}^{\circ}(s)$ has type $\D_m\times\D_m$. In this case $\mathcal{E}(\bG^F,s)$ is $\rK(\bG^F)$-invariant because $s$ and $sz$ are $\bG^{\star F^{\star}}$-conjugate for any $z \in \rZ(\bG^{\star F^{\star}})$, which is easily checked.
\end{pa}

\begin{thm}\label{thm:KG-action}
Assume $\bG = \SO(V)$ and $\dim(V) = 2n$ is even. If $s \in \bG^{\star F^{\star}}$ is a quasi-isolated semisimple element and $\chi \in \mathcal{E}(\bG^F,s)$ is a cuspidal character then:
%%%%
\begin{enumerate}
	\item $s \in \rO^{p'}(\bG^{\star F^{\star}})$ except when $n=2m$ is even, $q \equiv 1 \pmod{4}$, and $\SO(V_{-1}^{\star}(s))^{F^{\star}}$ is of type ${}^2\D_m(q)$,
	\item $\rK(\bG^F)_{\chi} = \rK(\bG^F)$ if and only if one of the following holds:
	%%%%
	\begin{itemize}
		\item $q^n \equiv -\epsilon \pmod{4}$,
		\item $s \in \rO^{p'}(\bG^{\star F^{\star}})$,
		\item $\dim(V_1^{\star}(s)) = \dim(V_{-1}^{\star}(s))$.
	\end{itemize}
\end{enumerate}
\end{thm}

\begin{proof}
(i). By the classification of unipotent characters, a quasisimple group of type $\D_n$, resp., ${}^2\D_n$, has a cuspidal unipotent character if and only $n$ is an even, resp., odd, square, see \cite[Table II]{lusztig:1978:representations-of-finite-chevalley-groups}. This is now a straightforward check using Lemma~\ref{lem:omega-inv}.

(ii). If $q^n \equiv -\epsilon \pmod{4}$ then $\rZ(\bG^F) \cap \rO^{p'}(\bG^F) = \{1\}$ so $\rK(\bG^F)$ is trivial. Assume this is not the case then the argument in \ref{pa:K-inv-char} shows that $\rK(\bG^F)_{\chi} = \rK(\bG^F)$ when $s \in \rO^{p'}(\bG^{\star F^{\star}})$. Hence, we can assume $s \not\in \rO^{p'}(\bG^{\star F^{\star}})$. Again by the argument in \ref{pa:K-inv-char} we need only show that $\rK(\bG^F)_{\chi} = \rK(\bG^F)$ when $\dim(V_1^{\star}(s)) = \dim(V_{-1}^{\star}(s))$.

Assume $\tau_{\alpha} \in \rK(\bG^F)$ is the unique nontrivial element so that $\alpha : \bG^F \to \rZ(\bG^F)$ factors through an isomorphism $\bG^F/\rO^{p'}(\bG^F) \to \rZ(\bG^F)$. By Proposition~\ref{prop:cent-chars} we have $\omega_{\chi}\circ\alpha$ is nontrivial so is equal to $\hat{z}$ with $z \in \rZ(\bG^{\star F^{\star}})$ the unique nontrivial element. As remarked above $\mathcal{E}(\bG^F,s)$ is $\rK(\bG^F)$-stable. Now any GGGC is clearly $\rK(\bG^F)$-invariant, as it is non-zero only on $p$-elements. Hence, the statement follows as in the proof of Proposition~\ref{prop:CquasiisolatedGGGR}.
\end{proof}

\begin{rem}
This shows that when $q \equiv -1 \pmod{4}$ and $s$ is quasi-isolated we have all cuspidal characters contained in $\mathcal{E}(\bG^F,s)$ are $\rK(\bG^F)$-invariant.
\end{rem}

\begin{pa}
Now, let $\bG=\GO(V)$ with $\dim(V)$ even.  We next wish to describe the action of $\rK(\bG^F)$ on $\Irr(\bG^F)$.  Toward this end, we begin with an analogue of Theorem~\ref{prop:GGGCext-auts} for $\rK(\bG^F)$. Recall from the proof of Lemma~\ref{lem:asogeny-orthogonal-fin} that $K(\bG^{\circ F}) \leqslant \Aut(\bG^{\circ F})$ is the image of $K(\bG^F) \leqslant \Aut(\bG^F)$ under the natural restriction map $\Aut(\bG^F) \to \Aut(\bG^{\circ F})$.

Note that since $\rO^{p'}(\bG^F)=\rO^{p'}(\bG^{\circ F})$ and $\rZ(\bG^F) = \rZ(\bG^{\circ F})$, we may identify $\rK({\bG^\circ}^F)$ with a subgroup of $\rK(\bG^F)$.  We keep the notation of \ref{pa:notation-1} and \ref{pa:notation-2}.  
\end{pa}

\begin{prop}\label{prop:GGGCext-KG}
Let $\bG=\GO(V)$ with $\dim(V)$ even.  Then for any unipotent element $u \in \bG^{\circ F}$, we have $\widehat{\Gamma}_{u,2}$ is fixed by $\rK({\bG}^F)$.  
\end{prop}
\begin{proof}
Let $\tau\in\rK({\bG}^F)$ be nontrivial.  Since $\bU(\lambda,-1)^F$ is comprised of $p$-elements and the Sylow $2$-subgroup $S_{u,\lambda}$ of $\rC_{\bG}(\lambda)_{\zeta_{u,\lambda}}^F$ contains $\rZ(\bG^F)$, we see that the subgroup $\bU(\lambda,-1)^F\cdot S_{u,\lambda}$ is preserved by $\tau$. Further, $\wh{\zeta}_{u,\la}^{\tau}$ is an extension of $\zeta_{u,\la}$ satisfying $o(\wh{\zeta}_{u,\la}^{\tau})=o({\wh{\zeta}_{u,\la}})=o({\zeta}_{u,\la})$. Hence $\wh{\zeta}_{u,\la}^{\tau}={\wh{\zeta}_{u,\la}}$ by the uniqueness of this property. Then we have
%%%%
\begin{equation*}
\widehat{\Gamma}_{u,2}^{\tau} = \Ind_{\bU(\lambda,-1)^FS_{u,\lambda}}^{\bG^F}(\widehat{\zeta}_{u,\lambda}^{\tau})= \Ind_{\bU(\lambda,-1)^FS_{u,\lambda}}^{\bG^F}({\widehat{\zeta}_{u,\lambda}})={\widehat{\Gamma}_{u,2}},
\end{equation*}
%%%%
as claimed.
\end{proof}

\begin{thm}\label{thm:cuspidal-KG}
Let $\bG=\GO(V)$ with $\dim(V)$ even, let $s\in {{\bG^\circ}^\star}^{F^\star}$ be a quasi-isolated semisimple element, and let $\chi\in\mathcal{E}({\bG^\circ}^F,s)$ be cuspidal such that $\chi$ is $\rK({\bG^\circ}^F)$-invariant.  Then $\chi$ extends to a $\rK({\bG}^F)$-invariant cuspidal character in $\Irr(\bG^F)$.
\end{thm}

\begin{proof}
This follows from Proposition \ref{prop:GGGCext-KG}, using the same argument as in Theorem \ref{prop:CquasiisolatedGGGR}.
\end{proof}

\begin{pa}
We next consider the action of $\rK({\bG^\circ}^F)$ on pairs $(\bL,\la)\in \Cusp_s(\bG^\circ, F)$ for $s$ quasi-isolated.  Note that since $\rZ(\bG^{\circ})^F\leq \bL^F$, the group $\bL^F$ is stabilised by $\rK(\bG^{\circ F})$.  Furthermore, recall that using Corollary~\ref{cor:K-series}, $\mathcal{E}({\bG^\circ}^F, s)$ is fixed by $\rK({\bG^\circ}^F)$ if and only if $s\in \rO^{p'}(\bG^F)$ or $\dim(V_1(s))=\dim(V_{-1}(s))$. Hence, we wish to determine when the series $\mathcal{E}({\bG^\circ}^F, \bL, \la)\subseteq \mathcal{E}({\bG^\circ}^F, s)$ are also preserved in this situation. We keep the notation of Section~\ref{sec:quasiisol}. In particular, $\bM_0$, $\bM_1$, and $\lambda = \lambda_1\boxtimes\psi$ are as in Lemma~\ref{lem:cusp-pair}.
\end{pa}

\begin{thm}
Let $\bG=\GO(V)$ with $\dim(V)$ even, and let $s\in {{\bG^\circ}^\star}^{F^\star}$ be a quasi-isolated semisimple element such that $\mathcal{E}({\bG^\circ}^F, s)$ is preserved by $\rK({\bG^\circ}^F)$.  Let $(\bL,\la)\in \Cusp_s(\bG^\circ, F)$.  Then
\begin{enumerate}
\item If $\epsilon=1$ or $\bL\neq \bT$, the character $\la$ is invariant under $\rK({\bG^\circ}^F)$ if and only if $s\in \rO^{p'}(\bG^F)$.
\item If $\epsilon=-1$ and $\bL=\bT$, the character $\la$ is invariant under $\rK({\bG^\circ}^F)$ if and only if  at least one of the following holds \begin{itemize} \item $s\in \rO^{p'}(\bG^F)$ or \item $q\equiv 1\pmod 4$ and $\la$ is trivial on $\bT_n^F$. \end{itemize}

\item Every character of $\mathcal{E}({\bG^\circ}^F, \bL, \la)$ is fixed by $\rK({\bG^\circ}^F)$ if and only if $\la$ is fixed by $\rK({\bG^\circ}^F)$.
\end{enumerate}
\end{thm}
\begin{proof}
Keep the notation of Section~\ref{sec:galoisSO}, so  $\la=\la_1\boxtimes\psi$, with $\la_1$ as in the proof of Proposition~\ref{prop:isolprincipalSO} or Theorem~\ref{mainthmgalSO}(i), depending on whether $\bL=\bT$ or $\bL\neq \bT$.  Let $\sigma\in \rK({\bG^\circ}^F)$ be nontrivial.

We begin by showing the ``only if" direction of (i) and (ii). For this, assume $s\not\in \rO^{p'}({\bG^\star}^{F^\star})$ and that $\dim(V_1(s))=\dim(V_{-1}(s))$.   Write $2b=\dim(V_1(s))=\dim(V_{-1}(s))$ and let $s=s_1s_0$ with $s_i\in {\bM_i^\star}^{F^\star}$.  For $i\in\{0,1\}$, write $2a_i:=\dim(V_1(s_i))$ and $2b_i:=\dim(V_{-1}(s_i))$.  Note that by Lemma~\ref{lem:omega-inv},  we have $q^b\equiv -\epsilon\pmod 4$. 

We claim that for each $g_0\in {\bM_0^\circ}^F,$ there is $g_1\in\bS_1^F$ such that $g_1g_0$ is not in $\rO^{p'}(\bG^F)$.  Note that we would otherwise have $\bS_1^F\subseteq \rO^{p'}(\bG^F)$.  If $q\equiv -1\pmod 4$, Theorem~\ref{thm:KG-action} yields that $s_0\in \rO^{p'}({\bM_0^\star}^{F^\star})$, so our assumption $s\not\in \rO^{p'}({\bG^\star}^{F^\star})$ implies $s_1\not\in\rO^{p'}({\bG^\star}^{F^\star})$.   Viewing $s_1$ instead as an element of $\bG^F\cong {\bG^\star}^{F^\star}$, we therefore see $\bS_1^F\not\subseteq \rO^{p'}(\bG^F)$. Similarly, if $q\equiv 1\pmod 4$, then $\epsilon=-1$ and by Lemma~\ref{lem:omega-inv}, the element $\mathbf{d}_m(-1)$, for example, is an element of $\bS_1^F$ but not $\rO^{p'}(\bG^F)$. This proves the claim.

First, assume $\bT\neq \bL$ and recall that $\la_1$ may be chosen to be trivial on $\bT_i^F$ for $1\leq i\leq a_1$ and $\varepsilon_{\bT_i^F}$ for $a_1+1\leq i\leq m$.  
Write $z, z_0,$ and $z_1$ for the generators of $Z(\bG^F)$, $Z(\bM_0^F)$, and $Z(\bM_1^F)$, respectively. 
Let $g_0\in {\bM_0^\circ}^F$ such that $\psi(g_0)\neq 0$, and let $g_1\in\bS_1^F$ such that $g=g_1g_0\in\bL^F$ is not a member of $\rO^{p'}(\bG^F)$. 
 Then
\[\la^\sigma(g)= \la(gz)= \psi(g_0z_0)\la_1(g_1z_1)=\la(g)\omega_{s_0}(z_0)(-1)^{(q-1)b_1/2}.\]
Using Lemma~\ref{lem:omega-inv}, we see 
%%%%
\begin{equation*}
\omega_{s_0}(z_0)= \begin{cases}
(-1)^{b_0(q-1)/2}   &\text{if }\epsilon=1,\\
(-1)^{b_0(q-1)/2+1} &\text{if }\epsilon=-1.
\end{cases}
\end{equation*}
%%%%
(Indeed, by Lemma~\ref{lem:omega-inv}, we have $\omega_{s_0}$ is trivial on $z_0$ if and only if $q^{b_0}\equiv \epsilon\pmod 4$, and we know $\omega_{s_0}(z_0)=-1$ otherwise.  So, considering the options of $q\pmod 4$ and the parities of $b_0$, we arrive at this statement.)
Then since $q^b\not\equiv\epsilon\pmod 4$ and $\la(g)\neq 0$, we have $\lambda^\sigma(g)\neq \lambda(g)$.

In the case $\bT=\bL$, we take $\la=\la_1$, $m=n$, $a_0=b_0=0$, and $a_1=b_1=b$. If $\epsilon=1$, $\la$ is therefore of the form above.  If $\epsilon=-1$, the restriction of $\la$ to $\bT_n^F$ may be either trivial or $\varepsilon_{\bT_n^F}$.  Let $g\in \bT^F$ such that $g\not\in\rO^{p'}(\bG^F)$.  Since $\la$ is linear, we know $\la(g)\neq 0$, and we have
%%%%
\begin{equation*}
\la^\sigma(g)= \la(gz)= \begin{cases}
\la(g)(-1)^{(q-1)(b-1)/2} & \text{if $\epsilon=-1$ and $\la|_{\bT_n^F}$ is trivial,}\\
\la(g)(-1)^{(q-1)(b-1)/2+(q-\epsilon)/2} & \hbox{otherwise.}
\end{cases}
\end{equation*}
%%%%
If $\epsilon=1$ or $\la$ is nontrivial on $\bT_n^F$, we therefore have $\la^\sigma\neq\la$, since $q^b\neq \epsilon\pmod 4$.  This completes the ``only if" direction of (i).  If $\epsilon=-1$ and $\la$ is trivial on $\bT_n^F$, we have $\la^\sigma=\la$ if and only if $q\equiv 1\pmod 4$ or $b$ is odd, giving (ii).

We now show the ``if" direction of (iii).  That is, we claim that if $\la^\sigma=\la$, then $\chi^\sigma=\chi$ for each $\chi\in\mathcal{E}({\bG^\circ}^F, \bL, \la)$.  Recall that $\chi$ may be written $\chi={R_\bL^{\bG^\circ}(\la)_\eta}$ for some $\eta\in\Irr(W^\circ(\la))$, as in \ref{pa:HLnotn}.  Let $\Lambda$ be an extension map with respect to $\bL^F\lhd \rN_{{\bG^\circ}^F}(\bL)$, and let $\delta_{\la,\sigma}\in\Irr(W^\circ(\la))$ be such that ${\Lambda(\la)}^\sigma= \delta_{\la,\sigma}\Lambda({\la}^\sigma)=\delta_{\la,\sigma}\Lambda({\la})$.  Recall that we may write $\Lambda(\lambda) = \Res_{\rN_{{\bG^\circ}^F}(\bL)_{\lambda}}^{\rN_{{\bG}^F}(\bL)_{\lambda}}\left(\Lambda_1(\lambda_1) \boxtimes \Lambda_0(\psi)\right)$ as in Lemma~\ref{extGO}.  Now, calculations exactly as those above for $\la$ yield that $\Lambda_1(\lambda_1) \boxtimes \Lambda_0(\psi)$ is fixed by $\sigma$ as well, so $\delta_{\la,\sigma}=1$.  Then by \cite[Theorem 4.6]{malle-spaeth:2016:characters-of-odd-degree}, we have $\left({R_\bL^\bG(\la)_\eta}\right)^\sigma=R_\bL^\bG({\la}^\sigma)_{{\eta}^\sigma\cdot \delta_{\la,\sigma}^{-1}}=R_\bL^\bG({\la})_{{\eta}^\sigma}.$  But note that $\eta^\sigma=\eta$ for all $\eta\in\Irr(W^\circ(\la))$, since $Z({\bG^\circ}^F)\leq \bL^F$, completing the claim.  
 
Recall that by Proposition~\ref{prop:cent-chars},  every $\chi \in \mathcal{E}({\bG^\circ}^F,s)$ is fixed by $\rK({\bG^\circ}^F)$ if $s\in \rO^{p'}(\bG^F)$.  Hence to complete the proof, it now suffices to show that $(\bL, \la)$ and $(\bL, \la)^\sigma$ define distinct Harish-Chandra series when $\la\neq\la^\sigma$.  Note that $(\bL, \la)$ is moved by $\sigma$ if and only if $\la$ is.  We claim $(\bL, {\la}^\sigma)\neq (\bL, \la)^w$ for any $w\in {\bW^\circ}^F$. Otherwise, we have $w\in \rN_{{\bG^\circ}^F}(\bL)$ and $w\sigma^{-1}$ fixes $\la$.  But from the structure of $\rN_{{\bG^\circ}^F}(\bL)$ and $\la$, we see the action of $w$ on $\la_1$ does not affect the number of $i$ such that $\Res_{\bT_i^F}^{\bS_1^F}{\la_1}$ is of the form $\varepsilon_{\bT_i^F}$. Considering the calculations above, we see that $\la^{w\sigma^{-1}}=\la$ is therefore impossible if $\la^{\sigma}\neq \la$.  It follows that the series $\mathcal{E}({\bG^\circ}^F, \bL, \la)$ is preserved by $\sigma$ if and only if $\la$ is fixed by $\sigma$.
\end{proof}

\section{Reductions to the Quasi-Isolated Case}\label{sec:red-qi-case}
\begin{pa}\label{pa:sa-pair}
The following discussion, aimed at clarifying the equivariance statement of \cite[Thm.~9.5]{taylor:2018:action-of-automorphisms-symplectic}, applies to any connected reductive algebraic group with no restriction on $p$. We also take this opportunity to upgrade this discussion to incorporate the action of the Galois group $\gal$. However, we first need to define a dual action of the Galois group on semisimple elements of $\bG^{\star}$.

Recall that we fixed an isomorphism $\imath : \mathbb{Z}_{(p)}/\mathbb{Z} \to \mathbb{F}^{\times}$ in \ref{pa:torus-borel} and an injective homomorphism $\jmath : \mathbb{Q/Z} \to \overline{\mathbb{Q}}^{\times}$ in \ref{pa:j-inj-hom}. The composition $\kappa = \jmath\circ\imath^{-1}$ is an injective group homomorphism $\mathbb{F}^{\times} \to \overline{\mathbb{Q}}^{\times}$. The image $\kappa(\mathbb{F}^{\times}) \leqslant \overline{\mathbb{Q}}^\times$ is simply the group of all roots of unity whose order is coprime to $p$. This is preserved by the action of $\gal$ on $\overline{\mathbb{Q}}$, so identifying $\mathbb{F}^{\times}$ with its image under $\kappa$ yields an automorphic action of $\gal$ on $\mathbb{F}^{\times}$.
\end{pa}

\begin{pa}\label{pa:gal-acts-s}
Let $\wc{X}^{\star}$ be the cocharacter group of the torus $\bT^{\star}$. Recall that we have a group isomorphism $\mathbb{F}^{\times}\otimes_{\mathbb{Z}}\wc{X}^{\star} \to \bT^{\star}$, as in \ref{pa:torus-borel}, that on simple tensors is given by $k\otimes \gamma \mapsto \gamma(k)$. Under this isomorphism, the action of $\gal$ on $\mathbb{F}^{\times}$ defines an automorphic action $\gal$ on $\bT^{\star}$. We denote the resulting homomorphism $\gal \to \Aut(\bT^{\star})$ by $\sigma \mapsto \sigma^{\star}$ and say $\sigma^{\star}$ is \emph{dual} to $\sigma$. To be specific, assume $s \in \bT^{\star}$ is an element of order $d > 0$ and let $\zeta_d = \tilde{\jmath}(\frac{1}{d})$ be a primitive $d$th root of unity. If $\sigma \in \gal$, then $\sigma(\zeta_d) = \zeta_d^k$ for some integer $k \in \mathbb{Z}$ coprime to $d$ and we have $\sigma^{\star}(s) = s^k$.

We will freely use the notation of \cite{taylor:2018:action-of-automorphisms-symplectic}, in particular that of \cite[\S8]{taylor:2018:action-of-automorphisms-symplectic}, but we maintain our preference here for right over left actions. Moreover, the reference tori $(\bT_0,\bT_0^{\star})$ of \cite{taylor:2018:action-of-automorphisms-symplectic} will simply be denoted by $(\bT,\bT^{\star})$ to be in keeping with our notation here. Now fix a semisimple element $s \in \bT^{\star}$ with $T_{\bW^{\star}}(s,F^{\star}) \neq \emptyset$ and a coset $a \in \mathcal{A}_{\bW^{\star}}(s,F^{\star}) = \bW^{\star \circ}(s) \setminus T_{\bW^{\star}}(s,F^{\star})$. After possibly replacing $s$ with a $\bW^{\star}$-conjugate, we will assume that the Levi cover of $\rC_{\bG^{\star}}(s)$, defined as in \cite[Def.~8.11]{taylor:2018:action-of-automorphisms-symplectic}, is standard of the form $\bL_I^{\star}$ for some subset $I \subseteq \Delta$.

Let $X^{\star}$ be the character group of $\bT^{\star}$ and set $\mathbb{R}X^{\star} = \mathbb{R}\otimes_{\mathbb{Z}}X$. Given a vector $v \in \mathbb{R}X^{\star}$, we denote by $[v] = \{kv \mid 0 < k \in \mathbb{R}\}$ the positive half-line spanned by $v$. Now let $\Aso^*(\bG,F) \leqslant \Aso(\bG,F)$ be the subgroup stabilising the pair $(\bT,\bB)$. If $\sigma^{\star} \in \Aso(\bG^{\star},F^{\star})$ is dual to $\sigma \in \Aso^*(\bG,F)$, then $\sigma^{\star}$ defines an automorphism of $\bW^{\star}$ and a permutation of the set of half-lines $[\Delta^{\star}]$ of the simple roots. To have a uniform notation, we assume that if $\sigma \in \gal$, then the dual $\sigma^{\star}$ induces the identity on $\bW^{\star}$ and $[\Delta^{\star}]$. We use the same conventions for $\sigma$ itself.
\end{pa}

\begin{pa}\label{pa:zsig-star}
We set $\mathscr{A}(\bG,F) = \Aso(\bG,F) \times \gal$ and $\mathscr{A}^*(\bG,F) = \Aso^*(\bG,F) \times \gal$. To each pair $(s,a)$, with $s \in \bT^{\star}$ an element such that $T_{\bW^{\star}}(s,F^{\star}) \neq \emptyset$ and $a \in \mathcal{A}_{\bW^{\star}}(s,F^{\star})$ a coset, we define a corresponding rational Lusztig series $\mathcal{E}_0(\bG^F,s,a)$ as in \cite[8.5]{taylor:2018:action-of-automorphisms-symplectic}. By the Lang--Steinberg Theorem there exists an element $g \in \bG^{\star}$ such that $g^{-1}F^{\star}(g) \in T_{\bG^{\star}}(s,F^{\star})$ maps onto $a^{\star} \in \mathcal{A}_{\bW^{\star}}(s,F^{\star})$. Letting $t = {}^gs \in \bG^{\star F^{\star}}$ we have $\mathcal{E}(\bG^F,t) = \mathcal{E}_0(\bG^F,s,a)$, see \cite[Lem.~8.6]{taylor:2018:action-of-automorphisms-symplectic}.

Under this identification, it is easy to translate statements about $\mathcal{E}(\bG^F,t)$ into statements about $\mathcal{E}_0(\bG^F,s,a)$. For instance, let us denote by $\mathscr{A}^*(\bG,F)_{s,a} \leqslant \mathscr{A}^*(\bG,F)$ the subgroup stabilising the series $\mathcal{E}_0(\bG^F,s,a)$. If $\sigma \in \mathscr{A}^*(\bG,F)$, then $\mathcal{E}_0(\bG^F,s,a)^{\sigma} = \mathcal{E}_0(\bG^F,\sigma^{\star}(s),\sigma^{\star}(a))$. Note that $\sigma^{\star}(a) = a$ if $\sigma \in \gal$. If $\sigma \in \Aso(\bG,F)$ then this follows by translating \cite[Prop.~7.2]{taylor:2018:action-of-automorphisms-symplectic} and if $\sigma \in \gal$ then this follows by translating \cite[Lem.~3.4]{schaeffer-fry-taylor:2018:on-self-normalising-Sylow-2-subgroups}.

It follows from this that if $\sigma \in \mathscr{A}^*(\bG,F)_{s,a}$ then there exists an element $x^{\star} \in \bW^{\star}$ such that
%%%%
\begin{equation*}
(s,a) = ({}^{x^{\star}}\sigma^{\star}(s), x^{\star}\sigma^{\star}(a)F^{\star}(x^{\star-1}))
\end{equation*}
%%%%
and a unique element $y^{\star} \in \bW_I^{\star}$ such that $[y^{\star}x^{\star}\sigma^{\star}(I^{\star})] = [I^{\star}]$, where $[I^{\star}]$ denotes the set of positive half-lines as in \ref{pa:gal-acts-s}. We set $z_{\sigma}^{\star} = y^{\star}x^{\star}$. The choice of $x^{\star}$ above is unique up to multiplying by an element of $\bW^{\star}(s) \leqslant \bW_I^{\star}$. Therefore, the element $z_{\sigma}^{\star}$ is uniquely determined by $\sigma$.

The coset $\bW_I^{\star}a$ contains a unique element $w_1^{\star} \in T_{\bW^{\star}}(I^{\star},F^{\star})$, and this element is fixed in the sense that $z_{\sigma}^{\star}\sigma^{\star}(w_1^{\star})F^{\star}(z_{\sigma}^{\star-1}) = w_1^{\star}$. Dually, we have $\sigma^{-1}(z_{\sigma}(I)) = I$ and $z_{\sigma}^{-1}F(\sigma^{-1}(w_1)z_{\sigma}) = F(w_1)$. As in \cite[Prop.~9.2]{taylor:2018:action-of-automorphisms-symplectic}, there exists an element $n_{\sigma} \in \rN_{\bG}(\bT_0)$ representing $z_{\sigma}$ such that $n_{\sigma}^{-1}F(\sigma^{-1}(n_{w_1})n_{\sigma}) = F(n_{w_1})$. The choice for $n_{\sigma}$ is unique up to multiplication by an element of $\bT^{Fw_1}$.

As $\Aso(\bG,F) = \Inn(\bG^F)\Aso^*(\bG,F)$ the map $\Aso^*(\bG,F) \to \Out(\bG,F)$ is surjective. We set $\mathscr{O}(\bG,F) = \mathscr{A}(\bG,F)/\Inn(\bG^F) \cong \Out(\bG,F) \times \gal$. A straightforward calculation shows that the map $\sigma \mapsto \sigma n_{\sigma}$ gives a well-defined group homomorphism $\mathscr{O}(\bG,F)_{s,a} \to \mathscr{O}(\bL_I,Fw_1)$, where $\mathscr{O}(\bG,F)_{s,a} \leqslant \mathscr{O}(\bG,F)$ is the subgroup preserving $\mathcal{E}_0(\bG^F,s,a)$. The following makes \cite[Thm.~9.5]{taylor:2018:action-of-automorphisms-symplectic} more precise.
\end{pa}

\begin{thm}\label{thm:equivariance}
Let $a_1 = aw_1^{\star -1} \in \mathcal{A}_{\bW_I^{\star}}(s,w_1^{\star}F^{\star})$. Then we have a bijection $(-1)^{l(w_1)}R_{I,w_1}^{\bG} : \mathcal{E}_0(\bL_I^{Fw_1},s,a_1) \to \mathcal{E}_0(\bG^F,s,a)$ which is $\mathscr{O}(\bG,F)_{s,a}$-equivariant with respect to the homomorphism $\mathscr{O}(\bG,F)_{s,a} \to \mathscr{O}(\bL_I,Fw_1)$ defined above. In particular, $R_{I,w_1}^{\bG}(\psi)^{\sigma} = R_{I,w_1}^{\bG}(\psi^{\sigma n_{\sigma}})$ for any $\psi \in \mathcal{E}_0(\bL_I^{Fw_1},s,a_1)$ and $\sigma \in \mathscr{A}(\bG,F)$.
\end{thm}

\begin{proof}
In the case of Galois automorphisms, we use \cite[Cor.~8.1.6, Prop.~9.1.6]{digne-michel:2020:representations-of-finite-groups-of-lie-type} to conclude that $R_{I,w_1}^{\bG}$ is $\gal$-equivariant. We can now argue as in \cite[Thm.~9.5]{taylor:2018:action-of-automorphisms-symplectic}.
\end{proof}

\subsection*{Proofs of Theorems~\ref{thm:completecharfields} and \ref{thm:equiv-J-decomp-no-D}}
\begin{pa}
We now consider the above discussion in the context of $\bG = \SL(V \mid \mathcal{B})$. Up to $\bW^{\star}$-conjugacy, we can assume that $s = s_1s_0 \in \bT^{\star}$ where, for some $1 \leqslant m \leqslant n$, we have $s_1 = \mathbf{d}_1(\xi_1)\cdots \mathbf{d}_m(\xi_m)$ and $s_0 = \mathbf{d}_{m+1}(\xi_{m+1})\cdots \mathbf{d}_n(\xi_n)$ with $\xi_i^2 \neq 1$ for all $1 \leqslant i \leqslant m \leqslant n$ and $\xi_i^2 = 1$ for all $m < i \leqslant n$. We let $\bM_1^{\star}\times \bM_0^{\star}$ be a group determined by the integer $m$, as in \ref{pa:subspace-subgrp}. Then the Levi cover of $\rC_{\bG^{\star}}(s)$ is a standard Levi subgroup $\bL_I^{\star}$ with $I = J_1 \sqcup J_0 \subseteq \Delta$ with these subsets defined such that $\bL_{J_1}^{\star} = \rC_{\bM_1^{\star}}(s)\cdot\bT^{\star}$ and $\bL_{J_0}^{\star} = \bM_0^{\star\circ}\cdot\bT^{\star}$.

It is clear that we have a bijection $T_{\bW_1^{\star}}(s_1,F^{\star}) \times T_{\bW_0^{\star}}(s_0,F^{\star}) \to T_{\bW^{\star}}(s,F^{\star})$, and in this way we get a bijection $\mathcal{A}_{\bW_0^{\star}}(s_0,F^{\star}) \to \mathcal{A}_{\bW^{\star}}(s,F^{\star})$. We have a unique decomposition $w_1^{\star} = v_1^{\star}v_0^{\star}$ with $v_i^{\star} \in T_{\bW_i^{\star}}(s_i,F^{\star})$ and we get a corresponding coset $b_i \in \mathcal{A}_{\bW_i^{\star}}(s_i,F^{\star})$ containing $v_i^{\star}$. It is easy to see that $v_i^{\star}$ is the unique element of the coset $\bW_{J_i}^{\star}v_i^{\star}$ contained in $T_{\bW_i^{\star}}(J_i^{\star},F^{\star})$. As in the proof of \cite[Thm.~16.2]{taylor:2018:action-of-automorphisms-symplectic}, we have a tensor product decomposition
%%%%
\begin{equation}\label{eq:tensor-decomp}
\mathcal{E}_0(\bL_I^{Fw_1},s,a_1) = \mathcal{E}_0(\bL_{J_1}^{Fv_1},s_1,b_1) \boxtimes \mathcal{E}_0(\bL_{J_0}^{Fv_0},s_0,b_0).
\end{equation}
%%%%
\end{pa}

\begin{pa}\label{pa:im-cur-Out}
To use Theorem~\ref{thm:equivariance}, we need to describe the image of the homomorphism $\mathscr{O}(\bG,F)_{s,a} \to \mathscr{O}(\bL_I,Fw_1)$, which we do now. For any $\sigma \in \mathscr{A}^*(\bG,F)_{s,a}$ we have $\sigma^{\star}(s_0) = s_0$. Hence, it is clear that the element $z_{\sigma}^{\star}$ defined in \ref{pa:zsig-star} is contained in the subgroup $\bW_1^{\star}$. Thus the element $n_{\sigma}$ may be chosen inside the group $\rN_{\bM_1}(\bS_1)$, where $\bS_1 \leqslant \bM_1$ is the maximal torus defined in \ref{pa:subspace-subgrp}. If $\bG = \SO(V)$, then the Weyl group lifts in $\bG$ and we may, in fact, choose $n_{\sigma} = n_{z_{\sigma}}$. If $\bG = \Sp(V)$ then this is not necessarily the case.

If $t \in \bG^{\star F^{\star}}$ is a semisimple element, then there exists $g \in \bG^{\star}$ and $s \in \bT^{\star}$ such that $t = {}^gs$. If $g^{-1}F^{\star}(g) \in T_{\bG^{\star}}(s,F^{\star})$ maps onto $a^{\star} \in \mathcal{A}_{\bW^{\star}}(s,F^{\star})$, then $\mathcal{E}(\bG^F,t) = \mathcal{E}_0(\bG^F,s,a)$, see \cite[Lem.~8.6]{taylor:2018:action-of-automorphisms-symplectic}.

Using the identification between the sets $\mathcal{E}_0(\bG^F,s,a)$ and $\mathcal{E}(\bG^F,t)$ described in \cite[Lem.~8.6]{taylor:2018:action-of-automorphisms-symplectic}, it is easy to translate the statements of Theorems~\ref{thm:completecharfields} and \ref{thm:equiv-J-decomp-no-D} into corresponding statements in the language used above. It is these translated statements that we prove below. We leave the translation to the reader.
\end{pa}

\begin{proof}[of Theorem~\ref{thm:completecharfields}]
One calculates easily using \ref{pa:im-cur-Out} that the image of $\gal_{s,a}:=\gal\cap \mathscr{O}(\bG,F)_{s,a}$ inside $\mathscr{O}(\bL_I,Fw_1)$ is contained in the subgroup $\mathscr{O}(\bL_{J_1},Fv_1)_{s_1,b_1} \leqslant \mathscr{O}(\bL_I,Fw_1)$. Theorem~\ref{thm:equivariance} implies that if $\chi \in \mathcal{E}_0(\bG^F,s,a)$, then $\chi = \pm R_{I,w_1}(\psi)$ for a unique $\psi = \psi_1\boxtimes \psi_0 \in \mathcal{E}_0(\bL_I^{Fw_1},s,a_1)$. The equivariance of this bijection implies that we need only calculate the stabiliser of $\psi$ under the image of $\gal_{s,a}$. The group $\bL_{J_1}$ is a direct product of general linear groups so a unicity of multiplicities argument shows that $\psi_1$ is $\mathscr{O}(\bL_{J_1},Fv_1)_{s_1,b_1}$-invariant, see \cite[Prop.~3.10]{schaeffer-fry-taylor:2018:on-self-normalising-Sylow-2-subgroups} for instance. It follows that $(\gal_{s,a})_{\chi} = (\gal_{s,a})_{\psi_0}$ so $\mathbb{Q}(\chi) = \mathbb{Q}_s(\psi_0)$. The statement now follows from Theorems~\ref{mainthmgalSO} and \ref{mainthmgalSp}, together with the fact that the fields $\mathbb{Q}_s$ and $\mathbb{Q}(\psi_0) \subseteq \mathbb{Q}(\tilde{\jmath}(\frac{1}{p}))$ are linearly disjoint since $p$ and $d$ are coprime.
\end{proof}

\begin{rem}\label{rem:gal-stab-L-series}
Let $\Cl_{\mathrm{ss}}(\bG^{\star F^{\star}})$ be the set of semisimple conjugacy classes of $\bG^{\star F^{\star}}$. We define an action of $\gal$ on $\Cl_{\mathrm{ss}}(\bG^{\star F^{\star}})$ by setting $\sigma^{\star}([{}^gs]) = [{}^g\sigma^{\star}(s)]$, where $s \in \bT^{\star}$ and $g \in \bG^{\star}$ are such that ${}^gs \in \bG^{\star F^{\star}}$. One easily checks this is well-defined, given the action $\sigma^\star(s)$ via power maps as in \ref{pa:gal-acts-s}. The stabiliser $\mathcal{G}_s \leqslant \mathcal{G}$ of $\mathcal{E}(\bG^F,s)$ is then the stabiliser of the class $[s] \in \Cl_{\mathrm{ss}}(\bG^{\star F^{\star}})$ under this action.
\end{rem}

\begin{proof}[of Theorem~\ref{thm:equiv-J-decomp-no-D}]
(i). If $s_0 = 1$, equivalently $s = s_1$, then the series $\mathcal{E}_0(\bG^F,s,a)$ is not $\gamma$-invariant because the $\bW^{\star}$-orbit of $s$ is not $\gamma^{\star}$-stable, so we can assume that $s_0 \neq 1$. In this case, we have $\gamma^{\star}(s) = s$. The set $\mathcal{A}_{\bW^{\star}}(s,F^{\star})$ has size $2$ so clearly $\gamma^{\star}(a) = a$ as $\gamma^{\star}$ fixes the coset $\bW^{\star\circ}(s) \in \mathcal{A}_{\bW^{\star}}(s,F^{\star})$.

The following argument can be extrapolated from \cite[\S6.21]{lusztig:1984:characters-of-reductive-groups}, see also \cite{digne-michel:1990:lusztigs-parametrization}. We defer to \cite{lusztig:1984:characters-of-reductive-groups,digne-michel:1990:lusztigs-parametrization} for the details. The coset $a_1 = aw_1^{\star -1}$ contains a unique element of minimal length, which we denote by $w_{a_1}^{\star} \in a_1$. The element $w_a^{\star} = w_{a_1}^{\star}w_1^{\star} \in a$ is then the corresponding minimal length element in $a$.

We pick a Jordan decomposition $J_{s,a_1}^{\bL_I, Fw_1} : \mathcal{E}_0(\bL_I^{Fw_1},s,a) \to \mathcal{E}_0(\rC_{\bL_I^{\star}}(s)^{w_{a_1}^{\star}F^{\star}},1)$ so that
%%%%
\begin{equation*}
J_{s,a_1}^{\bL_I, Fw_1}(R_{x^{\star}w_{a_1}^{\star}}^{\bL_I,Fw_1}(s)) = (-1)^{l(w_{a_1}^{\star})}R_{x^{\star}}^{C_{\bG^{\star}}(s),w_{a_1}^{\star}w_1^{\star}F^{\star}}(1)
\end{equation*}
%%%%
for all $x^{\star} \in \bW_I^{\star\circ}(s)$. We adorn the superscripts here with the relevant Frobenius endomorphisms for the given group. For the definition of the Deligne--Lusztig character in the disconnected group $C_{\bG^{\star}}(s)$ we refer to \cite[Prop.~1.3]{digne-michel:1990:lusztigs-parametrization}. There is then a unique bijection $J_{s,a}^{\bG,F} : \mathcal{E}_0(\bG^F,s,a) \to \mathcal{E}_0(\rC_{\bG^{\star}}(s)^{w_a^{\star}F^{\star}},1)$ satisfying the property
%%%%
\begin{equation*}
J_{s,a_1}^{\bL_I,Fw_1} = (-1)^{l(w_1^{\star})}(J_{s,a}^{\bG,F} \circ R_{I,w}^{\bG,F}).
\end{equation*}
%%%%
We claim this a Jordan decomposition. Firstly, by transitivity of induction, we have
%%%%
\begin{equation*}
R_{I,w_1}^{\bG,F}(R_{x^{\star}w_{a_1}^{\star}}^{\bL_I,Fw_1}(s)) = R_{x^{\star}w_{a_1}^{\star}w_1^{\star}}^{\bG,F}(s) = R_{x^{\star}w_a^{\star}}^{\bG,F}(s)
\end{equation*}
%%%%
for any $x^{\star} \in \bW_I^{\star \circ}(s) = \bW^{\star \circ}(s)$. This implies that
%%%%
\begin{equation*}
J_{s,a}^{\bG,F}(R_{x^{\star}w_a^{\star}}^{\bG}(s)) = (-1)^{l(w_1^{\star})}J_{s,a_1}^{\bL_I,Fw_1}(R_{x^{\star}w_{a_1}^{\star}}^{\bL_I,Fw_1}(s)) = (-1)^{l(w_a^{\star})}R_{x^{\star}}^{C_{\bG^{\star}}(s),w_a^{\star}F^{\star}}(1).
\end{equation*}
%%%%
That the signs agree is shown in \cite{lusztig:1984:characters-of-reductive-groups}. This shows the bijection is a Jordan decomposition.

After \ref{pa:im-cur-Out}, we see that the image of $\gamma$ in $\mathscr{O}(\bL_I,Fw_1)$ is the restriction of $\gamma$ to $\bL_I$. Using the tensor product decomposition, it suffices to prove the statement in the case $s = s_0$. But now the statement follows from Proposition~\ref{prop:unip-qi} and (iii) of Theorem~\ref{mainthmgalSO}.

(ii). Note that $\Aso^*(\bG,F)_{s,a} = \rD(\bG,\bT,F)\cdot\Gamma(\bG,F)_{s,a}$, where $\Gamma(\bG,F)_{s,a} \leqslant \Gamma(\bG,F)$ is the stabiliser of the series $\mathcal{E}_0(\bG^F,s,a)$ and $\rD(\bG,\bT,F) = \{\iota_g|_{\bG^F} \mid g \in \bT$ and $g^{-1}F(g) \in \rC_{\bG}(\bG^F) = \rZ(\bG)^F\}$. After \ref{pa:im-cur-Out}, we see that the image of $\Gamma(\bG,F)_{s,a}$ is contained in the subgroup $\Out(\bL_{J_1},Fv_1) \times \Gamma(\bL_{J_0},Fv_0) \leqslant \Out(\bL_I,Fw_1)$.

As $Z(\bL_J)$ is connected, $\Out(\bL_J,Fv_1)$ contains no diagonal automorphisms. Hence, the image of $\rD(\bG,\bT,F)$ is contained in the image of $\rD(\bL_{J_0},\bS_0,Fv_0)$ in the group $\Out(\bL_{J_0},Fv_0)$, where $\bS_0 \leqslant \bL_{J_0}$ is the maximal torus defined in \ref{pa:subspace-subgrp}. The statement now follows from the equivariance of the bijection in Theorem~\ref{thm:equivariance} and Theorem~\ref{mainthmgalSO}.
\end{proof}

\begin{pa}
We note that the proof of (ii) of Theorem~\ref{thm:equiv-J-decomp-no-D} given above also applies in the case $\bG = \Sp(V)$, clarifying somewhat the proof in \cite[Thm.~16.2]{taylor:2018:action-of-automorphisms-symplectic}. Now, we end with a discussion of the real-valued characters for these groups.
\end{pa}

\begin{cor}
Keep the hypothesis of Theorem~\ref{thm:completecharfields}.
%%%%
\begin{enumerate}
	\item If $\bG = \SO(V)$, then $\chi \in \mathcal{E}(\bG^F,s)$ is real-valued if and only if $s$ is ${\bG^{\star}}^{F^\star}$-conjugate to $s^{-1}$. 
	\item If $\bG = \Sp(V)$, then $\chi \in \mathcal{E}(\bG^F,s)$ is real-valued if and only if $s$ is ${\bG^{\star}}^{F^\star}$-conjugate to $s^{-1}$ and either $-1$ is not an eigenvalue of $s$ or $q\equiv 1\pmod 4$.
\end{enumerate}
\end{cor}

\begin{proof}
Let us continue the notation used in the proof of Theorem~\ref{thm:completecharfields} above. We have $\chi \in \mathcal{E}_0(\bG^F, s,a)$ is real if and only if $\Q(\psi_0)$ and $\Q_s$ are both real.  But note that $\Q_s$ is real if and only if the element of $\Gal(\Q(\zeta_d)/\Q)$ mapping $\zeta_d$ to $\zeta_d^{-1}$ is a member of $\gal_{s,a}$, and hence $(s,a)$ is in the same $\bW^{\star}$-orbit as $(s^{-1},a)$. Applying Remark~\ref{rem:gal-stab-L-series} completes the proof.
\end{proof}

\section*{Acknowledgements}
Part of this work was completed while the authors were in residence at the Mathematical
Sciences Research Institute in Berkeley, California during the Spring 2018 semester program
on Group Representation Theory and Applications, supported by the National Science Foundation 
under Grant No. DMS-1440140. The authors thank Radha Kessar for useful discussions concerning Section~\ref{sec:aut-uni-char} and Martin Liebeck for sharing a proof of Corollary~\ref{cor:sp-conj} whilst at the MSRI. The reduction to regular unipotent elements used in Section~\ref{sec:power-maps-2} was inspired by his proof.

Part of this work was also completed while the authors were in residence at the Isaac Newton Institute for Mathematical Sciences during the Spring 2020 programme Groups, Representations, and Applications: New Perspectives, supported by EPSRC grant number EP/R014604/1. 

Both authors thank the MSRI, INI, and the organisers of the programs for making their stays possible and providing a collaborative and productive work environment.

The authors also thank Gunter Malle for useful discussions and for his comments on an early preprint and Ryan Vinroot for pointing out an inaccuracy in an earlier version of this work. Finally, the authors thank the anonymous referees for their careful reading and invaluable comments.

\setstretch{0.96}
\renewcommand*{\bibfont}{\small}
\printbibliography
\end{document}